\documentclass{article} 
\usepackage[dvips]{graphicx}
\usepackage{subfigure}
\usepackage{pstricks, amssymb, xr, multind, multicol} 
\input rlepsf.tex
\newtheorem{proposition}[equation]{Proposition} 
 
\newtheorem{theorem}[equation]{Theorem} 
\newtheorem{exa}[equation]{Example} 
\newtheorem{ex}[equation]{Exercise} 
\newtheorem{s-ex}[equation]{Side-exercise} 
\newtheorem{exas}[equation]{Examples} 
\newtheorem{lemma}[equation]{Lemma} 
\newtheorem{sublemma}[equation]{Sublemma} 
\newtheorem{remar}[equation]{Remark} 
\newtheorem{remars}[equation]{Remarks} 
\newtheorem{nota}[equation]{Notation} 
\newtheorem{sremar}[equation]{Side-remark} 
\newtheorem{definitio}[equation]{Definition}

\newenvironment{remark}{\begin{remar} \rm }{\end{remar}} 
\newenvironment{remarks}{\begin{remars} \rm }{\end{remars}}

\newenvironment{definition}{\begin{definitio} \rm }{\end{definitio}} 
\newcommand{\IT}{\mbox{\rm Int}}

\newcommand{\CA}{{\cal A}}

\newcommand{\CF}{{\cal F}} 
 
\newcommand{\CH}{{\cal H}} 
\newcommand{\CI}{{\cal I}} 
 
\newcommand{\CL}{{\cal L}}

\newcommand{\ZZ}{\mathbb{Z}} 
\newcommand{\RR}{\mathbb{R}} 
\newcommand{\QQ}{\mathbb{Q}} 
\newcommand{\CC}{\mathbb{C}} 
\newcommand{\NN}{\mathbb{N}} 
\newcommand{\bp}{\noindent {\sc Proof: }} 
\newcommand{\eop}{\nopagebreak 
                        \hspace*{\fill}{$\diamond$} 
                        \medskip} 
\newcommand{\tata}{ \begin{pspicture}[0.2](0,0)(.5,.4) 
\pscircle(0.25,0.2){.25} 
\psline{*-*}(0.05,.2)(.45,.2) \end{pspicture}} 
\newcommand{\tatata}{\begin{pspicture}[.2](-.2,-.1)(.8,.6)
\psline{*-*}(.1,0)(.5,0)
\psline{*-*}(.1,.4)(.5,.4)
\pscurve(.1,0)(0,.2)(.1,.4)
\pscurve(.1,0)(.2,.2)(.1,.4)
\pscurve(.5,0)(.4,.2)(.5,.4)
\pscurve(.5,0)(.6,.2)(.5,.4)
\end{pspicture}}
\newcommand{\tetra}{\begin{pspicture}[.2](-.2,-.1)(.8,.6)
\psline{*-}(0,0)(.6,0)(.3,.2)
\psline{*-*}(.6,0)(.3,.5)(.3,.2)
\psline{*-}(.3,.5)(0,0)(.3,.2)
\end{pspicture}}
\begin{document} 
\thispagestyle{empty}
\begin{center}
\vskip2cm

\huge{Surgery formulae for finite type invariants of rational homology $3$-spheres}

\vskip1cm

\large{Christine Lescop}\\
\large{Institut Fourier, UJF Grenoble, CNRS}

\vskip1cm
\begin{large}
Pr\'epublication de l'Institut Fourier n$^{o}$ 697 (2007)
http://www-fourier.ujf-grenoble.fr/prepublications.html
\end{large}
\end{center}

\vskip 2cm

\begin{abstract}

We first present four graphic surgery formulae for the degree $n$ part $Z_n$ of the
Kontsevich-Kuperberg-Thurston
universal finite type invariant of rational homology spheres.

Each of these four formulae determines an alternate sum
of the form 
$$\sum_{I \subset N} (-1)^{\sharp I}Z_n(M_I)$$
where $N$ is a finite set of disjoint operations to be performed on a rational homology sphere $M$, and 
$M_I$ denotes the manifold resulting from the operations in $I$.
The first formula treats the case when $N$ is a set of $2n$ Lagrangian-preserving surgeries (a {\em Lagrangian-preserving surgery\/} 
replaces
a rational homology handlebody by another such without changing the linking numbers of curves in its
exterior).
In the second formula, $N$ is a set of $n$ Dehn surgeries on the components of a boundary link.
The third formula deals with the case of $3n$ surgeries on the components of an algebraically split link.
The fourth formula is for $2n$ surgeries on the components of an algebraically split link in which all Milnor triple
linking numbers vanish.
In the case of homology spheres, these formulae can be seen as a refinement
of the Garoufalidis-Goussarov-Polyak comparison of different filtrations of the rational vector space freely generated by oriented homology spheres (up to orientation-preserving homeomorphisms).

The presented formulae are then applied to the study of the variation of $Z_n$ under a $p/q$-surgery on a knot $K$. This variation is a degree $n$ polynomial in $q/p$ when the class of $q/p$ in $\QQ/\ZZ$ is fixed, and the coefficients of these polynomials are knot invariants, for which various
topological properties or topological definitions are given.
\vskip1cm


\end{abstract}

\newpage
\tableofcontents

\newpage
\section{Introduction}
In this article, we shall focus on the real finite type 
invariants of homology $3$-spheres 
in the sense of Ohtsuki, Goussarov and Habiro, and 
on the topological properties of the surgery formulae that these invariants satisfy.

M.~Kontsevich proposed a topological construction for an invariant $Z$ of 
oriented rational homology $3$-spheres using configuration space integrals.  
G.~Kuperberg and D.~Thurston proved that $Z$ is a universal finite type 
invariant for homology $3$-spheres,
in the sense of Ohtsuki, Goussarov and Habiro  \cite{kt,lesconst}.
Like the LMO invariant, the Kontsevich-Kuperberg-Thurston invariant $Z=(Z_n)_{n \in \NN}$
takes its values in a space of Jacobi diagrams $\CA=\prod_{n \in \NN}\CA_n$, and 
any real degree $n$ invariant $\nu$ of homology $3$-spheres is obtained from the
Kontsevich-Kuperberg-Thurston invariant $(Z_i)_{i \in \NN}$ by a composition 
with a linear form that kills the $Z_i$, for $i>n$.

We shall first prove four formulae for $Z_n$, for $n\in \NN$. 
Each of these four formulae will determine an alternate sum
of the form 
$$\sum_{I \subset N} (-1)^{\sharp I}Z_n(M_I)$$
where $N$ is a set of disjoint operations to be performed on $M$, and 
$M_I$ denotes the manifold resulting of the operations in $I$.
Our first formula, Theorem~\ref{thmflag}, will be a mere alternative statement of the main theorem of \cite{sumgen} and 
will treat the case when $N$ is a set of $2n$ Lagrangian-preserving surgeries (a {\em Lagrangian-preserving surgery\/} 
replaces
a rational homology handlebody by another such without changing the linking numbers of curves in its
exterior).
In our second formula, Theorem~\ref{thmfboun}, $N$ will be a set of $n$ rational surgeries on the components of a boundary link.
Our third formula, Theorem~\ref{thmfas}, will deal with the case of $3n$ surgeries on the components of an algebraically split link.
Our fourth formula, Theorem~\ref{thmfasmu}, will be for $2n$ surgeries on the components of an algebraically split link in which all Milnor triple
linking numbers vanish.

In the case of integral homology spheres, these results can be seen as refinements of the Garoufalidis-Goussarov-Polyak comparison of the filtrations 
of the vector space generated by homology spheres, with respect to algebraically
split links, boundary links or Lagrangian-preserving surgeries \cite{ggp,al}.

As it was proved by Garoufalidis in \cite{garouf}, a degree $n$
finite type invariant $\nu$ of homology spheres satisfies a surgery
formula that describes its variation under $1/r$-surgery on a knot $K$
as 
$$\nu(M(K;1/r))-\nu(M)=\sum_{i=1}^n \nu^{(i)}(K \subset M) r^i$$
where $\nu^{(i)}$ is a finite type knot invariant in the Vassiliev sense  for knots in $S^3$ as defined in \cite{bn}.

Since all the real finite type 
invariants of homology $3$-spheres can be obtained from the universal
LMO invariant by composition with a linear form on the
space of Jacobi diagrams \cite{lmo,le}, and since the LMO invariant is defined using the
Kontsevich integral and surgery presentations of manifolds, the knot
invariants  $\nu^{(i)}$ can be explicitly given in terms of the Kontsevich
integral of surgery presentations of the knots. See also
the {\AA}rhus construction \cite{Aa1,Aa2,Aa3}.

We seek for a better topological understanding of the invariants $\nu^{(i)}$, and we shall relate some of them to the topology of Seifert surfaces of the knots.
For example, for any degree $n$ invariant $\nu$,
we give a formula in terms of the entries of the Seifert matrix
and the weight system of $\nu$ for the leading coefficient $\nu^{(n)}$
of the surgery polynomial. See Theorem~\ref{thmpol}. 
We shall also prove that $\nu^{(i)}$ is of degree less than $2n$ for any $i<n$.
This specifies a result of Garoufalidis and Habegger who proved that 
$(\nu(M(K;1))-\nu(M))$ is a degree $2n$ knot invariant with the same weight system as a degree $2n$ knot invariant induced by the Alexander polynomial, by using the LMO invariant \cite{gh}.

Some of the results proved in this article can be refined in the extensively studied case of the Casson-Walker invariant
$\lambda=W_1 \circ Z_1$, where $W_1(\tata)=2$, that satisfies the well-known formula, for a knot $K$ in a homology sphere $M$,
$$\lambda(M(K;p/q))-\lambda(M)=\frac{q}{p}\lambda^{\prime}(K) + \lambda(L(p,-q))$$
where $\lambda^{\prime}(K)$ is half the second derivative of the normalized Alexander polynomial of $K$ evaluated at $1$ and $L(p,-q)$ is the lens space obtained by $p/q$-surgery on the unknot.
We shall prove some graphical formulae for $\lambda^{\prime}(K)$ and for its variation under surgeries on disjoint algebraically unlinked knots in Propositions~\ref{propcasknot}, \ref{propcasknotvar}, \ref{propcasknottwo}.

Next, we shall concentrate on the case of the degree $2$ invariant $\lambda_2= W_2 \circ Z_2$, where $W_2(\tetra)=1$ and $W_2(\tata \;\tata)=0$.
The invariants $\lambda$ and $\lambda_2$ generate the vector space of real-valued invariants of degree lower than $3$ that are additive under connected sum. We shall prove that
$\lambda_2$ satisfies the surgery formula
$$\lambda_2(M(K;p/q))-\lambda_2(M)=\lambda_2^{\prime \prime}(K) \left(\frac{q}{p}\right)^2 + w_3(K) \frac{q}{p} + c(q/p) \lambda^{\prime}(K) +\lambda_2(L(p;-q))$$
for a knot $K$ in a homology sphere $M$,
where $c(q/p)$ only depends on $q/p$ modulo $\ZZ$,
$\lambda_2^{\prime \prime}$ is explicitly given in Theorem~\ref{thmpol} and 
$w_3$ is a knot invariant, for which we shall prove various properties. These properties include a crossing change formula, Proposition~\ref{propvarwthree}, and a formula for genus one knots $K$, Proposition~\ref{propgenusone}.
For knots in $S^3$, $w_3$ is the degree $3$ knot invariant 
that changes sign under mirror image, and that maps the chord diagram
with three diameters to $(-1)$. In his thesis \cite{auc}, Emmanuel Auclair independently
obtained a formula for $w_3(K)$ in terms of topological invariants of curves of an arbitrary Seifert surface of $K$, that is fortunately equivalent to Proposition~\ref{propgenusone} in the genus one case.

The article is organized as follows. The main results are stated precisely
without proofs from Section~\ref{secstatelag} to Section~\ref{secstatewthree}. The proofs occupy the following sections. Questions and expected generalizations of the proved results are mentioned at the end.

\section{The Kontsevich-Kuperberg-Thurston universal finite type invariant $Z$}
\setcounter{equation}{0}
\subsection{Jacobi diagrams}
\label{subjac}
Here, a {\em {Jacobi diagram}\/} $\Gamma$ is a trivalent graph $\Gamma$ without simple loop like $\begin{pspicture}[.2](0,0)(.6,.4)
\psline{-*}(0.05,.2)(.25,.2)
\pscurve{-}(.25,.2)(.4,.05)(.55,.2)(.4,.35)(.25,.2)
\end{pspicture}$. The set of vertices of such a $\Gamma$ will be denoted by $V(\Gamma)$,
its set of edges will be denoted by $E(\Gamma)$.
A {\em {half-edge}\/} $c$ of $\Gamma$ is an element of
$$H(\Gamma)=\{c=(v(c);e(c)) | v(c) \in V(\Gamma); e(c) \in E(\Gamma);v(c) \in e(c)\}.$$
An {\em automorphism\/} of $\Gamma$ 
is a permutation $b$ of $H(\Gamma)$
such that for any $c,c^{\prime} \in H(\Gamma)$,
$$v(c)=v(c^{\prime}) \Longrightarrow v(b(c))=v(b(c^{\prime}))\;\;\;\mbox{and} \;\;\;e(c)=e(c^{\prime}) \Longrightarrow e(b(c))=e(b(c^{\prime})).$$
The number of automorphisms of $\Gamma$ is
denoted by $\sharp \mbox{Aut}(\Gamma)$. 
For example, $ \sharp \mbox{Aut}(\tata)=12$.
{\em An orientation\/} of a vertex of such a diagram $\Gamma$ 
 is a cyclic order of the three
half-edges that meet at that vertex.
A Jacobi diagram $\Gamma$ is {\em oriented\/} if all its vertices are oriented (equipped with an orientation).
The {\em degree} of such a diagram is 
half the number of its vertices. 

Let $\CA_n=\CA_n(\emptyset)$ 
denote the real vector space generated by the degree $n$ oriented Jacobi diagrams, quotiented out by the following relations AS and IHX:

$$ {\rm AS :}  \begin{pspicture}[.2](0,-.2)(.8,1)
\psset{xunit=.7cm,yunit=.7cm}
\psarc[linewidth=.5pt](.5,.5){.2}{-70}{15}
\psarc[linewidth=.5pt](.5,.5){.2}{70}{110}
\psarc[linewidth=.5pt]{->}(.5,.5){.2}{165}{250}
\psline{*-}(.5,.5)(.5,0)
\psline{-}(.1,.9)(.5,.5)
\psline{-}(.9,.9)(.5,.5)
\end{pspicture}
+
\begin{pspicture}[.2](0,-.2)(.8,1)
\psset{xunit=.7cm,yunit=.7cm}
\pscurve{-}(.9,.9)(.3,.7)(.5,.5)
\pscurve[border=2pt]{-}(.1,.9)(.7,.7)(.5,.5)
\psline{*-}(.5,.5)(.5,0)
\end{pspicture}=0,\;\;\mbox{and IHX :} 
\begin{pspicture}[.2](0,-.2)(.8,1)
\psset{xunit=.7cm,yunit=.7cm}
\psline{-*}(.1,1)(.35,.2)
\psline{*-}(.5,.5)(.5,1)
\psline{-}(.75,0)(.5,.5)
\psline{-}(.25,0)(.5,.5)
\end{pspicture}
+
\begin{pspicture}[.2](0,-.2)(.8,1)
\psset{xunit=.7cm,yunit=.7cm}
\psline{*-}(.5,.6)(.5,1)
\psline{-}(.8,0)(.5,.6)
\psline{-}(.2,0)(.5,.6)
\pscurve[border=2pt]{-*}(.1,1)(.3,.3)(.7,.2)
\end{pspicture}
+
\begin{pspicture}[.2](0,-.2)(.8,1)
\psset{xunit=.7cm,yunit=.7cm}
\psline{*-}(.5,.35)(.5,1)
\psline{-}(.75,0)(.5,.35)
\psline{-}(.25,0)(.5,.35)
\pscurve[border=2pt]{-*}(.1,1)(.2,.75)(.7,.75)(.5,.85)
\end{pspicture}
=0. 
$$

Each of these relations relate diagrams which can be represented by planar immersions that are identical outside the part of them represented in the pictures. Here, the orientation of vertices is induced by the counterclockwise order of the half-edges. For example, 
AS identifies the sum of two diagrams which only differ by the orientation
at one vertex to zero.
$\CA_0(\emptyset)$ is equal to $\RR$ generated by the empty diagram.

\subsection{The Kontsevich-Kuperberg-Thurston universal finite type invariant $Z$}

Let $\Lambda$ be $\ZZ$, $\ZZ/2\ZZ$ or $\QQ$.
A {\em $\Lambda$-sphere\/} is a  compact oriented 3-manifold $M$ such that $H_{\ast}(M;\Lambda)=H_{\ast}(S^3;\Lambda)$. A \index{TT}{Zsphere@$\ZZ$-sphere} 
$\ZZ$-sphere is also called a {\em homology sphere} while a 
{\em rational homology sphere} is a 
$\QQ$-sphere.
Following Witten, Axelrod, Singer, Kontsevich, Bott and Cattaneo, Greg Kuperberg and Dylan Thurston constructed invariants $Z_n=(Z_{KKT})_n$ of oriented $\QQ$-spheres valued in $\CA_n(\emptyset)$ and they proved that their invariants have the following property \cite{kt}. See also \cite{lesconst}.

\begin{theorem}[Kuperberg-Thurston \cite{kt}]
\label{thktone}
An  invariant $\nu$ of\/ $\ZZ$-spheres valued in a real vector space $X$ is
of degree $\leq n$ if and only if there exist
linear maps $$\phi_k(\nu): \CA_k(\emptyset) \longrightarrow X,$$ for any $k \leq n$, such that $$\nu=\sum_{k=0}^n
\phi_k(\nu) \circ Z_k.$$
\end{theorem}

A {\em real {finite type invariant}\/} of $\ZZ$-spheres is a topological invariant of 
$\ZZ$-spheres valued in a real vector space $X$ which is of degree less than
some natural integer $n$.
Theorem~\ref{thktone} can be used as a definition of degree $\leq n$
real-valued invariants of $\ZZ$-spheres. 

A degree $\leq n$ invariant $\nu$ is of {\em degree\/} $n$ if $\phi_n(\nu) \neq 0$. In this case,
$\phi_n(\nu)$ is the {\em {weight system}\/} of $\nu$ and is denoted by $W_{\nu}$.

Let $p^c \colon \CA_k(\emptyset) \rightarrow \CA_k(\emptyset)$ be the canonical linear projection of $\CA_k(\emptyset)$ onto its subspace $\CA^c_k(\emptyset)$ generated by the connected
diagrams, such that $p^c$ maps the non-connected diagrams to $0$ 
and the restriction of $p^c$ to $\CA^c_k(\emptyset)$ is of course the identity.
Then $Z_n^c=p^c \circ Z_n$ is additive under connected sum. Furthermore 
any real-valued degree $n$ invariant belongs to the algebra generated by
the $(\phi_{k,i} \circ Z_k^c)_{k \leq n}$ for linear forms $\phi_{k,i}$ generating the dual of $\CA^c_k(\emptyset)$.

\begin{remark} The above definition coincides with the Ohtsuki definition of real finite type invariants \cite{oht}. The Ohtsuki degree (that is always a multiple of $3$) is three times the above degree. See \cite{oht,ggp,hab,al} and references therein for more discussions about the various concepts of finite-type invariants.
\end{remark}

\section{Lagrangian-preserving surgeries} 
\setcounter{equation}{0}
\label{secstatelag}

\noindent{\bf Conventions:}
Unless otherwise mentioned, manifolds are compact and oriented.
Boundaries are always oriented with the outward normal first convention. 
The normal bundle $N(V)$ of an oriented submanifold $V$ in an oriented manifold $M$ is oriented so that the tangent bundle $T_xM$ of the ambient manifold $M$ at some $x \in V$ is oriented as  $T_xM = N_xV \oplus T_xV$.
If $V$ and $W$ are two oriented transverse submanifolds of an oriented manifold $M$, their intersection is oriented so that the normal bundle of $T_x(V \cap W)$ is the sum 
$N_xV \oplus N_xW$.
If the two manifolds are of complementary dimensions, then the sign of an intersection point is $+1$ if the orientation of its normal bundle coincides with the orientation of the ambient space that is equivalent to say that $T_xM=T_xV \oplus T_xW$. Otherwise, the sign is $-1$. If $V$ and $W$ are compact and if $V$ and $W$ are of complementary dimensions in $M$, their algebraic intersection is the sum of the signs of the intersection points, it is denoted by
$\langle V, W \rangle_M$.

Recall that the {\em linking number\/} $lk(J,K)$ of two disjoint knots $J$ and $K$ in a rational homology sphere $M$
is the algebraic intersection of $J$ with a surface $\Sigma_K$ bounded by $K$ if $K$ is null-homologous, that $lk(J,.)$ is linear on $H_1(M \setminus J)$, and that $lk(J,K)=lk(K,J)$.

The {\em Milnor triple linking number\/} $\mu(K_1,K_2,K_3)$ of three null-homologous knots $K_1, K_2, K_3$ that do not link each other algebraically in a rational homology sphere $M$ may be defined,
as the algebraic intersection of three Seifert surfaces $\Sigma_2$, $\Sigma_1$, $\Sigma_3$ of these knots in the complement of the other ones.
$$\mu(K_1,K_2,K_3)=-\langle \Sigma_1,\Sigma_2,\Sigma_3\rangle= -\langle \Sigma_1 \cap \Sigma_2, \Sigma_3\rangle=-lk(K_3,\Sigma_1 \cap \Sigma_2).$$

We now describe part of the behaviour of the $Z_n$ under the Lagrangian surgeries
defined below.

A {\em genus $g$ $\QQ$-handlebody}\/  is an (oriented, compact) 3-manifold $A$
with 
the same homology with rational coefficients as the standard (solid) handlebody
$H_g$ below.  
$$H_g = 
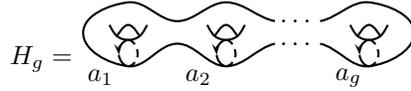
\begin{pspicture}[.4](0,-.5)(4.5,.95) 
\psset{xunit=.5cm,yunit=.5cm}   
\psecurve{-}(5.7,1.3)(5.2,1.3)(3.9,1.8)(2.6,1.3)(1.3,1.8)(.1,1)(1.3,.1)(2.6,.7)
(3.9,.1)(5.2,.7)(5.7,.7) 
\pscurve{-}(.8,1.2)(1,.9)(1.3,.8)(1.6,.9)(1.8,1.2) 
\pscurve{-}(1,.9)(1.3,1.2)(1.6,.9) 
\rput[r](.9,-.2){$a_1$} 
\psecurve{->}(1.6,.4)(1.3,.8)(1.05,.4)(1.3,.1) 
\psecurve{-}(1.3,.8)(1.05,.4)(1.3,.1)(1.6,.4) 
\psecurve[linestyle=dashed,dash=3pt 2pt](1,.4)(1.3,.8)(1.55,.4)(1.3,.1)(1,.4) 
\pscurve{-}(3.4,1.2)(3.6,.9)(3.9,.8)(4.2,.9)(4.4,1.2) 
\pscurve{-}(3.6,.9)(3.9,1.2)(4.2,.9) 
\rput[r](3.5,-.2){$a_2$} 
\psecurve{->}(4.2,.4)(3.9,.8)(3.65,.4)(3.9,.1) 
\psecurve{-}(3.9,.8)(3.65,.4)(3.9,.1)(4.2,.4) 
\psecurve[linestyle=dashed,dash=3pt
2pt](3.6,.4)(3.9,.8)(4.15,.4)(3.9,.1)(3.6,.4) 
\rput(5.8,1.3){\dots} 
\rput(5.8,.7){\dots} 
\psecurve{-}(5.8,1.3)(6.4,1.3)(7.7,1.8)(8.9,1)(7.7,.1)(6.4,.7)(5.8,.7)  
\pscurve{-}(7.2,1.2)(7.4,.9)(7.7,.8)(8,.9)(8.2,1.2) 
\pscurve{-}(7.4,.9)(7.7,1.2)(8,.9) 
\rput[l](6.8,-.2){$a_g$} 
\psecurve{->}(8,.4)(7.7,.8)(7.45,.4)(7.7,.1) 
\psecurve{-}(7.7,.8)(7.45,.4)(7.7,.1)(8,.4) 
\psecurve[linestyle=dashed,dash=3pt
2pt](7.4,.4)(7.7,.8)(7.95,.4)(7.7,.1)(7.4,.4) 
\end{pspicture}$$ 
Note that the boundary of such a $\QQ$-handlebody $A$ is homeomorphic to the 
boundary $(\partial H_g =\Sigma_g)$ of $H_g$. 

For a $\QQ$-handlebody $A$,  
${\cal L}_A$ denotes the kernel of the map induced by the inclusion: 
$$ H_1(\partial A;\QQ) \longrightarrow H_1( A;\QQ).$$ 
It is a Lagrangian of $(H_1(\partial A;\QQ),\langle,\rangle_{\partial A})$, we
call it the {\em {Lagrangian}}\/ of $A$.
A {\em Lagrangian-preserving surgery\/} or {\em LP--surgery} $(A,A^{\prime})$ is the replacement of a 
$\QQ$-handlebody $A$ embedded in a $3$-manifold by
another such $A^{\prime}$ with identical (identified via a homeomorphism) boundary and Lagrangian.

There is a canonical isomorphism 
$$\partial_{MV} \colon H_2(A \cup_{\partial A}
-A^{\prime};\QQ) \rightarrow {\cal L}_A$$
that maps the class of a closed surface in the closed $3$-manifold $(A \cup_{\partial A}
-A^{\prime})$ to the boundary of its intersection with $A$.
This isomorphism carries the algebraic triple intersection of surfaces to a trilinear antisymmetric form $\CI_{AA^{\prime}}$ on $\CL_A$.
$$\CI_{AA^{\prime}}(a_{i},a_{j},a_{k})=\langle \partial_{MV}^{-1}(a_i), \partial_{MV}^{-1}(a_j), \partial_{MV}^{-1}(a_k)\rangle_{A \cup
-A^{\prime}}$$

Let $(a_1,a_2,\dots,a_g)$ be a basis of $\CL_A$, and let $z_1,\dots,z_g$ be homology classes of $\partial A$, such that $(z_1,\dots,z_g)$ is dual to $(a_1,a_2,\dots,a_g)$ with respect to $\langle,\rangle_{\partial A}$ 
($\langle a_i,z_j \rangle_{\partial A}=\delta_{ij}$). Note that $(z_1,\dots,z_g)$
is a basis of $H_1(A;\QQ)$.
 
Represent $\CI_{AA^{\prime}}$
by the following combination $T(\CI_{AA^{\prime}})$ of tripods 
whose three univalent vertices form an ordered set:
$$T(\CI_{AA^{\prime}})=\sum_{\left\{\{i,j,k\}  \subset \{1,2,\dots ,g_A\} ; i <j <k\right\}}\CI_{AA^{\prime}}(a_{i},a_{j},a_{k})
\begin{pspicture}[0.2](-.2,-.1)(.5,.7) 
\psline{-}(0.05,.3)(.45,.6)
\psline{*-}(0.05,.3)(.45,.3) 
\psline{-}(0.05,.3)(.45,0)
\rput[l](.55,0){$z_i$}
\rput[l](.55,.3){$z_j$}
\rput[l](.55,.6){$z_k$}
\end{pspicture}$$
When $G$ is a graph with $2n$ trivalent vertices and with univalent vertices decorated by disjoint
curves of $M$,
define its contraction $\langle \langle G \rangle \rangle_n$ as the sum that runs over all the ways $p$ of gluing the 
univalent vertices two by two in order to produce a vertex-oriented Jacobi diagram $G_p$ $$\langle \langle G \rangle \rangle_n=\sum_p \ell(G_p)[G_p]$$ 
where $\ell(G_p)$ is the product over the pairs of glued univalent vertices, with respect to the {\em pairing\/} $p$, of the linking numbers of the corresponding curves. The contraction $\langle \langle . \rangle \rangle$ is linearly extended to linear combination
of graphs, and the disjoint union of combinations of graphs is bilinear.

A {\em $k$--component Lagrangian-preserving surgery datum\/} in a rational
homology sphere $M$ is a datum $(M;(A_i,A_i^{\prime})_{i \in \{1,\dots,k\}})$
of $k$ disjoint $\QQ$--handlebodies $A_i$, for $i \in \{1,\dots,k\}$, in $M$, and $k$ associated LP-surgeries $(A_i,A_i^{\prime})$.

\begin{theorem}
\label{thmflag}
Let  $$(M;(A_i,A_i^{\prime})_{i \in \{1,\dots,2n\}})$$
be a $2n$--component Lagrangian-preserving surgery datum in a rational homology sphere $\;M$. For $I \subset \{1,\dots,2n\}$, let $M_I$ denote the manifold obtained from $M$ by replacing $A_i$ by $A_i^{\prime}$ for all $i\in I$. Then
$$\sum_{I \subset \{1,\dots,2n\}} (-1)^{\sharp I}Z_n(M_I)=\langle \langle \bigsqcup_{i \in \{1,\dots,2n\}} T(\CI_{A_iA_i^{\prime}}) \rangle \rangle_n.$$
\end{theorem} 

We shall prove that this formula is equivalent to the formula of \cite{sumgen}
in Section~\ref{secprooflag}.

 

Let $\CF_0$ be the rational vector space freely generated by the oriented $\QQ$-spheres
viewed up to oriented homeomorphisms. 
For a $k$--component Lagrangian-preserving surgery datum $(M;(A_i,A_i^{\prime})_{i \in \{1,\dots,k\}})$
in a rational homology sphere $\;M$, define
$$[M;(A_i,A_i^{\prime})_{i \in \{1,\dots,k\}}]=\sum_{I \subset \{1,\dots,k\}} (-1)^{\sharp I}M_I \in \CF_0$$
and 
define $\CF_k$ as the subspace of $\CF_0$ generated by elements of $\CF_0$ of the above form. Then, it easily follows from the above theorem that $Z_n(\CF_{2n+1})=0$ where $Z_n$ is linearly extended to $\CF_0$.
For two elements $x$ and $y$ of $\CF_0$, we write $x \stackrel{n}{\equiv} y$
to say that $x-y \in \CF_{2n+1}$. Thus, if  $x \stackrel{n}{\equiv} y$, then
$Z_n(x)=Z_n(y)$.

The intersection of this filtration with the rational vector space freely generated by the oriented $\ZZ$-spheres is the Goussarov-Habiro filtration.
(The inclusion of the Goussarov-Habiro filtration $(\CF^{GH}_{k})_k$ in the intersection is obvious, the other one comes from the fact that $\CF^{GH}_{k}$ is the intersection of the kernels 
of the $Z_i$ for $2i < k$ because of the universality of $Z$.)

\section{Surgeries on algebraically split links}
\setcounter{equation}{0}

Let $L(p_i,-q_i)$ be the lens space obtained from $S^3$ by $p_i/q_i$-surgery on a trivial knot. (The standard conventions for surgery coefficients are recalled in the beginning of Section~\ref{seccasstate}.)
When $L=(K_i;p_i/q_i)_{i \in N}$ is a given link whose components are equipped with surgery coefficients in a rational homology sphere $\;M$, for $I\subset N$,
let
$$M_I=M_{(K_i;p_i/q_i)_{i \in I}}\sharp \sharp_{j \in N \setminus I} L(p_j,-q_j)$$
denote the connected sum of the manifold $M_{(K_i;p_i/q_i)_{i \in I}}=M\left((K_i;p_i/q_i)_{i \in I}\right)$ obtained from $M$ by surgery on $(K_i;p_i/q_i)_{i \in I}$ and the lens spaces $L(p_j,-q_j)$ for $j \notin I$.

Set 
$$[M;(K_i;p_i/q_i)_{i \in N}]=\sum_{I \subset N} (-1)^{\sharp I} M_I.$$

Note that the connected sums with lens spaces are trivial when the $p_i$ are $1$.

The invariant $Z_n$ is linearly extended to $\CF_0$.
By the additivity of the connected part $Z^c_n$ of $Z_n$ under
connected sum, if $N$ has more than one element,
$$Z^c_n\left([M;(K_i;p_i/q_i)_{i \in N}]\right)=Z^c_n\left(\sum_{I \subset N} (-1)^{\sharp I} M_{(K_i;p_i/q_i)_{i \in I}}\right)$$
and the connected sums with lens spaces do not appear in this case either.

In Section~\ref{secproofboun}, we shall see how Theorem~\ref{thmflag} easily implies the following surgery formula on $n$-component boundary links.
\begin{theorem}
\label{thmfboun} Let  $n$ and  $r$ be elements of  $\NN$.
Consider a link $(K_1,K_2, \dots, K_r)$ where all the $K_i$ bound disjoint
oriented surfaces $\Sigma^i$.
Let $p_i/q_i$ be a surgery coefficient for $K_i$, and let
$(x_j^i,y_j^i)_{j=1,\dots,g(\Sigma^i)}$ 
be a symplectic basis for the Seifert surface $\Sigma^i$.
Define $$I(\Sigma^i)= \sum_{(j,k) \in \{1,2,\dots, g(\Sigma^i)\}^2} 
\begin{pspicture}[0.2](-.5,-.1)(1.7,1) 
\psline{-}(0,0.05)(.15,.35)
\psline{*-*}(0,0.05)(1,.05) 
\psline{-}(0,0.05)(-.15,.35)
\psline{-}(1,0.05)(1.15,.35)
\psline{-}(1,0.05)(.85,.35)
\rput[br](-.15,.5){$x_j^i$}
\rput[b](.15,.5){$y_j^i$}
\rput[b](.85,.5){$(y_k^i)^+$}
\rput[lb](1.25,.5){$(x_k^i)^+$}
\end{pspicture}.$$
Then 
$$\begin{array}{lll}Z_n\left([M;(K_i;p_i/q_i)_{i \in \{1,\dots,r\}}]\right)&=0&\mbox{if}\; r>n
\\&=\frac{1}{2^n}\langle \langle \bigsqcup_{i \in \{1,\dots,n\}}  (-\frac{q_i}{p_i} I(\Sigma^i)) \;\;\rangle \rangle\;&\mbox{if}\; r=n.\end{array}$$
\end{theorem} 

\begin{definition}
A link $L$ in a $3$-manifold is said to be {\em  algebraically split} if 
any component of $L$ is null-homologous in the exterior of the other ones
(i.e. if any component of $L$ bounds a surface in the complement of the other components of $L$).
\end{definition}

An {\em edge-labelled\/} Jacobi diagram is a Jacobi diagram $\Gamma$ equipped with a bijection from $E(\Gamma)$ to $\{1,2,3,\dots, 3n\}$ for some integer $n$.
Let $D_{e,n}$ be the set of unoriented edge-labelled Jacobi diagrams of degree $n$. 
Let $L=(K_i)_{i \in \{1,2,3,\dots, 3n\}}$ be a $3n$--component algebraically split link. Let $\Gamma \in D_{e,n}$, orient $\Gamma$. 
To any vertex of $\Gamma$, whose incoming edges are labeled by $i,j,k$ with respect to the cyclic order induced by the orientation, associate the Milnor triple number $\mu(K_i,K_j,K_k)$. Then define $\mu_{\Gamma}(L)$ as the product over all the vertices of $\Gamma$ of the corresponding Milnor numbers of $L$. Note that $\mu_{\Gamma}(L)[\Gamma]$ does not depend on the orientation of $\Gamma$. Let $\theta(\Gamma)$ be the number of components of $\Gamma$ homeomorphic to $\tata$.

\begin{theorem}
\label{thmfas}
Let  $n$ and  $r$ be elements of  $\NN$.
Let $L=(K_i;p_i/q_i)_{i \in \{1,2,3,\dots, r\}}$ be a (rationally) framed $r$--component algebraically split link in a rational homology sphere $\;M$.
Then with the notation above,
$$\begin{array}{lll}Z_n\left([M;(K_i;p_i/q_i)_{i \in \{1,\dots,r\}}]\right)&=0&\mbox{if}\; r>3n
\\&= \prod_{i=1}^{3n}\frac{q_i}{p_i}\sum_{\Gamma \in D_{e,n}} \frac{\mu_{\Gamma}(L)}{2^{\theta(\Gamma)}}[\Gamma]&\mbox{if}\; r=3n.\end{array}$$ 
\end{theorem} 

A {\em $2/3$-labelled\/} Jacobi diagram is a degree $n$ Jacobi diagram $\Gamma$ equipped with an injection $\iota$ from $\{1,2,3,\dots, 2n\}$ to $E(\Gamma)$ such that at each vertex two edges of the image of $\iota$ meet one edge outside the image of $\iota$. 
Let $D_{2/3,n}$ be the set of unoriented $2/3$-labelled Jacobi diagrams of degree $n$. 
Let $(F_i)_{i=1, \dots 2n}$ be a collection of transverse oriented surfaces
that meet pairwise only inside their respective interiors,
such that $\langle F_i, F_j, F_k \rangle=0$ for any $\{i,j,k\} \subset \{1,2,3,\dots, 2n\}$.
Let $\Gamma \in D_{2/3,n}$, orient $\Gamma$. 
For a vertex of $\Gamma$, whose half-edges belong to edges labelled by $(i,j,$ nothing$)$, with respect to the cyclic order induced by the orientation,
assign the intersection curve $F_i \cap F_j$ to the unlabelled half-edge.
To any unlabelled edge $e$ that is now equipped with intersection curves $F_i \cap F_j$ and $F_k \cap F_{\ell}$ associate the linking number $\ell((F_i)_{i=1, \dots 2n};\Gamma;e)$
of $F_i \cap F_j$ and $F^+_k \cap F^+_{\ell}$, where $F^+_k$ and $F^+_{\ell}$ are
parallel copies of $F_k$ and $F_{\ell}$.

Note that there is no need to push the intersection curves by using parallels if $F_i$, $F_j$, $F_k$ and $F_{\ell}$ are distinct, to define this linking number.
If $\{i,j\}=\{k,\ell\}$, this linking number is the self-linking 
of the intersection curve that is framed by the surface, up to sign. 
Now, note that $lk(F_i \cap F_j, F^+_i \cap F^+_{\ell})= lk(F_i \cap F_j, F^+_i \cap F_{\ell})$ and that $$lk(F^+_i \cap F_j, F_i \cap F_{\ell})-lk(F_i \cap F_j, F^+_i \cap F_{\ell})=\pm \langle F_i, F_j, F_{\ell} \rangle.$$
Therefore if the cardinality of $\{i,j\} \cap \{k,\ell\}$ is $1$, the linking number $\ell((F_i)_{i=1, \dots 2n};\Gamma;e)$
is well-defined, too.
Define $\ell((F_i)_{i=1, \dots 2n};\Gamma)$ as the product over all the unlabelled edges of $\Gamma$ of the $\ell((F_i)_{i=1, \dots 2n};\Gamma;e)$.
Note that $\ell((F_i)_{i=1, \dots 2n};\Gamma)[\Gamma]$ is independent of the 
orientation of $\Gamma$.
Let $\sharp \mbox{Aut}_{2/3}(\Gamma)$ be the number of automorphisms of $\Gamma$ that preserve its $2/3$-labelling.

\begin{theorem}
\label{thmfasmu}
Let  $n$ and  $r$ be elements of  $\NN$.
Let $L=(K_i;p_i/q_i)_{i \in \{1,2,3,\dots, r\}}$ be a  framed $r$--component algebraically split link in a rational homology sphere $\;M$
such that for any $\{i,j,k\} \subset \{1,2,3,\dots, r\}$, $\mu(K_i,K_j,K_k)=0$.

Let $(F_i)_{i \in \{1,2,3,\dots, r\}}$ be a collection of transverse Seifert surfaces for the $K_i$ where $F_i$ does not meet
the $K_j$ for $i \neq j$. 

Then with the notation above
$$\begin{array}{lll}Z_n\left([M;(K_i;p_i/q_i)_{i \in \{1,\dots,r\}}]\right)&=0&\mbox{if}\; r>2n\\ &=\prod_{i=1}^{2n}\frac{q_i}{p_i}\sum_{\Gamma \in D_{2/3,n}} \frac{\ell((F_i)_{i=1, \dots 2n};\Gamma)}{\sharp \mbox{Aut}_{2/3}(\Gamma)}[\Gamma]&\mbox{if}\; r=2n\end{array}$$ 
where the sum runs over all $2/3$-labelled unoriented Jacobi diagrams $\Gamma$.
\end{theorem} 

When $M$ is a $\ZZ$-sphere, when the $p_i$ are equal to $1$, and when $r$ is greater or equal, than
$n$ for Theorem~\ref{thmfboun}, than $2n$ for Theorem~\ref{thmfasmu}, and than $3n$ for Theorem~\ref{thmfas},
the left-hand sides of the equalities of these theorems are in $\CF_{2n}^{GH}$. Since the degree $n$ part of the LMO invariant coincides with  $Z_n$ on $\CF_{2n}^{GH}$, these three theorems hold for the LMO invariant as well, in these cases.

Theorems~\ref{thmfas} and \ref{thmfasmu} will be proved in Section~\ref{secprooffas}.
Their proofs will rely on some clasper calculus performed in Section~\ref{secclasper}, that will also lead to the following proposition.
\begin{proposition}
\label{proprealmil}
Let $L=(K_i)_{i \in \{1,2,3,\dots, r\}}$ be an $r$--component algebraically split link in a rational homology sphere $\;M$. Then there exist transverse Seifert surfaces $\Sigma_i$ in $M \setminus \left(\cup_{j \neq i}K_j \right)$ for each component $K_i$ of $L$,
such that, for any triple $(K_i,K_j,K_k)$ of components of $L$, the geometric triple intersection of the surfaces $\Sigma_i$, $\Sigma_j$ and $\Sigma_k$ is made of $|\mu(K_i,K_j,K_k)|$ points.
\end{proposition}
Section~\ref{secprooffas} also contains an equivalent definition of the Matveev Borromean surgery 
(or surgery on a $Y$-graph), see Proposition~\ref{propborlag}.

\section{On the polynomial form of the knot surgery formula}
\setcounter{equation}{0}

Recall that for any rational homology sphere $M$, $Z_0(M)=1$.
Theorem~\ref{thmfboun} implies that for any knot $K$ that bounds a surface $F$ in a rational homology sphere $M$ and for any two coprime integers $p$ and $q$ such that $p\neq 0$, 
$$Z_1(M(K;\frac{p}{q}))-Z_1(M)=\frac12 \langle \langle  I(F)\rangle \rangle \frac{q}{p}  +Z_1(L(p,-q)).$$
We shall see in Section~\ref{secproofpol} that Theorem~\ref{thmfboun} also easily implies the following theorem. The first part of this theorem
is essentially \cite[Prop. 4.1]{garouf}.
\begin{theorem}
\label{thmpol}
Let $p$ and $q$ be coprime integers such that $p\neq 0$. Let $n \in \NN$, $n\geq 1$.
Let $K$ be a knot that bounds a Seifert surface $F$ in a rational homology sphere $M$.
Let $F^i$ be parallel copies of $F$ for $i \in \{1,\dots,n\}$, and let $L_i$ denote the framed
link made of the boundary components of $\cup_{j=1}^i F_j$, where each component is framed by $1$. 

Then
$$Z_n(M(K;\frac{p}{q+rp}))-Z_n(M)=\sum_{i=0}^n Y_{n,q/p}^{(i)}(K \subset M) (r+\frac{q}{p})^i$$
for any $r \in \ZZ$
where the coefficients $Y_{n,q/p}^{(i)}(K)$ satisfy the following properties.
\begin{itemize}
\item
$$Y_{n,q/p}^{(n)}(K)=\frac{(-1)^n}{n!}Z_n([M;L_n])= \frac{1}{n!2^n} \langle \langle \bigsqcup_{i \in \{1,\dots,n\}}  I(F^i) \;\;\rangle \rangle\; ,$$
\item if $n\geq 2$,
$$Y_{n,q/p}^{(n-1)}(K)=\frac{(-1)^{n-1}}{(n-1)!}\left(\left(\frac{n-1}2+\frac{q}{p}\right) Z_n\left([M;L_n]\right) +Z_n\left([M(K;\frac{p}{q});L_{n-1}]\right)\right),$$
\item if $n\geq 2$, $p^c(Y_{n,q/p}^{(n-1)})= Y_{n,q/p}^{(n-1)c}$ does not depend on $p$ and $q$, 
\item if $i \leq n-1$, $Y_{n,q/p}^{(i)}$ only depends on $q/p$ mod $\ZZ$,
\item if $U$ bounds a disk in $M$, then 
$Y_{n,q/p}^{(i)}(U \subset M)=0$ if $i>0$ and 
$$Y_{n,q/p}^{(0)}(U \subset M)=Z_n(M \sharp L(p,-q))-Z_n(M),$$
\item
$$Y_{n,0}^{(0)}(K \subset M)=0,$$
\item $$Y_{n,q/p}^{(i)}(K \subset M)=(-1)^{i+n}Y_{n,-q/p}^{(i)}(K \subset - M).$$
\end{itemize}
\end{theorem} 

A {\em singular knot} is an immersion of $S^1$ in a $3$-manifold whose only multiple points are transverse double points like
\begin{pspicture}[0.4](0,0)(.5,.4) 
\psline{->}(0.05,0.05)(0.45,0.45)
\psline{->}(0.05,0.45)(0.45,0.05)
\psdot(0.25,.25) \end{pspicture}.

Such a double point can be removed in a positive way \begin{pspicture}[0.4](0,0)(.5,.4) 
\psline{->}(0.05,0.05)(0.45,0.45)
\psline[border=1pt]{->}(0.05,0.45)(0.45,0.05)
\end{pspicture}
or in a negative way \begin{pspicture}[0.4](0,0)(.5,.4) 
\psline{->}(0.05,0.45)(0.45,0.05)
\psline[border=1pt]{->}(0.05,0.05)(0.45,0.45)
\end{pspicture}.

Let $K^s$ be a singular primitive knot in a rational homology sphere
with $k$ double points. Fix a bijection from $\{1,\dots,k\}$ to its set of double points. For $I \subset \{1,\dots,k\}$, let $K_I$ be the desingularisation
of $K^s$ such that the singular points in the image of $I$ have been removed in a negative way, and the
singular points outside the image of $I$ have become positive. If $y$ is a knot invariant valued in an abelian group, set
$$y(K^s)=\sum_{I \subset \{1,\dots,k\}}(-1)^{\sharp I}y(K_I).$$

\begin{remark}
It may happen that we do not know whether $Z_n(M(K_I;\frac{p}{q+pr}))$ is polynomial in $r$ for a given $I$, but that we know that 
$$\sum_{I \subset \{1,\dots,k\}}(-1)^{\sharp I}Z_n(M(K_I;\frac{p}{q+pr}))$$
is. Then the definition of $Y_{n,q/p}^{(i)}(K^s)$ extends in an obvious way.
\end{remark}

\begin{proposition}
\label{propzsing}
For any singular knot $K^s$ in a rational homology sphere
with $k$ double points, for any integers $n$, $i$, $q$ and $p$, with $0 \leq i \leq n$,\\
$Y_{n,q/p}^{(i)}(K^s)=0$ if $k > 2n$,\\
$Y_{n,q/p}^{(i)}(K^s)=0$ if $k > 2n-1$ and if $i < n$.\\
In other words, $Y_{n,q/p}^{(i)}$ is a knot invariant of degree at most $2n$
with respect to the crossing changes, and if $i<n$, $Y_{n,q/p}^{(i)}$ is a knot invariant of degree at most $(2n-1)$
with respect to the crossing changes.
\end{proposition}

Two disjoint pairs of points in $S^1$ are said to be {\em unlinked\/} if they
bound disjoint intervals in $S^1$. Otherwise, they are said to be {\em linked.\/}
Two double points of a singular knot are said to be {\em linked\/} if their
preimages are linked.

Associate a symmetric linking matrix $[\ell_{ij}(K^s)]_{i,j \in \{1,2,\dots,k\}}$ to a singular knot $K^s$ with $k$ pairwise unlinked double points numbered from $1$ to $k$ in the following way. 
Each double point $i$ \begin{pspicture}[0.2](0,0)(.5,.5) 
\psline{->}(0.05,0.05)(0.45,0.45)
\psline{->}(0.05,0.45)(0.45,0.05)
\psdot(0.25,.25) \end{pspicture} can be {\em smoothed\/} to transform the knot into two 
oriented singular knots $K^{s\prime}_i$ and $K^{s\prime\prime}_i$.
\begin{center}
\begin{pspicture}[0.2](0,0)(1,.7) 
\psecurve{->}(-0.05,-0.05)(0.05,0.05)(0.25,.15)(0.45,0.05)(0.55,-0.05)
\psecurve{->}(-0.05,.55)(0.05,0.45)(0.25,.35)(0.45,0.45)(0.55,.55)
\rput[l](.55,.45){$K^{s\prime}_i$}
\rput[l](.55,.05){$K^{s\prime\prime}_i$}
\end{pspicture}
\end{center}
Set $$\ell_{ii}(K^s)=lk(K^{s\prime}_i,K^{s\prime\prime}_i).$$
If $i$ and $j$ label two unlinked double points, let $K^{s,j}_i$ be the curve
among $K^{s\prime}_i$ and $K^{s\prime\prime}_i$ that does not contain the double point labeled by $j$, then
$\ell_{ij}(K^s)= lk(K^{s,j}_i,K^{s,i}_j)$ if $i \neq j$.

\begin{proposition}
\label{propzsingun}
For any singular knot $K^s$ in a rational homology sphere
with $k$ pairwise unlinked double points, for any integers $n$, $i$, $q$ and $p$, with $0 \leq i \leq n$, and $n\geq 1$,\\
if $k > n$, $Y_{n,q/p}^{(i)}(K^s)=0$,\\
if $k=n$, 
$Y_{n,q/p}^{(i)}(K^s)=Y_{n,0}^{(i)}(K^s)$ is an explicit homogeneous polynomial of degree $i$ in the coefficients of the linking matrix of $K^s$,
and $Y_{n,q/p}^{(0)}(K^s)=0$.
\end{proposition}

Proposition~\ref{propvarztwo} will give explicit examples of computations of the above homogeneous polynomials.
\section{A few formulae for the Casson-Walker invariant}
\setcounter{equation}{0}
\label{seccasstate}

Set $\lambda=W_1 \circ Z_1$ where $W_1(\tata)=2.$
According to \cite{sumgen}, $\lambda$ is the Casson-Walker invariant as normalized by Casson for $\ZZ$--spheres (see \cite{akmc,gm,mar}), $\lambda$ is half the Walker invariant as normalized in \cite{wal}, and $\lambda$ coincides with $\frac{\overline{\lambda}}{|H_1(M)|}$ where $|H_1(M)|$ denotes the cardinality of $H_1(M;\ZZ)$ for $\QQ$--spheres, and $\overline{\lambda}$ is the extension of $|H_1(M)|\lambda$ to oriented closed $3$-manifolds that is denoted by $\lambda$ in \cite{pup}.

A {\em rationally algebraically split link\/} is a link whose components
do not link each other.
The following proposition gives formulae that generalize Theorem~\ref{thmfboun},
Theorem~\ref{thmfas} and Theorem~\ref{thmfasmu} in the degree $1$ case ($n=1$) 
for rationally algebraically split links.

The {\em order\/} of a knot $K$ in a rational homology sphere is the smallest positive integer $O_K$ such that $O_K K$ is null-homologous.
A {\em primitive curve\/} on a torus $S^1 \times S^1$ is a non-separating simple closed curve on the torus.
A {\em primitive satellite\/} of a knot is a primitive curve on the boundary $\partial N(K)$ of its tubular neighborhood.
A surgery on a knot $K$ is determined by a primitive satellite $\mu$ (oriented arbitrarily) of the knot
that will bound a disk inside the surgered torus after surgery. If $m_K$ is the meridian of $K$, the isotopy class of such a curve
is determined by the pair $$(p_K =lk(\mu,K),q_K =\langle m_K, \mu \rangle_{\partial N(K)})$$
and the {\em surgery coefficient\/} is $p_K/q_K$.

For any order $d$ component $K$ of a rationally algebraically split link $L$, there exists an embedded surface $\Sigma$ in the complement of $L$ whose boundary $\partial \Sigma$ is made of essential parallel curves of the boundary $\partial N(K)$ of the tubular neighborhood $N(K)$ of $K$ such that $\partial \Sigma$ is homologous to $d$ parallels of $K$ in $N(K)$. Let $H_1(\Sigma)/H_1(\partial \Sigma)$ denote the quotient of $H_1(\Sigma)$ by the image of $H_1(\partial \Sigma)$ under the map induced by the inclusion. Let $B_s=(x_i,y_i)_{i \in \{1,\dots, g\}}$ be a symplectic basis of $H_1(\Sigma)/H_1(\partial \Sigma)$, define $$I(\Sigma)= \sum_{(j,k) \in \{1,2,\dots, g\}^2}
\begin{pspicture}[0.2](-.5,-.1)(1.7,1) 
\psline{-}(0,0.05)(.15,.35)
\psline{*-*}(0,0.05)(1,.05) 
\psline{-}(0,0.05)(-.15,.35)
\psline{-}(1,0.05)(1.15,.35)
\psline{-}(1,0.05)(.85,.35)
\rput[br](-.15,.5){$x_j$}
\rput[b](.15,.5){$y_j$}
\rput[b](.85,.5){$y_k^+$}
\rput[lb](1.25,.5){$x_k^+$}
\end{pspicture}.$$

If $W_n \colon \CA_n \rightarrow \QQ$ is a linear form, then 
$W_n\left( \langle \langle \cdot \rangle \rangle\right)$ will also be denoted by $\langle \langle \cdot \rangle \rangle_{W_n}$.

For example, 
$$\langle \langle I(\Sigma) \rangle \rangle_{W_1}= 2 \sum_{(j,k) \in \{1,2,\dots, g\}^2} \left(lk(x_j,x_k^+)lk(y_j,y_k^+)-lk(x_j,y_k^+)lk(y_j,x_k^+)\right).$$

\begin{proposition}
\label{propcasknot}
Let $n$ be an integer. Set $N=\{1,\dots,n\}$.
Let $L=(K_i;p_i/q_i)_{i \in N}$ be a framed rationally algebraically split link in a rational homology sphere $M$.
Let $d_i$ be the order of $K_i$ in $H_1(M)$, let $\Sigma_i$ be a surface of $M\setminus L$ whose boundary is made of essential parallel curves of $\partial N(K_i)$ and is homologous to $d_iK_i$ in $N(K_i)$.
If $n=1$, assume that the $\QQ/\ZZ$--self-linking number of $K_1$ is zero.

Then 
$$\sum_{I \subset N} (-1)^{\sharp I} \lambda \left(M_{(K_i;p_i/q_i)_{i \in I}}\sharp \sharp_{j \in N \setminus I} L(p_j,-q_j)\right)
=(-1)^n \left(\prod_{i=1}^n\frac{q_i}{p_i}\right) \lambda^{\prime}(L)$$
where
$$\begin{array}{lll}\lambda^{\prime}(L)
&= \frac{\langle \langle I(\Sigma_1) \rangle \rangle_{W_1}}{2d_1^2}+ \frac{1}{12} - \frac{1}{12d_1^2} \; &\mbox{if}\; n=1\\
&=\frac{\langle \langle I(\Sigma_1) \subset M(K_2;1)\rangle \rangle_{W_1}}{2d_1^2} - \frac{\langle \langle I(\Sigma_1) \subset M \rangle \rangle_{W_1}}{2d_1^2}= - \frac{lk(\Sigma_1 \cap \Sigma_2,(\Sigma_1 \cap \Sigma_2)_{\parallel})}{d_1^2d_2^2} \; &\mbox{if}\; n=2\\&= \frac{\langle \Sigma_1,\Sigma_2,\Sigma_3 \rangle^2}{d_1^2d_2^2d_3^2} \; &\mbox{if}\; n=3\\
&=0 \; &\mbox{if}\; n \geq 4\end{array}$$
and, if $n>1$,
$$\sum_{I \subset N} (-1)^{\sharp I} \lambda \left(M_{(K_i;p_i/q_i)_{i \in I}}\sharp \sharp_{j \in N \setminus I} L(p_j,-q_j)\right)= \sum_{I \subset N} (-1)^{\sharp I} \lambda \left(M_{(K_i;p_i/q_i)_{i \in I}}\right).$$
\end{proposition}

This proposition is proved in Section~\ref{secproofcasone}.
Under its hypotheses, we then obviously have the following equalities 
 $$\lambda^{\prime}(K_1 \subset M(K_2,p/q))-\lambda^{\prime}(K_1 \subset M)
=\frac{q}{p}\lambda^{\prime}(K_1,K_2)$$
 and
 $$\lambda^{\prime}\left(K_1  \subset M((K_2;p_2/q_2),(K_3;p_3/q_3))\right)-\lambda^{\prime}\left(K_1  \subset M(K_2;p_2/q_2)\right) -\lambda^{\prime}\left(K_1  \subset M(K_3;p_3/q_3)\right)$$
$$+\lambda^{\prime}(K_1 \subset M) =\frac{q_2q_3}{p_2p_3}\lambda^{\prime}(K_1,K_2,K_3).$$
 Then
the variation of linking numbers under surgery recalled in Lemma~\ref{lemvarlk} easily implies the following proposition (see also the proof of Lemma~\ref{lemvarlambdaprime}).

\begin{proposition}
\label{propcasknotvar}
Let $(K_1,K_2,K_3)$ be a rationally algebraically split link in a rational homology sphere $M$.
Let $d_i$ be the order of $K_i$ in $H_1(M)$, let $\Sigma_i$ be a surface of $M\setminus L$ whose boundary is made of essential parallel curves of $\partial N(K_i)$ and is homologous to $d_iK_i$ in $N(K_i)$. Then

$$\begin{array}{ll}\lambda^{\prime}(K_1,K_2)& = -\frac{1}{4} \langle \langle \frac{1}{d^2_1} I(\Sigma_1) \; \begin{pspicture}[0.2](-.7,-.2)(1.1,.6)  
\psline{-}(.35,0.05)(-.05,.05)
\rput[rt](-.1,.25){$K_2$}
\rput[lt](.4,.25){$K_{2\parallel}$}
\end{pspicture}  \rangle \rangle_{W_1}\\
& = -\frac{1}{4} \langle \langle \frac{1}{d^2_2} I(\Sigma_2) \; \begin{pspicture}[0.2](-.7,-.2)(1.1,.6)  
\psline{-}(.35,0.05)(-.05,.05)
\rput[rt](-.1,.25){$K_1$}
\rput[lt](.4,.25){$K_{1\parallel}$}
\end{pspicture} \rangle \rangle_{W_1},
\end{array}$$

$$\lambda^{\prime}(K_1,K_2,K_3) = \frac{1}{8d_1^2} \langle \langle  I(\Sigma_1) \; \begin{pspicture}[0.2](-.7,-.2)(2.8,.6)  
\psline{-}(.35,0.05)(-.05,.05)
\rput[rt](-.1,.25){$K_2$}
\rput[lt](.4,.25){$K_{2\parallel}$}
\psline{-}(2.05,0.05)(1.65,.05)
\rput[rt](1.6,.25){$K_3$}
\rput[lt](2.1,.25){$K_{3\parallel}$}
\end{pspicture}  \rangle \rangle_{W_1}.$$
\end{proposition}
\eop

\begin{proposition}
\label{propcasknottwo}
If $K^s$ is a singular knot with one double point, then
$$\lambda^{\prime}(K^s)=\ell_{11}(K^s).$$ 
\end{proposition}
The easy proof of this well-known proposition is also given in Section~\ref{secproofcasone}.

\section{On the knot surgery formula for the degree $2$ invariant $\lambda_2$}
\setcounter{equation}{0}
\label{secstatewthree}
Consider the degree $2$ invariant
$$\lambda_2= W_2 \circ Z_2^c$$
where $W_2\left(\begin{pspicture}[.2](-.2,-.1)(.8,.6)
\psline{*-}(0,0)(.6,0)(.3,.2)
\psline{*-*}(.6,0)(.3,.5)(.3,.2)
\psline{*-}(.3,.5)(0,0)(.3,.2)
\end{pspicture}\right)=1$
and therefore $W_2\left(\tatata\right)=2$.
The invariant $\lambda_2$ is invariant under orientation change and additive under connected sum.

\begin{theorem}
\label{thmw3}
There exists a function $c$ from $\QQ/\ZZ$ to $\QQ$ such that $c(0)=0$, $c(q/p)=c(-q/p)$ and the following assertions hold.
Let $r=q/p\in \QQ \setminus \{0\}$, where $p$ and $q$ are coprime integers.
Let $K$ be a knot that bounds a Seifert surface $F$ in a rational homology sphere $M$.
Let $F^1$ and $F^2$ be two parallel copies of $F$.
Then
$$\lambda_2(M(K;1/r))-\lambda_2(M)=\lambda_2^{\prime \prime}(K) r^2 + w_3(K)r +  C(K;q/p)+ \lambda_2(L(p;-q))$$
where $$\lambda_2^{\prime \prime}(K)= \frac{1}{8} \langle \langle \bigsqcup_{i \in \{1,2\}}  (I(F^i)) \;\;\rangle \rangle_{W_2}$$
$$w_3(K \subset M)=-w_3(K \subset (-M))$$
and $C(.;q/p)$ is an invariant of null-homologous knots that only depends on $q/p$ mod $\ZZ$, such that:\\
$C(K;0)=0$, and, if $K$ bounds a surface whose $H_1$ vanishes in $H_1(M)$, then  $C(.;q/p)=c(q/p)\lambda^{\prime}(K)$.

Furthermore,
if $K^s$ is a singular knot with two unlinked double points, then
$$w_3(K^s)=-\frac{\ell_{12}(K^s)}{2}\;\;\mbox{and}\;\; C(K^s;q/p)=0.$$ 
\end{theorem} 

Like all the statements in this section, the above theorem will be proved
in Section~\ref{secprooflambdatwo}.

\begin{proposition}
\label{propvarwthree}
Let $K^s$ be a singular 
knot with one double point in a rational homology sphere.
Let $K^+$ and $K^-$ be its two desingularisations, and let $K^{\prime}$ and $K^{\prime \prime}$
be the two knots obtained from $K^s$ by smoothing the double point. Assume that
$K^{\prime}$ and $K^{\prime \prime}$ are null-homologous, then
 $$w_3(K^+)-w_3(K^-)=\frac{\lambda^{\prime}(K^{\prime}) + \lambda^{\prime}(K^{\prime \prime})}{2} -\frac{\lambda^{\prime}(K^+) + \lambda^{\prime}(K^-)+lk^2(K^{\prime},K^{\prime \prime})}{4}.$$ 
\end{proposition}

Let \begin{pspicture}[.4](-.1,-.1)(.6,.6) 
\psframe(0,0)(.5,.5)
\rput(.25,.25){$x$}
\end{pspicture} denote a two strand braid 
with $|x|$ vertical juxtapositions of the motive 
\begin{pspicture}[.2](.65,.2)(1.35,.8)
\psecurve(1,.1)(1.15,.25)(1,.4)(.85,.55)(1,.7)
\psecurve[border=1pt](1,.7)(1.15,.55)(1,.4)(.85,.25)(1,.1)
\end{pspicture} if $x>0$ and
$|x|$ vertical juxtapositions of the motive 
\begin{pspicture}[.2](.65,.2)(1.35,.8)
\psecurve(1,.1)(.85,.25)(1,.4)(1.15,.55)(1,.7)
\psecurve[border=1pt](1,.7)(.85,.55)(1,.4)(1.15,.25)(1,.1)
\end{pspicture} if $x<0$.

Let $x$, $y$ and $z$ be three odd numbers.
Let $K(x,y,z)$ be the following pretzel knot that bounds a genus
one Seifert surface $\Sigma$ whose thickening $H$ coincides with the thickening
of the twice punctured disk next to it.
$H$ is a genus two handlebody whose boundary is equipped with
curves $X$, $Y$ and $Z$ that bound disks in its exterior.

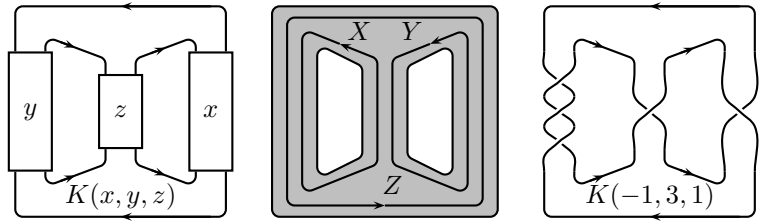
\begin{center}
\begin{pspicture}[0.4](-.2,-1.5)(3.2,1.5)
\psline[linearc=.1](2.9,.8)(2.9,1.4)(.1,1.4)(.1,.8)
\psline[linearc=.1](.5,.8)(.5,1)(1.3,.7)(1.3,.5)
\psline[linearc=.1](1.7,.5)(1.7,.7)(2.5,1)(2.5,.8)
\psline{->}(2.7,1.4)(1.5,1.4)
\psline{->}(.7,.925)(.9,.85)
\psline{->}(1.9,.775)(2.1,.85)
\psframe(1.2,-.5)(1.8,.5)
\psframe(0,-.8)(.6,.8)
\psframe(2.4,-.8)(3,.8)
\rput(.3,0){$y$}
\rput(1.5,0){$z$}
\rput(2.7,0){$x$}
\psline[linearc=.1](2.9,-.8)(2.9,-1.4)(.1,-1.4)(.1,-.8)
\psline[linearc=.1](.5,-.8)(.5,-1)(1.3,-.7)(1.3,-.5)
\psline[linearc=.1](1.7,-.5)(1.7,-.7)(2.5,-1)(2.5,-.8)
\psline{->}(2.7,-1.4)(1.5,-1.4)
\psline{->}(.7,-.925)(.9,-.85)
\psline{->}(1.9,-.775)(2.1,-.85)
\rput[b](1.5,-1.3){$K(x,y,z)$}
\end{pspicture}
\begin{pspicture}[0.4](-.2,-1.5)(3.2,1.5)
\pspolygon[linearc=.1,fillstyle=solid,fillcolor=lightgray](0,-1.4)(3,-1.4)(3,1.4)(0,1.4)
\pspolygon[linearc=.1,fillstyle=solid,fillcolor=white](.6,-.9)(1.2,-.6)(1.2,.6)(.6,.9)
\pspolygon[linearc=.1,fillstyle=solid,fillcolor=white](2.4,-.9)(1.8,-.6)(1.8,.6)(2.4,.9)
\rput[bl](1,.95){$X$}
\psline[linearc=.1]{->}(.9,.9)(.4,1.1)(.4,-1.1)(1.4,-.7)(1.4,.7)(.9,.9)
\rput[br](2,.95){$Y$}
\psline[linearc=.1]{->}(2.1,.9)(1.6,.7)(1.6,-.7)(2.6,-1.1)(2.6,1.1)(2.1,.9)
\rput[b](1.6,-1.1){$Z$}
\psline[linearc=.1]{->}(1.5,-1.25)(2.8,-1.25)(2.8,1.25)(.2,1.25)(.2,-1.25)(1.5,-1.25)
\end{pspicture}
\begin{pspicture}[0.4](-.2,-1.5)(3.2,1.5)
\psline[linearc=.1](2.9,.8)(2.9,1.4)(.1,1.4)(.1,.8)
\psline[linearc=.1](.5,.8)(.5,1)(1.3,.7)(1.3,.5)
\psline[linearc=.1](1.7,.5)(1.7,.7)(2.5,1)(2.5,.8)
\psline{->}(2.7,1.4)(1.5,1.4)
\psline{->}(.7,.925)(.9,.85)
\psline{->}(1.9,.775)(2.1,.85)
\psecurve(1.7,-.7)(1.7,-.5)(1.7,-.2)(1.5,0)(1.3,.2)(1.3,.5)(1.3,.7)
\psecurve[border=1pt](1.7,.7)(1.7,.5)(1.7,.2)(1.5,0)(1.3,-.2)(1.3,-.5)(1.3,-.7)
\psecurve(.3,0)(.5,.2)(.3,.4)(.1,.6)(.1,.8)(.1,1.4)
\psecurve[border=1pt](.5,1.4)(.5,.8)(.5,.6)(.3,.4)(.1,.2)(.3,0)
\psecurve(.3,-.4)(.5,-.2)(.3,0)(.1,.2)(.3,.4)
\psecurve[border=1pt](.3,.4)(.5,.2)(.3,0)(.1,-.2)(.3,-.4)
\psecurve(.5,-1.4)(.5,-.8)(.5,-.6)(.3,-.4)(.1,-.2)(.3,0)
\psecurve[border=1pt](.3,0)(.5,-.2)(.3,-.4)(.1,-.6)(.1,-.8)(.1,-1.4)
\psecurve(2.5,-1.4)(2.5,-.8)(2.5,-.2)(2.7,0)(2.9,.2)(2.9,.8)(2.9,1.4)
\psecurve[border=1pt](2.5,1.4)(2.5,.8)(2.5,.2)(2.7,0)(2.9,-.2)(2.9,-.8)(2.9,-1.4)
\psline[linearc=.1](2.9,-.8)(2.9,-1.4)(.1,-1.4)(.1,-.8)
\psline[linearc=.1](.5,-.8)(.5,-1)(1.3,-.7)(1.3,-.5)
\psline[linearc=.1](1.7,-.5)(1.7,-.7)(2.5,-1)(2.5,-.8)
\psline{->}(2.7,-1.4)(1.5,-1.4)
\psline{->}(.7,-.925)(.9,-.85)
\psline{->}(1.9,-.775)(2.1,-.85)
\rput[b](1.5,-1.3){$K(-1,3,1)$}
\end{pspicture}
\end{center}

Note that any genus one knot that bounds a genus one surface,
whose
$H_1$ goes to $0$ in $H_1(M)$, may be written as the image of $K(x,y,z)$ under
an embedding $\phi$ of $H$ into $M$ that maps $X$ and $Y$ to $0$ in  
$H_1(M \setminus \phi(H))$.

\begin{proposition}
\label{propgenusone}
Let $\phi$ be an embedding of $H$ in a rational homology sphere such that
$\phi(X)$ and $\phi(Y)$ are null homologous in the exterior of $\phi(H)$. Then
$$w_3(\phi(K(x,y,z)))=w_3(K(x,y,z)) - \frac{x}{2} \lambda^{\prime}(\phi(X)) - \frac{y}{2}  
\lambda^{\prime}(\phi(Y))
- \frac{z}{2} \lambda^{\prime}(\phi(Z))+\frac{3}{2}\lambda^{\prime}(\phi(X),\phi(Y)) $$ and
$$w_3(K(x,y,z))=\frac{x^2(y+z)+y^2(x+z)+z^2(x+y)}{32} +\frac{xyz}{8} +\frac{x +y +z}{16}$$
where $\lambda^{\prime}(\phi(X))$ and $\lambda^{\prime}(\phi(X),\phi(Y))$ are defined in several equivalent ways in Section~\ref{seccasstate}. 
\end{proposition}

\section{Proof of the Lagrangian-preserving surgery formula}
\setcounter{equation}{0}
\label{secprooflag}
In this section, we prove Theorem~\ref{thmflag} by proving that its formula is equivalent to the formula of \cite{sumgen} (or \cite{al} for the case of integral homology spheres). We first rewrite the right-hand side of the formula
of Theorem~\ref{thmflag}.

Let $g(i)$ be the genus of $A_i$. Let $(a^i_1,a^i_2,\dots,a^i_{g(i)})$ be a basis of $\CL_{A_i}$, and let $z^i_1,\dots,z^i_{g(i)}$ be homology classes of $\partial A_i$, such that $\langle a^i_j,z^i_k \rangle_{\partial A}=\delta_{jk}$.
Let $F$ be the set of maps $f$ from $\{1,\dots,2n\} \times \{1,2,3\}$ to $\NN$
such that $1 \leq f(i,1) < f(i,2) <f(i,3) \leq  g(i)$.
Let $P$ be the set of pairings $p$ 
of the disjoint union $G^0$ of the following $2n$ tripods, that pair a univalent vertex of some tripod to a univalent vertex of a different tripod.
 
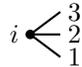
\begin{center}
\begin{pspicture}[0.2](-.2,-.1)(.8,.7) 
\psline{-}(0.05,.3)(.45,.6)
\psline{*-}(0.05,.3)(.45,.3) 
\psline{-}(0.05,.3)(.45,0)
\rput[l](.55,0){$1$}
\rput[l](.55,.3){$2$}
\rput[l](.55,.6){$3$}
\rput[r](-.1,.3){$i$}
\end{pspicture}
\end{center}
Let $p \in P$.
The half-edges of $G^0_p$ are naturally labeled in $\{1,\dots,2n\} \times \{1,2,3\}$. Assume that some $(f \in F)$ is given.
With a half-edge of $G^0_p$ labeled by $(i,j)$ that belongs to the tripod $i$, associate the curve $z_{f(i,j)}^i$ of $\partial A_i$. Then to an edge of 
$G^0_p$, associate the linking number of the curves associated to its two half-edges, and define $lk(p;f)$ as the product over the edges of these linking numbers. Set
$$c(p;f)=lk(p;f) \prod_{i=1}^{2n}\CI_{A_iA_i^{\prime}}(a_{f(i,1)}^i,a_{f(i,2)}^i,a_{f(i,3)}^i),$$
and
$$c(p)=\sum_{f \in F}c(p;f).$$

Then $$\langle \langle \bigsqcup_{i \in \{1,\dots,2n\}} T(\CI_{A_iA_i^{\prime}}) \rangle \rangle_n= \sum_{p\in P}c(p)[G^0_p]. $$
Let $D$ be the set of unoriented Jacobi diagrams of degree $n$. Consider a Jacobi diagram $\Gamma$ of $D$.
Let $P(\Gamma)$ be the set of the pairings $p$ of $P$ such that $G^0_p$ is isomorphic to $\Gamma$ as an unoriented Jacobi diagram.
Then $$\langle \langle \bigsqcup_{i \in \{1,\dots,2n\}} T(\CI_{A_iA_i^{\prime}}) \rangle \rangle_n= \sum_{\Gamma \in D} \sum_{p\in P(\Gamma)}c(p)[G^0_p]. $$

Fix $\Gamma$ in $D$.
Let $B(\Gamma)$ be the set of bijections $b$ from the set $H(\Gamma)$ of half-edges of $\Gamma$ to $\{1,\dots, 2n\} \times \{1,2,3\}$ 
that map any half-edge $c$ of a vertex $v(c)$ to three images with the same first coordinate $b_1(c)=b_1(v(c))$.
An element $b$ of $B(\Gamma)$ determines a pairing $p(b)$ of $P(\Gamma)$, and the number of elements of $B(\Gamma)$ that determine
the same pairing is the number of automorphisms of $\Gamma$.

$$\sum_{p\in P(\Gamma)}c(p)[G^0_p]= \sum_{b \in B(\Gamma)}\frac{c(p(b))}{\sharp \mbox{Aut}(\Gamma)}[G^0_{p(b)}]= \sum_{b \in B(\Gamma), f \in F}\frac{c(p(b);f)}{\sharp \mbox{Aut}(\Gamma)}[G^0_{p(b)}]. $$
Let $G(\Gamma)$ be the set of injections $g$ from the set $H(\Gamma)$ of half-edges of $\Gamma$ to $$\{(i,j) \in \{1,\dots, 2n\} \times \NN; 1 \leq j \leq g(i)\}$$
that map the three half-edges of a vertex to three images with the same first coordinate, and that induce a bijection from $V(\Gamma)$ to $\{1,\dots, 2n\}$. An injection $g$ of $G(\Gamma)$ provides a natural bijection $b(g)$
of $B(\Gamma)$ and a map $f(g)$ of $F$ such that $g(c)=(b_1(c),f(g) \circ b(g)(c))$. Furthermore, such a $g$ orders the three half-edges of a vertex, and hence provides an orientation $o(g)$ of $\Gamma$.
$$\sum_{p\in P(\Gamma)}c(p)[G^0_p]=\sum_{g \in G(\Gamma)}\frac{c(p(b(g));f(g))}{\sharp \mbox{Aut}(\Gamma)}[(\Gamma,o(g))].$$
Let $g \in G(\Gamma)$, its first coordinate $b_1(g)$ induces a bijection from $V(\Gamma)$ to $\{1,\dots, 2n\}$. Number the three half-edges of any vertex $w$
of $\Gamma$
with a bijection $b(w) \colon v^{-1}(w) \rightarrow \{1,2,3\}$, arbitrarily.
This orients $\Gamma$ and equips each injection $g \in G(\Gamma)$ with a sign that is $+1$
if $o(g)$ coincides with this orientation of $\Gamma$ (except for an even number of vertices) and $(-1)$ otherwise. Furthermore, $g$ provides summands of 
$$\CI(A_i,A_i^{\prime})=\sum_{g_i \colon\{1,2,3\} \rightarrow \{1,2, \dots, g(i)\}}\CI_{A_iA_i^{\prime}}(a_{g_i(1)}^i,a_{g_i(2)}^i,a_{g_i(3)}^i)z_{g_i(1)}^i \otimes z_{g_i(2)}^i \otimes z_{g_i(3)}^i$$
where $g(b(b_1(g)^{-1}(i))^{-1}(j))=(i,g_i(j))$.
Note that the sign of an injection $g$ is $+1$ if the number of vertices $b_1(g)^{-1}(i)$ where the cyclic order induced by $g_i$ does  not coincide with the cyclic order induced
by $b(b_1(g)^{-1}(i))$ is even, and $(-1)$, otherwise.
This shows that for any bijection $\sigma$ from $V(\Gamma)$ to $\{1,\dots, 2n\}$,
$$\sum_{g \in G(\Gamma);b_1(g)=\sigma}{c(p(b(g));f(g))}[(\Gamma,o(g))]=lk((A_i,A_i^{\prime})_{i=1, \dots, 2n}; \Gamma ;\sigma)[\Gamma]$$
with the notation of \cite{al} or \cite{sumgen}.
\eop

\section{A direct proof of the formula for boundary links}
\setcounter{equation}{0}
\label{secproofboun}

\subsection{A Lagrangian-preserving surgery associated to a Seifert surface}
\label{sublagboun}

Let $\Sigma$ be an oriented Seifert surface of a knot $K$ in a manifold $M$. 
Consider an annular neighborhood $[-3,0] \times K$ of $(\{0\} \times K)=K =\partial \Sigma$ in $\Sigma$,
a small disk $D$ inside $]-2,-1[ \times K$, and an open disk $d$ in the interior of $D$.
Let $F=\Sigma \setminus d$. Let $h_{F}$ be the composition of the two left-handed Dehn twists on $F$
along $c=\partial D$ and $K_2=\{-2\} \times K$ with the right-handed one along $K_1=\{-1\} \times K$.

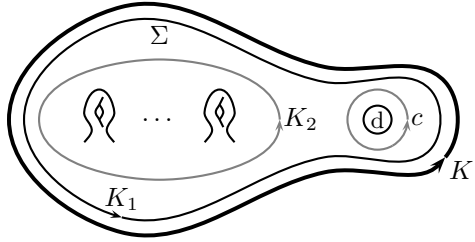
\begin{center}
\begin{pspicture}[.4](-.5,-.3)(6.3,3.3)
\psecurve[linewidth=1.5pt]{->}(4.2,.8)(5.8,1)(5.8,2)(4.2,2.2)(1.5,3)(0,1.5)(1.5,0)(4.2,.8)(5.8,1)(5.8,2)
\psecurve{->}(.2,1.5)(1.5,.2)(4.2,1)(5.6,1.1)(5.6,1.9)(4.2,2)(1.5,2.8)(.2,1.5)(1.5,.2)(5.6,1)
\psecurve[linecolor=gray]{->}(2,.5)(3.6,1.5)(2,2.3)(.4,1.5)(2,.7)(3.6,1.5)(2,2.5)
\pscircle[fillcolor=gray](4.9,1.5){.2}
\rput(4.9,1.5){\small d}
\psarc[linecolor=gray](4.9,1.5){.4}{0}{300}
\psarc[linecolor=gray]{->}(4.9,1.5){.4}{-60}{0}
\pscurve(1,1.2)(1.1,1.35)(1,1.6)(1.2,1.9)(1.4,1.6)(1.3,1.35)(1.4,1.2)
\pscurve(1.2,1.7)(1.25,1.6)(1.2,1.5)
\pscurve(1.25,1.8)(1.2,1.7)(1.15,1.6)(1.2,1.5)(1.25,1.4)
\pscurve(2.6,1.2)(2.7,1.35)(2.6,1.6)(2.8,1.9)(3,1.6)(2.9,1.35)(3,1.2)
\pscurve(2.8,1.7)(2.85,1.6)(2.8,1.5)
\pscurve(2.85,1.8)(2.8,1.7)(2.75,1.6)(2.8,1.5)(2.85,1.4)
\rput(2,1.5){\dots}
\rput[l](3.65,1.5){$K_2$}
\rput[l](5.35,1.5){$c$}
\rput[lt](5.85,.95){$K$}
\rput[b](1.5,.3){$K_1$}
\rput[b](2,2.5){$\Sigma$}
\end{pspicture}
\end{center}

See $F$ as $F \times \{0\}$ in the boundary of the handlebody 
$A_{F}=F \times [-1,0]$ of $M$.
Extend $h_F$
to a homeomorphism $h_A$ of $\partial A_{F}$ by defining it as the identity
outside $F \times \{0\}$.

Let $A_{F}^{\prime}$ be a copy of $A_{F}$.
Identify $\partial A^{\prime}_F$ with $\partial A_{F}$ with 
$$h_{A} \colon \partial A^{\prime}_F \rightarrow \partial A_{F}.$$

Define the {\em surgery associated to $\Sigma$\/} as the surgery 
associated with $(A_{F}, A_{F}^{\prime})$ (or $(A_{F}, A_{F}^{\prime}; h_{A})$).
If $\iota$ denotes the embedding from $\partial A_{F}$ to $M$. This surgery
replaces 
$$M =\left( M \setminus \mbox{Int}(A_{F}) \right) \cup_{\iota}  A_{F}$$
by
$$M_F=\left(M \setminus \mbox{Int}(A_{F}) \right) \cup_{\iota h_{A}}  A_{F}^{\prime}.$$

\begin{proposition}
With the notation above, the surgery $(A_{F}, A_{F}^{\prime})$ associated to $\Sigma$ is a Lagrangian-preserving surgery with the following properties.
There is a homeomorphism from $M_F$ to $M$ 
\begin{itemize}
\item that extends the identity of
$$ M \setminus \left([-3,0] \times K \times [-1,0] \right),$$
\item that transforms a curve going through $d \times [-1,0]$ by a band sum with $K$,
\item that transforms a $0$-framed meridian $m$ of $K$ going through $d \times [-1,0]$
into a $0$-framed copy of $K$ isotopic to the framed curve $h_A^{-1}(m)$ of the following figure.
\end{itemize}
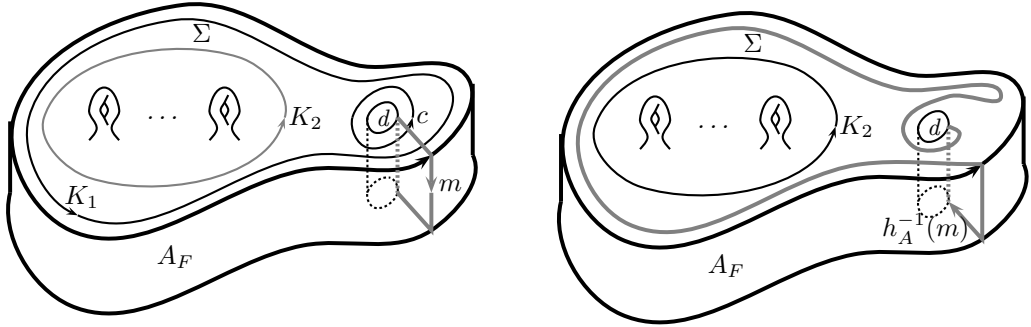
\begin{center}
\begin{pspicture}[.4](-.2,-1)(6.2,3.1)
\psecurve[linewidth=1.5pt]{->}(3.85,.8)(5.55,1)(6.05,2)(4.55,2.2)(2.25,3)(0,1.5)(.75,0)(3.85,.8)(5.55,1)(6.05,2)
\psecurve{->}(.2,1.5)(.85,.2)(3.95,1)(5.4,1.1)(5.8,1.9)(4.45,2)(2.15,2.8)(.2,1.5)(.85,.2)(3.95,1)
\psecurve[linecolor=gray]{->}(1,.7)(3.6,1.5)(3,2.3)(.4,1.5)(1,.7)(3.6,1.5)(3,2.3)
\psccurve[fillcolor=gray](5.1,1.5)(5,1.7)(4.7,1.5)(4.8,1.3)
\rput(4.9,1.5){\small d}
\psecurve(4.7,1.1)(5.3,1.5)(5.1,1.9)(4.5,1.5)(4.7,1.1)
\psecurve{->}(5.1,1.9)(4.5,1.5)(4.7,1.1)(5.3,1.5)(5.1,1.9)
\pscurve(1,1.2)(1.1,1.35)(1,1.6)(1.2,1.9)(1.4,1.6)(1.3,1.35)(1.4,1.2)
\pscurve(1.2,1.7)(1.25,1.6)(1.2,1.5)
\pscurve(1.25,1.8)(1.2,1.7)(1.15,1.6)(1.2,1.5)(1.25,1.4)
\pscurve(2.6,1.2)(2.7,1.35)(2.6,1.6)(2.8,1.9)(3,1.6)(2.9,1.35)(3,1.2)
\pscurve(2.8,1.7)(2.85,1.6)(2.8,1.5)
\pscurve(2.85,1.8)(2.8,1.7)(2.75,1.6)(2.8,1.5)(2.85,1.4)
\rput(2,1.5){\dots}
\rput[l](3.65,1.5){$K_2$}
\rput[l](5.35,1.5){$c$}
\rput[b](.9,.3){$K_1$}
\rput[b](2.5,2.5){$\Sigma$}
\psline[linewidth=1.5pt](6.1,1.95)(6.1,1)
\psline[linewidth=1.5pt](-.05,1.45)(-.05,.5)
\psccurve[linestyle=dashed, dash=1pt 1pt](5.1,.5)(5,.7)(4.7,.5)(4.8,.3)
\psline[linestyle=dashed, dash=1pt 1pt](4.7,1.5)(4.7,.5)
\psecurve[linewidth=1.5pt](2.25,2)(-.05,.5)(.75,-1)(3.85,-.2)(5.55,0)(6.1,1)(4.55,1.2)
\rput[b](2.15,-.5){$A_F$}
\psline[linewidth=1.5pt,linecolor=gray, linestyle=dashed, dash=1pt 1pt](5.1,1.5)(5.1,.5)
\psline[linewidth=1.5pt,linecolor=gray]{->}(5.1,1.5)(5.55,1)(5.55,.5)
\psline[linewidth=1.5pt,linecolor=gray]{-}(5.55,.5)(5.55,0)(5.1,.5)
\rput[l](5.65,.6){$m$}
\end{pspicture}
\begin{pspicture}[.4](-1,-.8)(6.2,3.1) 
\psecurve[linewidth=1.5pt]{->}(3.85,.8)(5.55,1)(6.05,2)(4.55,2.2)(2.25,3)(0,1.5)(.75,0)(3.85,.8)(5.55,1)(6.05,2)
\psecurve{->}(1,.7)(3.6,1.5)(3,2.3)(.4,1.5)(1,.7)(3.6,1.5)(3,2.3)
\psccurve[fillcolor=gray](5.1,1.5)(5,1.7)(4.7,1.5)(4.8,1.3)
\rput(4.9,1.5){\small d}
\pscurve(1,1.2)(1.1,1.35)(1,1.6)(1.2,1.9)(1.4,1.6)(1.3,1.35)(1.4,1.2)
\pscurve(1.2,1.7)(1.25,1.6)(1.2,1.5)
\pscurve(1.25,1.8)(1.2,1.7)(1.15,1.6)(1.2,1.5)(1.25,1.4)
\pscurve(2.6,1.2)(2.7,1.35)(2.6,1.6)(2.8,1.9)(3,1.6)(2.9,1.35)(3,1.2)
\pscurve(2.8,1.7)(2.85,1.6)(2.8,1.5)
\pscurve(2.85,1.8)(2.8,1.7)(2.75,1.6)(2.8,1.5)(2.85,1.4)
\rput(2,1.5){\dots}
\rput[l](3.65,1.5){$K_2$}
\rput[b](2.5,2.5){$\Sigma$}
\psline[linewidth=1.5pt](6.1,1.9)(6.1,1)
\psline[linewidth=1.5pt](-.05,1.4)(-.05,.5)
\psccurve[linestyle=dashed, dash=1pt 1pt](5.1,.5)(5,.7)(4.7,.5)(4.8,.3)
\psline[linestyle=dashed, dash=1pt 1pt](4.7,1.5)(4.7,.5)
\psecurve[linewidth=1.5pt](2.25,2)(-.05,.5)(.75,-1)(3.85,-.2)(5.55,0)(6.1,1)(4.55,1.2)
\rput[b](2.15,-.5){$A_F$}
\psline[linewidth=1.5pt, linestyle=dashed, dash=1pt 1pt,linecolor=gray](5.1,1.5)(5.1,.5)
\pscurve[linewidth=1.5pt,linecolor=gray](5.1,1.5)(5.25,1.4)(4.7,1.2)(4.5,1.5)(5.1,1.9)(5.8,1.9)(4.45,2)(2.15,2.8)(.2,1.5)(.85,.2)(3.95,1)(5.4,1)(5.55,1)
\psline[linewidth=1.5pt,linecolor=gray]{->}(5.55,1)(5.55,0)(5.1,.5)
\rput(4.8,.15){$h_A^{-1}(m)$}
\end{pspicture}
\end{center}
\end{proposition}
\bp
Observe that $h_F$ extends to $\Sigma \times [-1,0]$ as 
$$\begin{array}{llll} h \colon &\Sigma \times [-1,0] & \rightarrow & \Sigma \times [-1,0]\\
&(\sigma,t) & \mapsto & h(\sigma,t)=(h_t(\sigma),t)
\end{array}$$
where $h_0$ is the extension of $h_F$ by the identity on $d$ that is isotopic to the identity,\\
$h_{-1}$ is the identity of $\Sigma$,\\
$h_t$ coincides with the identity outside $[-5/2,-1/2]\times K(S^1)$,\\
and $h_t$ is defined as follows on $[-5/2,-1/2]\times K(S^1)$.\\
\noindent $\bullet$ When $t \leq -1/2$, then $h_t$ describes the following isotopy between $(h_{-1}=\mbox{identity})$ and 
the composition $h_{-1/2}$ of the left-handed Dehn twist along $K_2$ located on $[-5/2,-2]\times K(S^1)$ and the right-handed Dehn twist along $K_1$ located on $[-1,-1/2]\times K(S^1)$,
$$\begin{array}{lll}
h_t(u,K(z)) & = \left(u,K\left(z \exp\left(i(2t+2)(4\pi (u+5/2))\right)\right)\right) & \mbox{if} \; u \leq -2 \\
h_t(u,K(z)) & = \left(u,K\left(z \exp\left(i(2t+2)(2\pi)\right)\right)\right) & \mbox{if} \; -2 \leq u \leq -1 \\
h_t(u,K(z)) & = \left(u,K\left(z \exp\left(-i(2t+2)(4\pi (u+1/2))\right)\right)\right) & \mbox{if} \; u \geq -1.
\end{array}$$
\noindent $\bullet$ When $t \geq -1/2$, then $h_t$ coincides with $h_{-1/2}$ outside the disk $D$ whose elements will be written as $D(z\in \CC)$, with $|z|\leq 1$. The elements of $d$ will be the $D(z)$ for $|z|<1/2$. On $D$, $h_t$ will describe the isotopy between the identity and  
the composition $h_{0}$ of the left-handed Dehn twist along $\partial D$ located on $\{D(z);1/2 \leq |z| \leq 1\}$ and a negative twist of $d$.
$$\begin{array}{lll}
h_t(u,K(z)) & = \left(u,K\left(z \exp\left(i(4\pi (u+5/2))\right)\right)\right) & \mbox{if} \; u \leq -2 \\
h_t(u,K(z)) & = \left(u,K(z)\right) & \mbox{if} \; -2 \leq u \leq -1, (u,K(z)) \notin D \\
h_t(u,K(z)) & = \left(u,K\left(z \exp\left(-i(4\pi (u+1/2))\right)\right)\right) & \mbox{if} \; u \geq -1 \\
h_t(z \in D) & = z \exp\left(i\pi(2t+1)4(|z|-1)\right) & \mbox{if} \; |z| \geq 1/2 \\
h_t(z \in D) & = z \exp\left(-2i\pi(2t+1)\right) & \mbox{if} \; |z| \leq 1/2.\\
\end{array}$$
Now, $M_F$ is naturally homeomorphic to 
$$\left(M \setminus \mbox{Int}(\Sigma \times [-1,0]) \right) \cup_{h_{|\partial(\Sigma \times [-1,0])}}  (\Sigma \times [-1,0])$$
that maps to $M$ by the identity outside $\Sigma \times [-1,0]$ and by $h$ on $\Sigma \times [-1,0]$, homeomorphically. Therefore, we indeed have a homeomorphism from $M_F$ to $M$ that is the identity outside $[-3,0]\times K \times [-1,0]$
and that maps $d \times [-1,0]$ to a cylinder that runs along $K$ before being negatively twisted. In particular, looking at the action of the homeomorphism on a framed arc
$x \times [-1,0]$ where $x$ is on the boundary of $d$ shows that the meridian $m$ with its framing induced by the boundary of $A_F$ is mapped to a curve isotopic to 
$h_A^{-1}(m)$ in a tubular neighborhood of $K$ with the framing induced by the boundary of $A_F$.

Now, $H_1(\partial A_{F})$ is generated by the generators of $H_1(\Sigma)\times \{0\}$, the generators of $H_1(\Sigma)\times \{-1\}$,
and the homology classes of $c=\partial D$ and $m$. Among them, only the class of $m$ could be affected by $h_A$, and it is not.
Therefore $h_A$ acts trivially on $H_1(\partial A_{F})$, and the defined surgery is an $LP$--surgery.
\eop

Let $F \times [-1,2]$ be an extension of the previous neighborhood of $F$, and let $B_F=F \times [1,2]$. 
Define the homeomorphism $h_B$ of $\partial B_F$ as the identity anywhere except 
on $F \times \{1\}$ where it coincides with the homeomorphism $h_F$ of $F$ with
the obvious identification.

Let $B_{F}^{\prime}$ be a copy of $B_{F}$.
Identify $\partial B^{\prime}_{F}$ with $\partial B_{F}$ with 
$$h_B \colon \partial B^{\prime}_{F} \rightarrow \partial B_{F}.$$

Define the {\em inverse surgery associated to $\Sigma$\/} as the surgery 
associated with $(B_{F}, B_{F}^{\prime})$ (or $(B_{F}, B_{F}^{\prime}; h_B)$).
Note that the previous study can be used for this surgery by using the central symmetry of $[-1,2]$.

Then, we have the following obvious lemma that justifies the terminology.

\begin{lemma}
With the notation above, performing the two surgeries $(B_{F}, B_{F}^{\prime})$ and $(A_{F}, A_{F}^{\prime})$ affects neither $M$ nor the curves in the complement
of $F \times [-1,2]$, while performing one of them changes a $0$-framed meridian of $K$ going through $d \times [-1,2]$
into a $0$-framed copy of $\pm K$.
\end{lemma} 
\eop

\begin{lemma}
\label{lemtriplag}
Let $(x_i, y_i)_{i=1, \dots, g}$ be a symplectic basis of $\Sigma$, then the tripod combination $T(\CI_{A_FA_F^{\prime}})$ associated
to the surgery
$(A_{F}, A_{F}^{\prime})$
is $$T(\CI_{A_FA_F^{\prime}})=-\sum_{i=1}^g 
\begin{pspicture}[0.4](-.05,-.1)(.9,.7) 
\psline{-}(0.05,.3)(.45,.6)
\psline{*-}(0.05,.3)(.45,.3) 
\psline{-}(0.05,.3)(.45,0)
\rput[l](.55,0){$c$}
\rput[l](.55,.3){$x_i$}
\rput[l](.55,.6){$y_i$}
\end{pspicture}.$$

For a curve $c$ of $F$, let $c^+$ denote $c \times \{1\}$.
The tripod combination $T(\CI_{B_FB_F^{\prime}})$ associated
to the surgery
$(B_{F}, B_{F}^{\prime})$
is $$T(\CI_{B_FB_F^{\prime}})=\sum_{i=1}^g 
\begin{pspicture}[0.4](-.05,-.1)(1.45,.7) 
\psline{-}(0.05,.3)(.45,.6)
\psline{*-}(0.05,.3)(.85,.3) 
\psline{-}(0.05,.3)(.45,0)
\rput[l](.55,0){$c^+$}
\rput[l](.95,.3){$x_i^+$}
\rput[l](.55,.6){$y_i^+$}
\end{pspicture}.$$
\end{lemma}
\bp For a curve $c$ of $F$, $c^-$ denotes $c \times \{-1\}$. 
Use the basis $\left(m,(x_i - x_i^-,y_i-y_i^-)_{i=1, \dots, g}\right)$
of the Lagrangian of $A_F$ 
to compute the intersection form of $(A_F \cup -A^{\prime}_F)$. Its dual basis is $\left(c,(y_i, - x_i)_{i=1, \dots, g}\right)$.
Note that the only curve of the Lagrangian basis that is modified by $h_A$ is $m$, and that $h_A(m)=mK_2^{-1}$.
The isomorphism $\partial_{MV}^{-1}$ from $\CL_{A_F}$ to $H_2(A_F \cup -A^{\prime}_F)$ satisfies
$$\begin{array}{ll}
\partial_{MV}^{-1}(x_i - x_i^-) &= S(x_i)=-(x_i \times [-1,0]) \cup (x_i \times [-1,0] \subset A^{\prime}_F)\\
\partial_{MV}^{-1}(y_i - y_i^-) &=S(y_i)=-(y_i \times [-1,0]) \cup (y_i \times [-1,0] \subset A^{\prime}_F) \\
\partial_{MV}^{-1}(m) &=S_A(m)=D_m -(\Sigma \setminus (]-2,0] \times K))\cup (-D_m \subset A^{\prime}_F) 
\end{array}.$$

Since $x_i$ intersects only $y_i$, $S(x_i)$ intersects only $S(y_i)$ and $S_A(m)$. The algebraic intersection of
$S(x_i)$, $S(y_i)$ and $S_A(m)$ is $-1$.

For the surgery $(B_{F}, B_{F}^{\prime})$, $S_B(m)=D_m + \Sigma \setminus (]-2,0] \times K) \cup (-D_m \subset B^{\prime}_F)$, and the algebraic intersection of
$S(x_i)$, $S(y_i)$ and $S_B(m)$ is $1$.
\eop

\subsection{Proof of Theorem~\ref{thmfboun}}

\begin{remark}
For this proof, I could also have used the strategy of Section~\ref{secprooffas}.
But I prefer this  self-contained proof.
\end{remark}

First recall the following easy lemma that will be used several times.
\begin{lemma}
\label{lemvarlk}
The variation of the linking number of two knots $J$ and $K$ after a $p/q$-surgery
on a knot $V$ in a rational homology sphere $M$ is given by the following formula.
$$lk_{M_{(V;p/q)}}(J,K)=lk_M(J,K)-\frac{q}{p}lk_M(V,J)lk_M(V,K).$$
\end{lemma}
\eop

Let $(K_1,K_2, \dots, K_n)$ be a link where all the $K_i$ bound disjoint
oriented surfaces $\Sigma^i$. 
Consider an embedding of $\coprod_{i=1}^r \Sigma^i \times [-1,2]$.
Let $N=\{1,2, \dots, n\}$.
For $i\in N$, associate surfaces $F^i= \Sigma^i \setminus d^i$ and 
$LP$--surgeries $(A_i,A_i^{\prime})=(A_{F^i}, A_{F^i}^{\prime})$ and $(B_i,B_i^{\prime})=(B_{F^i}, B_{F^i}^{\prime})$
as in Subsection~\ref{sublagboun}. Let $U_i$ be a meridian of $K_i$
going through $d^i \times [-1,2]$, so that performing one of the two surgeries
transforms $U_i$ into $\pm K_i$ and performing both or none of them leaves $U_i$ unchanged.
Then $$\begin{array}{ll}[M;(K_i;p_i/q_i)]&=M_{(U_i;p_i/q_i)}-M_{(K_i;p_i/q_i)}\\
&=\frac12[M_{(U_i;p_i/q_i)};(A_i,A_i^{\prime}),(B_i,B_i^{\prime})].\end{array}$$
More generally, for $J \subset \{(A_i,A_i^{\prime}),(B_i,B_i^{\prime})\}_{i=1,\dots,n}$, 
$$(M_{(U_i;p_i/q_i)_{i \in N}})_J=M_{(K_i;p_i/q_i)_{i \in I(J)}} \sharp \sharp_{j \notin I(J)} L(p_j,-q_j)=M_{I(J)}$$ 
where $I(J)$ is the set of elements $i$ of $N$ such that $\sharp \left(J \cap \{(A_i,A_i^{\prime}),(B_i,B_i^{\prime})\}\right)$ is one.
Note that $(-1)^{\sharp J}= (-1)^{\sharp I(J)}$ and that for any subset $I$ of $N$ there are $2^n$ subsets $J$ of the
set of LP-surgeries such that $I(J)=I$. Thus
$$[M;(K_i;p_i/q_i)_{i \in N}]=\frac1{2^n}[M_{(U_i;p_i/q_i)_{i \in N}};(A_i,A_i^{\prime})_{i \in N},(B_i,B_i^{\prime})_{i \in N}].$$
In particular, we can apply Theorem~\ref{thmflag} to compute $Z_n([M;(K_i;p_i/q_i)_{i \in N}])$.

According to Lemma~\ref{lemtriplag},
the tripods associated
to the surgery
$(A_{F^i}, A_{F^i}^{\prime})$ and to the surgery
$(B_{F^i}, B_{F^i}^{\prime})$
are $-\sum_{i=1}^{g^i} 
\begin{pspicture}[0.4](-.05,-.1)(1.3,.8) 
\psline{-}(0.05,.3)(.45,.6)
\psline{*-}(0.05,.3)(.85,.3) 
\psline{-}(0.05,.3)(.45,0)
\rput[l](.55,0){$c^i$}
\rput[l](.95,.3){$x_j^i$}
\rput[l](.55,.6){$y_j^i$}
\end{pspicture}$
and $\sum_{i=1}^{g^i} 
\begin{pspicture}[0.4](-.05,-.1)(1.45,.8) 
\psline{-}(0.05,.3)(.45,.6)
\psline{*-}(0.05,.3)(.85,.3) 
\psline{-}(0.05,.3)(.45,0)
\rput[l](.55,0){$c^{i+}$}
\rput[l](.95,.3){$x_j^{i+}$}
\rput[l](.55,.6){$y_j^{i+}$}
\end{pspicture}$, respectively.
The only curve that links $c^i$ algebraically in $M_{(U_i;p_i/q_i)_{i\in N}}$ among those appearing in all the tripods
is $c^{i+}$ with a linking number $-q_i/p_i$. Therefore, these two must be paired together with this coefficient.
Theorem~\ref{thmfboun} follows when $r=n$. The case $r>n$ can be either deduced from the case $r=n$ or proved directly, it is easy.
\eop

\section{Some clasper calculus}
\setcounter{equation}{0}
\label{secclasper}

The proofs of Theorems~\ref{thmfas} and \ref{thmfasmu} will be given in Section~\ref{secprooffas}. They will rely on the current section, where we recall some known clasper calculus and where we show how to present algebraically split links $L=(K_i)_{i=1,\dots,n}$ by claspers so that the associated Seifert surfaces $\Sigma_i$ of the components $K_i$ in $M \setminus \left(\cup_{j \neq i}K_j \right) $ have minimal triple intersection, namely so that for any triple $(K_i,K_j,K_k)$ of components of $L$, the geometric triple intersection of the transverse surfaces $\Sigma_i$, $\Sigma_j$ and $\Sigma_k$ is made of $|\mu(K_i,K_j,K_k)|$ points. (This shows Proposition~\ref{proprealmil} that will be a direct corollary of Lemma~\ref{lemintclas} and Proposition~\ref{propclaspmu}.)

\subsection{Two ways of seeing surgeries on $Y$-graphs}
\label{subsectwoways}
Let $\Lambda$ be the graph embedded in the surface $\Sigma(\Lambda)$
shown below.
In the $3$--handlebody $(N=\Sigma(\Lambda) \times [-1,1])$,
the edges of $\Lambda$ are framed by a vector field normal
to $\Sigma(\Lambda)=\Sigma(\Lambda) \times\{0\}$.  
$\Sigma(\Lambda)$ is called a {\em framing surface\/} for $\Lambda$.

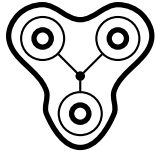
\begin{center}
\begin{pspicture}[.4](-2,0)(4,2)
\psccurve[linewidth=2pt](0,1.95)(-.32,1.82)(-.45,1.5)(.05,.8)(.05,.5)(.5,.05)(.95,.5)(.95,.8)(1.45,1.5)(1.32,1.82)(1,1.95)(.5,1.75)
\pscircle[linewidth=2pt](.5,.5){.15}
\pscircle[linewidth=2pt](1,1.5){.15}
\pscircle[linewidth=2pt](0,1.5){.15}
\pscircle(.5,.5){.3}
\pscircle(1,1.5){.3}
\pscircle(0,1.5){.3}
\psline{*-}(.5,1)(.5,.8)
\psline(.21,1.29)(.5,1)(.79,1.29)
\end{pspicture}
\end{center}

A {\em $Y$-graph\/} in $M$ is the isotopy class of an 
embedding $\phi$ of $N$ (or $\Sigma(\Lambda)$) into $M$. 
Such an isotopy class is determined by the framed image of the framed
unoriented graph $\Lambda $ under $\phi$. 
A {\em leaf} of a $Y$-graph $\phi$ 
is the image under $\phi$ of a simple loop of our graph $\Lambda$.
An {\em edge} of $\phi$ 
is an edge of $\phi(\Lambda)$ that is not a leaf. 
With this terminology, a $Y$-graph  has  
three 
edges and three leaves.

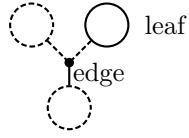
\begin{center}
\begin{pspicture}[.4](-2,0)(4,1.9)
\pscircle[linestyle=dashed, dash=2pt 1pt](.5,.4){.3}
\pscircle(1,1.5){.3}
\pscircle[linestyle=dashed, dash=2pt 1pt](0,1.5){.3}
\psline{*-}(.5,1)(.5,.7)
\psline[linestyle=dashed, dash=2pt 1pt](.21,1.29)(.5,1)(.79,1.29)
\rput[l](.55,.85){edge}
\rput[l](1.5,1.5){leaf}
\end{pspicture}
\end{center}

The surgery on such a $Y$-graph can be defined in several equivalent ways.

Originally, it was defined by Matveev in \cite{matv} and named {\em Borromeo transformation\/} 
as the effect of the surgery
on the following $6$-component framed link in the framed neighborhood 
of the $Y$-graph.

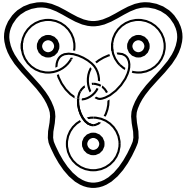
\begin{center}
\begin{pspicture}[.4](-2,-.2)(4,2.2)
\psccurve[linewidth=2pt](-.2,2.25)(-.52,2.12)(-.65,1.8)(-.05,.85)(-.05,.4)(.5,-.15)(1.05,.4)(1.05,.85)(1.65,1.8)(1.52,2.12)(1.2,2.25)(.5,2)
\pscircle[linewidth=2pt](.5,.4){.15}
\pscircle[linewidth=2pt](1.1,1.7){.15}
\pscircle[linewidth=2pt](-.1,1.7){.15}
\psarc(.5,.4){.35}{90}{0}
\psarc(1.1,1.7){.35}{-145}{180}
\psarc(-.1,1.7){.35}{-30}{180}
\psarc(.5,1){.2}{0}{180}
\psecurve(.5,1.2)(.3,1)(.5,.65)(.7,1)(.5,1.2)
\psccurve[border=1pt](.5,1.45)(.45,1.125)(.6,1)(.95,1.55)
\psarc[border=1pt](.5,.4){.35}{0}{90}
\psarc[border=1pt](1.1,1.7){.35}{180}{225}
\psccurve[border=1pt](.5,1.45)(.55,1.125)(.3,1)(.05,1.55)
\psarc[border=1pt](-.1,1.7){.35}{180}{-30}
\psarc[border=1pt](.5,1){.2}{60}{90}
\psarc[border=1pt](.5,1){.2}{150}{180}
\psarc(-.1,1.7){.35}{170}{190}
\psarc(1.1,1.7){.35}{175}{190}
\psecurve[border=1pt](.5,1.2)(.3,1)(.5,.65)(.7,1)
\psecurve(.5,1.2)(.3,1)(.5,.65)(.7,1)
\end{pspicture}
\end{center}
The framing of the link is induced by the framing of the surface.

We shall prove the following proposition.

\begin{proposition}
\label{propborlag}
The above surgery
is equivalent to the surgery $(A_F,A_F^{\prime})$ associated to the following subsurface $F$ of 
$\Sigma(\Lambda) \times [-1,1]$, with respect to the notation of Subsection~\ref{sublagboun}. 

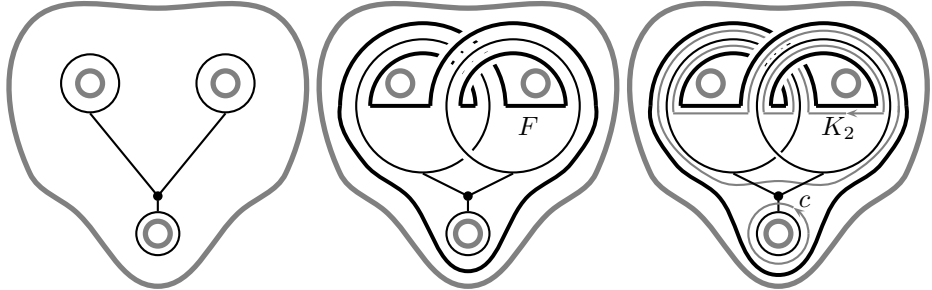
\begin{center}
\begin{pspicture}[.4](-.1,-1)(3.7,3)
\psccurve[linewidth=2pt, linecolor=gray](0,2.7)(0,.7)(.6,.2)(1.7,-.8)(2.8,.2)(3.4,.7)(3.4,2.7)(1.7,2.9)
\psline(1.7,.4)(.8,1.5)
\psline(1.7,.4)(2.6,1.5)
\psline{*-}(1.7,.4)(1.7,.2)
\pscircle(2.6,1.9){.4}
\pscircle(.8,1.9){.4}
\pscircle[linewidth=2pt, linecolor=gray](2.6,1.9){.2}
\pscircle[linewidth=2pt, linecolor=gray](.8,1.9){.2}
\pscircle[linewidth=2pt, linecolor=gray](1.7,-.1){.2}
\pscircle(1.7,-.1){.3}
\end{pspicture}
\begin{pspicture}[.4](-.3,-1)(3.7,3)
\psarc[linewidth=1.5pt](1.1,1.6){1.1}{0}{180}
\psarc[linewidth=1.5pt](1.1,1.6){.7}{0}{180}
\pscircle(1.1,1.6){.9}
\psarc[linewidth=1.5pt,border=2pt](2.3,1.6){1.1}{0}{180}
\psarc[linewidth=1.5pt,border=2pt](2.3,1.6){.7}{0}{180}
\pscircle[border=2pt](2.3,1.6){.9}
\psline[linewidth=1.5pt](.4,1.6)(1.2,1.6)
\psline[linewidth=1.5pt](1.6,1.6)(1.8,1.6)
\psline[linewidth=1.5pt](2.2,1.6)(3,1.6)
\psecurve[linewidth=1.5pt](1.1,2.7)(0,1.6)(.33,.83)(.9,.4)(1.7,-.6)(2.5,.4)(3.07,.83)(3.4,1.6)(2.3,2.7)
\psccurve[linewidth=2pt, linecolor=gray](0,2.7)(0,.7)(.6,.2)(1.7,-.8)(2.8,.2)(3.4,.7)(3.4,2.7)(1.7,2.9)
\psline(1.7,.4)(1.1,.7)
\psline(1.7,.4)(2.3,.7)
\psline{*-}(1.7,.4)(1.7,.2)
\pscircle[linewidth=2pt, linecolor=gray](2.6,1.9){.2}
\pscircle[linewidth=2pt, linecolor=gray](.8,1.9){.2}
\pscircle[linewidth=2pt, linecolor=gray](1.7,-.1){.2}
\pscircle(1.7,-.1){.3}
\rput[t](2.5,1.45){$F$}
\end{pspicture}
\begin{pspicture}[.4](-.3,-1)(3.5,3)
\psarc[linewidth=1.5pt](1.1,1.6){1.1}{0}{180}
\psarc[linewidth=1.5pt](1.1,1.6){.7}{0}{180}
\psarc[linecolor=gray](1.1,1.6){1}{-5}{185}
\psarc[linecolor=gray](1.1,1.6){.8}{-6}{186}
\pscircle(1.1,1.6){.9}
\psarc[linewidth=1.5pt,border=2pt](2.3,1.6){1.1}{0}{180}
\psarc[linewidth=1.5pt,border=2pt](2.3,1.6){.7}{0}{180}
\psarc[linecolor=gray,border=1pt](2.3,1.6){1}{-5}{185}
\psarc[linecolor=gray,border=1pt](2.3,1.6){.8}{-6}{186}
\pscircle[border=1pt](2.3,1.6){.9}
\psline[linewidth=1.5pt](.4,1.6)(1.2,1.6)
\psline[linewidth=1.5pt](1.6,1.6)(1.8,1.6)
\psline[linewidth=1.5pt](2.2,1.6)(3,1.6)
\psline[linecolor=gray](.3,1.5)(1.3,1.5)
\psline[linecolor=gray](1.5,1.5)(1.9,1.5)
\psline[linecolor=gray](2.6,1.5)(2.1,1.5)
\psline[linecolor=gray]{->}(3.1,1.5)(2.6,1.5)
\rput[t](2.5,1.45){$K_2$}
\psecurve[linewidth=1.5pt](1.1,2.7)(0,1.6)(.33,.83)(.9,.4)(1.5,-.6)(1.9,-.6)(2.5,.4)(3.07,.83)(3.4,1.6)(2.3,2.7)
\psecurve[linecolor=gray](1.1,2.6)(.1,1.6)(.9,.6)(1.7,.6)(2.5,.6)(3.3,1.6)(2.3,2.6)
\psccurve[linewidth=2pt, linecolor=gray](0,2.7)(0,.7)(.6,.2)(1.7,-.8)(2.8,.2)(3.4,.7)(3.4,2.7)(1.7,2.9)
\psline(1.7,.4)(1.1,.7)
\psline(1.7,.4)(2.3,.7)
\psline{*-}(1.7,.4)(1.7,.2)
\pscircle[linewidth=2pt, linecolor=gray](2.6,1.9){.2}
\pscircle[linewidth=2pt, linecolor=gray](.8,1.9){.2}
\pscircle[linewidth=2pt, linecolor=gray](1.7,-.1){.2}
\pscircle(1.7,-.1){.3}
\psarc[linecolor=gray](1.7,-.1){.4}{60}{360}
\psarc[linecolor=gray]{->}(1.7,-.1){.4}{0}{60}
\rput[b](2.05,.25){$c$}
\end{pspicture}
\end{center}
\end{proposition}

Let $G\subset M$ be a \mbox{$Y$-graph.} A leaf $l$ of a $Y$-component of $G$
is {\em trivial} if $l$  bounds an 
embedded disc that induces the framing of $l$,
in $M\setminus G$. It is easy to see that with both definitions, 
performing the surgery on a $Y$-graph with a trivial leaf does not change the
ambient manifold. More precisely, the following lemma is proved in \cite{ggp}, for
the first definition.

\begin{lemma}[\mbox{\cite[Lemma 2.1]{ggp}}] \label{lemysur}
Let $M$ be an oriented $3$--manifold (with possible boundary). Let $G$ be a $Y$-graph in 
$M$ with a trivial leaf 
that bounds a disc $D$ in $M\setminus G$. Then
\begin{itemize}
\item for any framed graph $T_0$ in $M\setminus G$ that does not meet
$D$, the pair $(M_G,T_0)$ is diffeomorphic to the pair $(M,T_0)$.
\item If $T$ is a framed graph 
in $M\setminus G$ that meets $\IT(D)$ at exactly one
point, then the pair $(M_G,T)$ is diffeomorphic to the pair $(M,T_G)$,
where $T_G$ is the framed graph in $M$ below. 
\end{itemize}
\begin{center}
\begin{pspicture}[.2](-.2,.2)(2.4,2) 
\psline[linewidth=1.5pt]{<-}(0,.5)(2,.5)
\psarc[border=2pt](1,.5){.2}{-160}{160}
\psarc(.5,1.6){.3}{140}{95}
\psarc[linestyle=dotted](.5,1.6){.3}{95}{140}
\psarc(1.5,1.6){.3}{85}{40}
\psarc[linestyle=dotted](1.5,1.6){.3}{40}{85}
\psline(.5,1.3)(1,1)(1.5,1.3)
\psline{*-}(1,1)(1,.7)
\rput[l](1.1,.9){$G$}
\rput[rb](-.05,.55){$T$} 
\end{pspicture}
$\longrightarrow$
\begin{pspicture}[.1](-2,.2)(2.4,2)
\psset{xunit=.6cm,yunit=.6cm} 
\psarc[linewidth=1.5pt](1.1,1.6){.66}{0}{95}
\psarc[linewidth=1.5pt,linestyle=dotted](1.1,1.6){.66}{95}{140}
\psarc[linewidth=1.5pt](1.1,1.6){.66}{140}{180}
\psarc[linewidth=1.5pt](1.1,1.6){.42}{0}{95}
\psarc[linewidth=1.5pt,linestyle=dotted](1.1,1.6){.42}{95}{140}
\psarc[linewidth=1.5pt](1.1,1.6){.42}{140}{180}
\psarc[linewidth=1.5pt,border=2pt](2.3,1.6){.66}{0}{40}
\psarc[linewidth=1.5pt,border=2pt,linestyle=dotted](2.3,1.6){.66}{40}{85}
\psarc[linewidth=1.5pt,border=2pt](2.3,1.6){.66}{85}{180}
\psarc[linewidth=1.5pt,border=2pt](2.3,1.6){.42}{0}{40}
\psarc[linewidth=1.5pt,border=2pt,linestyle=dotted](2.3,1.6){.42}{40}{85}
\psarc[linewidth=1.5pt,border=2pt](2.3,1.6){.42}{85}{180}
\psecurve[linewidth=1.5pt]{->}(1.1,2.7)(0,1.6)(.33,.93)(1.1,.6)(.9,.5)(.3,.5)(-.2,.5)
\psecurve[linewidth=1.5pt](2.3,2.7)(3.4,1.6)(3.07,.93)(2.3,.6)(2.5,.5)(3.1,.5)(3.7,.5)
\rput[rb](.25,.55){$T_G$} 
\psline[linewidth=1.5pt](.4,1.6)(1.2,1.6)
\psline[linewidth=1.5pt](1.6,1.6)(1.8,1.6)
\psline[linewidth=1.5pt](2.2,1.6)(3,1.6)
\end{pspicture}
\end{center}

\end{lemma}

Now, it is proved in \cite[Proof of Lemma 4.6]{al}, that this property fully determines the surgery.
Therefore, since this property is also true for the second definition, the two definitions coincide and Proposition~\ref{propborlag} is proved.
In particular, the second definition has the same symmetries as the first one obviously has.

This definition does not depend on the orientation of $\Sigma(\Lambda)$. 
Nevertheless, we shall sometimes need orientations of our $Y$-graphs. An {\em orientation\/} of a $Y$-graph
is an orientation of its three leaves, together with a cyclic order on the $3$--element set they form, induced by an orientation of $\Sigma(\Lambda)$ as in the figure (everything turns counterclockwise).

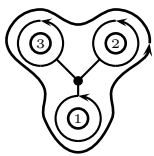
\begin{center}
\begin{pspicture}[.4](-2,0)(3,2)
\psecurve[linewidth=1pt]{->}(.95,.9)(1.45,1.5)(1.32,1.82)(1,1.95)(.5,1.75)(0,1.95)(-.32,1.82)(-.45,1.5)(.05,.9)(.05,.5)(.5,.05)(.95,.5)(.95,.9)(1.45,1.5)(1.32,1.82)
\psarc{->}(.5,.5){.3}{-90}{90}
\psarc{->}(1,1.5){.3}{-90}{90}
\psarc{->}(0,1.5){.3}{-90}{90}
\psarc(.5,.5){.3}{90}{-90}
\psarc(1,1.5){.3}{90}{-90}
\psarc(0,1.5){.3}{90}{-90}
\pscircle[linewidth=1pt](.5,.5){.15}
\pscircle[linewidth=1pt](1,1.5){.15}
\pscircle[linewidth=1pt](0,1.5){.15}
\psline{*-}(.5,1)(.5,.8)
\psline(.21,1.29)(.5,1)(.79,1.29)
\rput(.5,.5){\tiny 1}
\rput(1,1.5){\tiny 2}
\rput(0,1.5){\tiny 3}
\end{pspicture}
\end{center}

An {\em $n$--component $Y$-link\/}
$G\subset M$ is an embedding of the disjoint  union of $n$ copies of $N$ 
into $M$ up to isotopy. The  
$Y$-surgery along a $Y$-link $G$ is defined as the surgery
along each  $Y$-component of $G$. The resulting manifold is 
denoted by $M_G$.

\subsection{Some clasper calculus}

Recall the following equivalences between surgeries inside handlebodies -that can be themselves 
embedded in any $3$-manifold in an arbitrary way-.
The first one is move $Y_3$ in \cite{ggp}, as rectified by Emmanuel Auclair in his thesis \cite{auc}.
\begin{lemma}[\cite{auc}]
\label{lemY3}
The surgeries on the following two $Y$-links are equivalent.
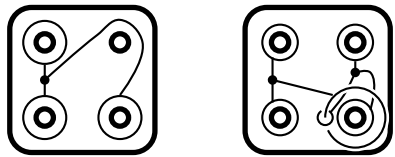
\begin{center}
\begin{pspicture}[.4](-.5,-.2)(2.5,2.2)
\psframe[linewidth=2pt, framearc=.3](0,0)(2,2)
\pscircle[linewidth=2pt](.5,.5){.15}
\pscircle[linewidth=2pt](.5,1.5){.15}
\pscircle[linewidth=2pt](1.5,.5){.15}
\pscircle[linewidth=2pt](1.5,1.5){.15}
\pscircle(.5,.5){.3}
\pscircle(.5,1.5){.3}
\pscircle(1.5,.5){.3}
\psline(.5,.8)(.5,1.2)
\pscurve{*-}(.5,1)(1.2,1.6)(1.5,1.8)(1.8,1.5)(1.5,.8)
\end{pspicture} 
\begin{pspicture}[.4](-.5,-.2)(2.5,2.2)
\psframe[linewidth=2pt, framearc=.3](0,0)(2,2)
\pscircle[linewidth=2pt](.5,.5){.15}
\pscircle[linewidth=2pt](.5,1.5){.15}
\pscircle[linewidth=2pt](1.5,.5){.15}
\pscircle[linewidth=2pt](1.5,1.5){.15}
\pscircle(1.5,1.5){.25}
\pscircle(.5,.5){.25}
\pscircle(.5,1.5){.25}
\pscircle(1.5,.5){.25}
\psline{-*}(1.5,1.25)(1.5,1.1)
\pscurve(1.5,1.1)(1.35,.8)(1.17,.57)
\pscurve(1.5,1.1)(1.7,1.1)(1.7,.6)
\psarc[border=1pt](1.5,.5){.4}{-155}{180}
\psline(.4,.7)(.4,1.3)
\psarc(1.1,.5){.1}{120}{60}
\psline{*-}(.4,1)(1.2,.8)
\end{pspicture}
\end{center}
\end{lemma}

\begin{lemma}[Theorem 3.1, move $Y_4$ in \cite{ggp}]
\label{lemY4}
The surgeries on the following two $Y$-links are equivalent.
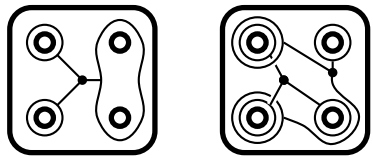
\begin{center}
\begin{pspicture}[.4](-.2,-.5)(2.5,2.2)
\psframe[linewidth=2pt, framearc=.3](0,0)(2,2)
\pscircle[linewidth=2pt](.5,.5){.15}
\pscircle[linewidth=2pt](.5,1.5){.15}
\pscircle[linewidth=2pt](1.5,.5){.15}
\pscircle[linewidth=2pt](1.5,1.5){.15}
\pscircle(.5,1.5){.25}
\pscircle(.5,.5){.25}
\psccurve(1.25,1)(1.2,1.5)(1.5,1.8)(1.8,1.5)(1.75,1)(1.8,.5)(1.5,.2)(1.2,.5)
\psline{-*}(1.25,1)(1,1)
\psline(.67,.67)(1,1)(.67,1.33)
\end{pspicture} 
\begin{pspicture}[.4](-.2,-.5)(2.5,2.2)
\psframe[linewidth=2pt, framearc=.3](0,0)(2,2)
\pscircle[linewidth=2pt](.5,1.5){.15}
\pscircle[linewidth=2pt](1.5,.5){.15}
\pscircle[linewidth=2pt](1.5,1.5){.15}
\pscircle(1.5,1.5){.25}
\pscircle(.5,.5){.25}
\pscircle(.5,1.5){.25}
\pscircle(1.5,.5){.25}
\pscircle(.5,.5){.35}
\psline(1.33,.67)(.85,1)(.67,1.33)
\pscircle[border=1pt](.5,1.5){.35}
\psline[border=1pt]{*-}(.85,1)(.67,.67)
\psline(.85,1.5)(1.5,1.1)(1.5,1.25)
\pscurve{*-}(1.5,1.1)(1.5,.95)(1.85,.5)(1.5,.15)(1.1,.4)(.85,.5)
\pscircle[linewidth=2pt](.5,.5){.15}
\end{pspicture}
\end{center}
\end{lemma}

The two equivalences above easily imply the following one.
\begin{lemma}
\label{lemY5}
The surgeries on the following two $Y$-links are equivalent.
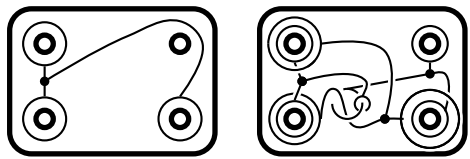
\begin{center}
\begin{pspicture}[.4](-.5,-.2)(3,2.2)
\psframe[linewidth=2pt, framearc=.3](0,0)(2.8,2)
\pscircle[linewidth=2pt](.5,.5){.15}
\pscircle[linewidth=2pt](.5,1.5){.15}
\pscircle[linewidth=2pt](2.3,.5){.15}
\pscircle[linewidth=2pt](2.3,1.5){.15}
\pscircle(.5,.5){.3}
\pscircle(.5,1.5){.3}
\pscircle(2.3,.5){.3}
\psline(.5,.8)(.5,1.2)
\pscurve{*-}(.5,1)(1.6,1.6)(2.3,1.8)(2.6,1.5)(2.3,.8)
\end{pspicture} 
\begin{pspicture}[.4](-.2,-.2)(3,2.2)
\psframe[linewidth=2pt, framearc=.3](0,0)(2.8,2)
\pscurve(2.3,1.1)(2.5,1.1)(2.5,.6)
\psline(2.3,1.1)(1.2,.9)
\pscurve[border=1pt](.6,1)(1.45,1)(1.4,.8)
\pscircle[border=1pt](2.3,.5){.4}
\pscircle[linewidth=2pt](.5,1.5){.15}
\pscircle[linewidth=2pt](2.3,.5){.15}
\pscircle[linewidth=2pt](2.3,1.5){.15}
\pscurve[border=1pt]{-*}(.85,1.5)(1.7,1.4)(1.7,.5)
\pscircle(.5,.5){.35}
\psline(.6,1)(.5,1.25)
\psarc[border=1pt](.5,1.5){.35}{180}{320}
\psarc(.5,1.5){.35}{-60}{180}
\psline[border=1pt]{*-}(.6,1)(.5,.75)
\pscurve(.85,.5)(1,.9)(1.3,.4)(1.7,.5)(1.95,.5)
\pscircle(2.3,.5){.4}
\psarc[border=1pt](1.2,.7){.2}{180}{5}
\psarc(1.2,.7){.2}{45}{105}
\psarc(1.4,.7){.1}{-65}{240}
\pscircle(2.3,1.5){.25}
\pscircle(2.3,.5){.25}
\pscircle(.5,.5){.25}
\pscircle(.5,1.5){.25}
\pscircle[linewidth=2pt](.5,.5){.15}
\psline{-*}(2.3,1.25)(2.3,1.1)
\end{pspicture}
\end{center}
\end{lemma}
\eop

As a consequence of Lemma~\ref{lemY4}, we also have the following lemma that provides an inverse for a $Y$-graph.
A mark $\begin{pspicture}[.4](0,0)(.3,.3) 
\psline(.15,0)(.15,.3)
\psline[border=1pt]{-}(.05,.05)(.25,.25)
\end{pspicture}$ on an edge indicates a positive half-twist of this edge. 

\begin{lemma}[Theorem 3.2 in \cite{ggp}]
\label{lemY6}
The surgery on the following $Y$-link is trivial.
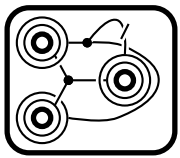
\begin{center}
\begin{pspicture}[.4](-.2,-.2)(2.5,2.2)
\psline(1.35,1)(.85,1)
\pscircle[border=1pt](1.6,1){.35}
\psline(.85,1)(.67,1.33)
\pscircle[border=1pt](.5,1.5){.35}
\pscircle(.5,.5){.35}
\psline[border=1pt]{*-}(.85,1)(.67,.67)
\pscurve{*-}(1.1,1.5)(1.6,1.45)(2.05,1)(1.5,.5)(.85,.5)
\pscurve[border=1pt](1.1,1.5)(1.5,1.8)(1.6,1.35)
\psline[border=1pt](1.55,1.55)(1.65,1.75)
\psframe[linewidth=2pt, framearc=.3](0,0)(2.3,2)
\pscircle[linewidth=2pt](.5,1.5){.15}
\pscircle[linewidth=2pt](1.6,1){.15}
\pscircle(1.6,1){.25}
\pscircle(.5,.5){.25}
\pscircle(.5,1.5){.25}
\psline(.85,1.5)(1.1,1.5)
\pscircle[linewidth=2pt](.5,.5){.15}
\psdot(1.1,1.5)
\end{pspicture}
\end{center}
\end{lemma}

\subsection{A clasper presentation of algebraically split links}
A leaf $\ell$ of a $Y$-link $G$ is a {\em meridional leaf\/} or is a  
{\em meridian\/}  of a link $L$, if it is trivial, and
if it bounds a meridian disk of some link component whose interior intersects  
 $G\cup L$ at exactly one point of $L$.

Say that a $Y$-link $G$ {\em laces\/} the trivial $r$-component link $U^{(r)}$ of a  
connected $3$-manifold if
\begin{itemize}
\item each of the $Y$-link components contains a meridional leaf of $U^{(r)}$,
\item The components $U_i$ of $U^{(r)}$ bound
disjoint disks $(D_i)_{i=1,\dots,r}$ ($U_i = \partial D_i$)
so that $D_i \cap G$ is inside the meridional leaves of $U_i$ (and contains one point per meridional leaf of $U_i$),
\item no component of $G$ contains more than one meridional leaf of a given component
$U_i$.
\end{itemize}

Performing the surgery on such a $G$ transforms $U^{(r)}$ into the
link $(K_1, \dots, K_r)=U^{(r)}_G$ in $M$
that is {\em presented by $(G,U^{(r)})$\/}.

Since any null-homologous knot bounds an oriented Seifert surface, by  
Lemma~\ref{lemysur},
it is easy to see that any null-homologous knot is presented by a pair  
$(G,U_1)$,
where $G$ is a $Y$-link that laces the trivial knot.
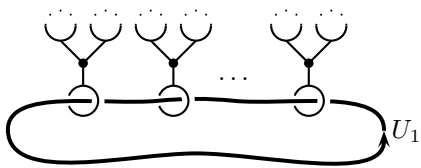
\begin{center}
\begin{pspicture}[.2](0,-.7)(5,1.6) 
\psecurve[linewidth=1.5pt]{->}(2.5,-.5)(5,-.2)(3,.2)(2,.2)(0,-.2)(2.5,-.5)(5,-.2)(2.5,.1)
\psarc[border=2pt](1,.2){.2}{-160}{160}
\psarc(.7,1.2){.2}{180}{0}
\psarc[linestyle=dotted](.7,1.2){.2}{0}{180}
\psarc(1.3,1.2){.2}{180}{0}
\psarc[linestyle=dotted](1.3,1.2){.2}{0}{180}
\psline(.7,1)(1,.7)(1.3,1)
\psline{*-}(1,.7)(1,.4)
\psarc[border=2pt](2.2,.2){.2}{-160}{160}
\psarc(1.9,1.2){.2}{180}{0}
\psarc[linestyle=dotted](1.9,1.2){.2}{0}{180}
\psarc(2.5,1.2){.2}{180}{0}
\psarc[linestyle=dotted](2.5,1.2){.2}{0}{180}
\psline(1.9,1)(2.2,.7)(2.5,1)
\psline{*-}(2.2,.7)(2.2,.4)
\rput[l](5.1,-.2){$U_1$}
\rput(3,.5){$\dots$}
\psarc[border=2pt](4,.2){.2}{-160}{160}
\psarc(3.7,1.2){.2}{180}{0}
\psarc[linestyle=dotted](3.7,1.2){.2}{0}{180}
\psarc(4.3,1.2){.2}{180}{0}
\psarc[linestyle=dotted](4.3,1.2){.2}{0}{180}
\psline(3.7,1)(4,.7)(4.3,1)
\psline{*-}(4,.7)(4,.4)
\end{pspicture}
\end{center}

In a connected oriented compact $3$-manifold $M$ such that $H_2(M;\ZZ)=0$, the {\em linking number of a null-homologous knot $K$
with a knot  $C$ in its complement\/} is well-defined as the algebraic intersection of $C$
with a surface bounded by $K$.
The {\em Milnor triple linking number\/} $\mu(K_1,K_2,K_3)$ of three null-homologous knots $K_1, K_2, K_3$ that do not link each other is also well-defined,
as the algebraic intersection of three Seifert surfaces of these knots in the complement of the other ones
with the sign $\mu(K_1,K_2,K_3)=-\langle \Sigma_1,\Sigma_2,\Sigma_3\rangle$.

Let $G$ be a $Y$-link that laces the trivial link $U^{(r)}$ of $M$.
Let $m_i$ denote the homology class of the oriented meridian of $U_i$.
Say that a component of $G$ is {\em of type $(\varepsilon_i m_i, \varepsilon_j m_j,f)$\/} if its leaves are one meridian of $U_i$, one meridian of $U_j$, and another oriented framed leaf $f$ and if it can be oriented so that the homology classes of its oriented leaves read $\varepsilon_i m_i$, $\varepsilon_j m_j$ and $[f]$ with respect to the cyclic order induced by the orientation, with $\varepsilon_i, \varepsilon_j \in \{-1,1\}$. Similarly, say that a component of $G$ is {\em of type $(\varepsilon_i m_i, \varepsilon_j m_j,\varepsilon_k m_k)$\/} if its leaves are one meridian of $U_i$, one meridian of $U_j$, and one meridian of $U_k$, and if it can be oriented so that the homology classes of its oriented leaves read $\varepsilon_i m_i$, $\varepsilon_j m_j$ and $\varepsilon_k m_k$ with respect to the cyclic order induced by the orientation, with $\varepsilon_i, \varepsilon_j, \varepsilon_k \in \{-1,1\}$.

\begin{lemma}
\label{lemintclas}
Let $G$ be a $Y$-link that laces the trivial link $U^{(r)}$ of an oriented connected $3$--manifold $M$.
Let $L=(K_1, \dots, K_r)=U^{(r)}_G$ be the link presented by $G$.
Then $L$ is algebraically split, and the $K_i$ bound surfaces $\Sigma_i$ such that
\begin{itemize}
\item for any $\{i,j\} \subset \{1,2,\dots,r\}$, $\Sigma_i \cap \Sigma_j$ is the union over all the components {\em of type $(\varepsilon_i m_i, \varepsilon_j m_j,f)$\/} of the framed oriented leaves $\varepsilon_i\varepsilon_j f$,
\item for any $\{i,j,k\} \subset \{1,2,\dots,r\}$, the oriented intersection $\Sigma_i \cap \Sigma_j \cap \Sigma_k$ is 
a union over all the components {\em of type $(\varepsilon_i m_i, \varepsilon_j m_j,\varepsilon_k m_k)$\/} of points
with sign $\varepsilon_i\varepsilon_j \varepsilon_k$.
\end{itemize}
In particular, if $H_2(M;\ZZ)=0$, then $\mu(K_i,K_j,K_k)$ is the sum over all the components {\em of type $(\varepsilon_i m_i, \varepsilon_j m_j,\varepsilon_k m_k)$\/} of the contributions $(-\varepsilon_i\varepsilon_j \varepsilon_k)$.
\end{lemma}
\bp Define the {\em index\/} of a component $Y$ of $G$ as the smallest $i$ such that $Y$ has a meridional leaf of $U_i$.
Realize the surgeries on the components of index $i$ of $G$ by applying Lemma~\ref{lemysur} to the trivial meridional leaf $\ell$ of $U_i$ and to the part of $U_i$ going through $\ell$.
These surgeries transform $U^{(r)}$ into $L$ and allow us to see each $K_i$ as the boundary of a surface $\tilde{\Sigma}_i$ whose $1$-handles are thickenings of the framed leaves that are not meridians of $U_i$ of the components of index $i$.

So far, $\tilde{\Sigma}_k$ may intersect the $K_i$ with $i<k$ (but not the $K_i$ with $i>k$). More precisely, if $i<k$, each component  of index $i$ of type $(m_i,\varepsilon m_k, f)$ or $(m_i,f,-\varepsilon m_k)$ gives rise to an arc of intersection of $\tilde{\Sigma}_i \cap \tilde{\Sigma}_k$.
Tubing $\tilde{\Sigma}_k$ along the part of $K_i$ between the two extremities of the intersection arc that is contained in the surgery picture transforms this arc of intersection into $\varepsilon f$ and removes the intersection of $K_i$ with $\tilde{\Sigma}_k$. 

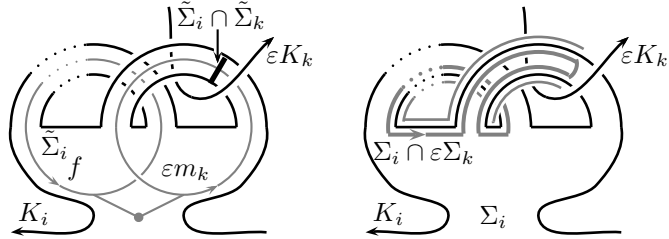
\begin{center}
\begin{pspicture}[.4](-.2,.2)(4.4,3.3)
\psarc[linewidth=1pt](1.1,1.6){1.1}{0}{95}
\psarc[linewidth=1pt,linestyle=dotted](1.1,1.6){1.1}{95}{140}
\psarc[linewidth=1pt](1.1,1.6){1.1}{140}{180}
\psarc[linewidth=1pt](1.1,1.6){.7}{0}{95}
\psarc[linewidth=1pt,linestyle=dotted](1.1,1.6){.7}{95}{140}
\psarc[linewidth=1pt](1.1,1.6){.7}{140}{180}
\psarc[linestyle=dotted,linecolor=gray](1.1,1.6){.9}{95}{140}
\psarc[linecolor=gray]{->}(1.1,1.6){.9}{140}{240}
\psarc[linecolor=gray](1.1,1.6){.9}{240}{95}
\psecurve[linewidth=1pt](2,4)(2.1,3.2)(2.4,2)(3.5,2.8)
\psarc[linewidth=1pt,border=2pt](2.3,1.6){1.1}{0}{180}
\psarc[linewidth=1pt,border=2pt](2.3,1.6){.7}{0}{180}
\psarc[border=1pt,linecolor=gray](2.3,1.6){.9}{-60}{180}
\psarc[border=1pt,linecolor=gray]{->}(2.3,1.6){.9}{180}{300}
\psecurve[linewidth=1pt,border=2pt]{->}(2.1,3.2)(2.4,2)(3.5,2.8)(4.2,3.5)
\psecurve[linewidth=1pt]{->}(1.1,2.7)(0,1.6)(.33,.83)(1.1,.4)(0,.2)(-1,.2)
\psecurve[linewidth=1pt](2.3,2.7)(3.4,1.6)(3.07,.83)(2.3,.4)(3.4,.2)(4.4,.2)
\rput[rb](2.65,.9){$\varepsilon m_k$} 
\rput[lb](.75,.9){$f$} 
\rput[lb](0.1,.3){$K_i$}
\rput[lt](3.4,2.7){$\varepsilon K_k$}
\psline[linewidth=2pt](2.65,2.18)(2.85,2.52)
\rput[b](2.8,2.9){$\tilde{\Sigma}_i \cap \tilde{\Sigma}_k$}
\rput[tl](.4,1.55){$\tilde{\Sigma}_i$}
\psline[border=1pt]{->}(2.75,2.9)(2.75,2.45)
\psline[linecolor=gray]{*-}(1.7,.4)(1.1,.7)
\psline[linecolor=gray](1.7,.4)(2.3,.7)
\psline[linewidth=1pt](.4,1.6)(1.2,1.6)
\psline[linewidth=1pt](1.6,1.6)(1.8,1.6)
\psline[linewidth=1pt](2.2,1.6)(3,1.6)
\end{pspicture}
\begin{pspicture}[.4](-.2,.2)(4.4,3.3)
\psarc[linewidth=1pt](1.1,1.6){1.1}{0}{95}
\psarc[linewidth=1pt,linestyle=dotted](1.1,1.6){1.1}{95}{140}
\psarc[linewidth=1pt](1.1,1.6){1.1}{140}{180}
\psarc[linewidth=1pt](1.1,1.6){.7}{0}{95}
\psarc[linewidth=1pt,linestyle=dotted](1.1,1.6){.7}{95}{140}
\psarc[linewidth=1pt](1.1,1.6){.7}{140}{180}
\psarc[linewidth=1.5pt,linecolor=gray](1.1,1.6){.8}{-6}{95}
\psarc[linewidth=1.5pt,linecolor=gray,linestyle=dotted](1.1,1.6){.8}{95}{140}
\psarc[linewidth=1.5pt,linecolor=gray](1.1,1.6){.8}{140}{186}
\psarc[linewidth=1pt,linecolor=gray](1.1,1.6){.6}{5}{95}
\psarc[linewidth=1pt,linecolor=gray,linestyle=dotted](1.1,1.6){.6}{95}{140}
\psarc[linewidth=1pt,linecolor=gray](1.1,1.6){.6}{140}{173}
\psecurve[linewidth=1pt](2,4)(2.1,3.2)(2.4,2)(3.5,2.8)
\psarc[linewidth=1.5pt,linecolor=gray,border=1pt](2.3,1.6){1}{60}{185}
\psarc[linewidth=1.5pt,linecolor=gray,border=1pt](2.3,1.6){.8}{60}{186}
\psarc[linewidth=1pt,linecolor=gray,border=1pt](2.3,1.6){.6}{60}{173}
\psarc[linewidth=1pt,linecolor=gray,border=1pt](2.3,1.6){1.2}{60}{176}
\psarc[linewidth=1pt,border=1pt](2.3,1.6){1.1}{0}{180}
\psarc[linewidth=1pt,border=1pt](2.3,1.6){.7}{0}{180}
\psecurve[linewidth=1pt,border=2pt]{->}(2.1,3.2)(2.4,2)(3.5,2.8)(4.2,3.5)
\psecurve[linewidth=1pt]{->}(1.1,2.7)(0,1.6)(.33,.83)(1.1,.4)(0,.2)(-1,.2)
\psecurve[linewidth=1pt](2.3,2.7)(3.4,1.6)(3.07,.83)(2.3,.4)(3.4,.2)(4.4,.2) 
\rput[lb](0.1,.3){$K_i$}
\rput[lt](3.4,2.7){$\varepsilon K_k$}
\psecurve[linewidth=1.5pt,linecolor=gray](2.3,2.46)(2.7,2.26)(2.8,2.5)(2.4,2.64)
\psline[linewidth=1pt](.4,1.6)(1.2,1.6)
\psline[linewidth=1pt](1.6,1.6)(1.8,1.6)
\psline[linewidth=1pt](2.2,1.6)(3,1.6)
\psline[linewidth=1pt,linecolor=gray](.5,1.7)(1.1,1.7)
\psline[linewidth=1.5pt,linecolor=gray](1.5,1.5)(1.9,1.5)
\rput(1.7,.4){${\Sigma}_i$}
\psline[linewidth=1.5pt,linecolor=gray]{->}(.3,1.5)(.8,1.5)
\rput[lt](.1,1.4){${\Sigma}_i \cap \varepsilon {\Sigma}_k$}
\psline[linewidth=1.5pt,linecolor=gray](.8,1.5)(1.3,1.5)
\end{pspicture}
\end{center}

If the third leaf is a meridian of 
$K_j$ for $i<j<k$ then perform the tubing along this leaf inside the tubing of $\tilde{\Sigma}_j$ along the meridional
leaf of $m_k$.
Let $\Sigma_k$ denote the surface obtained after all these tubings.
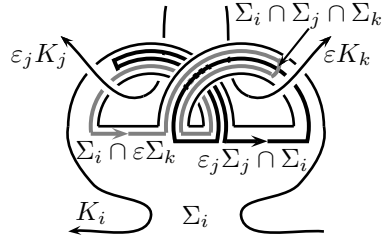
\begin{center}
\begin{pspicture}[.4](-.2,.2)(4.4,3.3)
\psecurve[linewidth=1pt](1.4,4)(1.3,3.2)(1,2)(-.1,2.8)
\psarc[linewidth=1pt,border=2pt](1.1,1.6){1.1}{0}{180}
\psarc[linewidth=1pt,border=2pt](1.1,1.6){.7}{0}{180}
\psarc[linewidth=1.5pt,border=1pt](1.1,1.6){.9}{-14}{120}
\psarc[linewidth=1.5pt](1.1,1.6){1}{-13}{120}
\psline[linewidth=1.5pt](.7,2.26)(.6,2.5)
\psecurve[linewidth=1pt](2,4)(2.1,3.2)(2.4,2)(3.5,2.8)
\psarc[linewidth=1.5pt,linecolor=gray,border=1pt](1.1,1.6){.8}{-6}{186}
\psecurve[linewidth=1pt,border=2pt]{->}(1.3,3.2)(1,2)(-.1,2.8)(-.8,3.5)
\psarc[linewidth=1.5pt,linecolor=gray,border=1pt](2.3,1.6){1}{60}{185}
\psarc[linewidth=1.5pt,linecolor=gray,border=1pt](2.3,1.6){.8}{60}{186}
\psarc[linewidth=1pt,border=1pt](2.3,1.6){1.1}{0}{180}
\psarc[linewidth=1.5pt](2.3,1.6){.9}{-14}{194}
\psarc[linewidth=1pt,border=1pt](2.3,1.6){.7}{0}{180}
\psecurve[linewidth=1pt,border=2pt]{->}(2.1,3.2)(2.4,2)(3.5,2.8)(4.2,3.5)
\psecurve[linewidth=1pt]{->}(1.1,2.7)(0,1.6)(.33,.83)(1.1,.4)(0,.2)(-1,.2)
\psecurve[linewidth=1pt](2.3,2.7)(3.4,1.6)(3.07,.83)(2.3,.4)(3.4,.2)(4.4,.2) 
\rput[lb](0.1,.3){$K_i$}
\rput[lt](3.4,2.7){$\varepsilon K_k$}
\rput[rt](0,2.7){$\varepsilon_j K_j$}
\psecurve[linewidth=1.5pt,linecolor=gray](2.3,2.46)(2.7,2.26)(2.8,2.5)(2.4,2.64)
\rput[b](3.2,2.9){${\Sigma}_i \cap {\Sigma}_j \cap {\Sigma}_k$}
\rput(1.7,.4){${\Sigma}_i$}
\psline[border=1pt]{->}(3.2,2.9)(2.8,2.4)
\psline[linewidth=1pt](.4,1.6)(1.2,1.6)
\psline[linewidth=1pt](1.6,1.6)(1.8,1.6)
\psline[linewidth=1pt](2.2,1.6)(3,1.6)
\psline[linewidth=1.5pt](1.4,1.4)(2,1.4)
\psline[linewidth=1.5pt]{->}(2.1,1.4)(2.65,1.4)
\rput[rt](3.2,1.3){$ \varepsilon_j{\Sigma}_j \cap \Sigma_i$}
\psline[linewidth=1.5pt](2.65,1.4)(3.2,1.4)
\psline[linewidth=1.5pt,linecolor=gray]{->}(.3,1.5)(.8,1.5)
\rput[lt](.1,1.4){${\Sigma}_i \cap \varepsilon {\Sigma}_k$}
\psline[linewidth=1.5pt,linecolor=gray](.8,1.5)(1.3,1.5)
\psline[linewidth=1.5pt,linecolor=gray](1.5,1.5)(1.9,1.5)
\end{pspicture}
\end{center}
It is left to the reader to check that the surfaces have the announced properties.
\eop

Say that a $Y$-link $G$ {\em $\mu$-laces\/} the trivial $r$-component link $U^{(r)}$ of a  
connected $3$-manifold if it laces it, and if for any triple $\{i,j,k\}$ of integers in $\{1,\dots,r\}$, there are
exactly $|\mu(K_i,K_j,K_k)|$ components with one leaf that
links $U_i$, one leaf that links $U_j$ and one leaf that links $U_k$.

It is known that any algebraically split link can be presented by a $Y$-link $G$  
that laces the trivial link $U^{(r)}$
\cite[Lemma 5.6]{ggp}, \cite{matv,murnak}.
We prove the following proposition that refines this result, and that, together with Lemma~\ref{lemintclas}, proves Proposition~\ref{proprealmil}. 

\begin{proposition}
\label{propclaspmu}
For any algebraically split link $L=(K_1,\dots,K_r)$ in a connected
$3$-manifold $M$ such that $H_2(M;\ZZ)=0$, there exists a $Y$-link $G$  
that $\mu$-laces the trivial link $U^{(r)}$ of $M$
such that $(G,U^{(r)})$ presents $L$.
\end{proposition}

This proposition will be a direct corollary from the slightly more general proposition~\ref{propclaspmubased} below, that may be used for the study of homology handlebodies.

Here, an {\em $r$-component based link\/} is an embedding $\Gamma_L$
of the following graph with $r$ loops, up to isotopy. Its underlying link is the restriction of the embedding to its $r$ loops. 
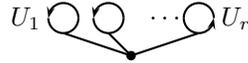
\begin{center} 
\begin{pspicture}[.2](0,.6)(3.2,1.5) 
\rput[r](.4,1.2){$U_1$}
\psarc(.7,1.2){.2}{180}{0}
\psarc{->}(.7,1.2){.2}{0}{180}
\psarc(1.3,1.2){.2}{180}{0}
\psarc{->}(1.3,1.2){.2}{0}{180}
\psline(.7,1)(1.6,.7)(1.3,1)
\psline{*-}(1.6,.7)(2.5,1)
\rput(2.05,1.2){$\dots$}
\psarc{->}(2.5,1.2){.2}{180}{0}
\psarc(2.5,1.2){.2}{0}{180}
\rput[l](2.8,1.2){$U_r$}
\end{pspicture}
\end{center}

The {\em trivial $r$-component based link\/} $\Gamma_U^{(r)}$ is the above $r$-component based link seen in a $3$--ball.

\begin{proposition}
\label{propclaspmubased}
For any based $r$-component link $\Gamma_L$, whose underlying link $L=(K_1,\dots,K_r)$ is algebraically split, in a connected
$3$-manifold $M$ such that $H_2(M;\ZZ)=0$, there exists a $Y$-link $G$  
in $M \setminus \Gamma_U^{(r)}$ that $\mu$-laces the trivial link $U^{(r)}$ of $M$
such that $(G,\Gamma_U^{(r)})$ presents $\Gamma_L$.
\end{proposition}
\bp
For any sublink $L^{\prime}$ of $L$, there is a canonical subgraph $\Gamma_{L^{\prime}}$ of $\Gamma_L$ that is a based link with underlying link 
$L^{\prime}$.
We prove Proposition~\ref{propclaspmubased} by proving the following statement by
induction on the number $r$ of components of $L$.

\noindent {\em Induction hypothesis}\\
Let $M$ be a connected $3$-manifold such that $H_2(M;\ZZ)=0$.
Let $\Gamma_{L \cup L^{\prime}}$ be a based algebraically split link in $M$ where $L$ has  
$r$ components. 
Let $\Gamma_{U^{(r)} \cup L^{\prime}}$ be the based link obtained from $\Gamma_{L \cup L^{\prime}}$ by replacing $\Gamma_L$ by $\Gamma_U^{(r)}$ so that each component of $U^{(r)}$ bounds a disk $D_i$ whose interior does not meet $\Gamma_{U^{(r)} \cup L^{\prime}}$.
Then there exists a $Y$-link $G$ in of $M \setminus \Gamma_{L^{\prime}}$ such that the following set of properties
$H(G,\Gamma_L,\Gamma_{L^{\prime}})$ is satisfied.

\begin{itemize}
\item $G \subset M \setminus \Gamma_{U^{(r)} \cup L^{\prime}}$,
\item $G$ {\em $\mu$-laces\/} the trivial  
link $U^{(r)}$ of $M \setminus L^{\prime}$
\item
$(G,\Gamma_{U^{(r)} \cup L^{\prime}})$ presents $\Gamma_{L \cup L^{\prime}}$ in $M$,
\item the only leaves of $G$ that
link $L^{\prime}$ algebraically
are meridional leaves of  
$L^{\prime}$,
\item no component of $G$ contains more than one meridional leaf of a
given component
of $L^{\prime}$,
\item For any triple $\{I,J,K\}$ of components of $L \cup L^{\prime}$  
with at least one component
in $L$, there are
exactly $|\mu(I,J,K)|$ components of $G$ with one leaf that
links $I$, one leaf that links $J$ and one leaf that links $K$.
\end{itemize}

This statement is obviously true for $0$-component links.

Assume that it is true for $(r-1)$-component links,
we wish to prove it for $(L=(K_1,\dots,K_r), L^{\prime})$.
Let $U^{(r-1)}=(U_1, \dots, U_{r-1})$ denote the trivial $(r-1)$-component link that bounds a disjoint union of disks
$(D_i)_{i=1,\dots, r-1}$.
By induction, there exists $(G_1 \subset M  
\setminus \Gamma_{L^{\prime} \cup K_r})$ such that $H(G_1,\Gamma_{K_1,\dots,K_{r-1}}, \Gamma_{K_r \cup L^{\prime}})$ is satisfied.

Consider a two-dimensional disk $D$ that meets $K_r$ along an arc $\alpha$ of its boundary around which all the meridional leaves of $K_r$ are, and such that $D$ intersects all the meridional leaves, so that
$$K^{\prime}_r = (K_r \setminus \stackrel{\circ}{\alpha}) \cup (-\partial D \setminus \alpha)$$ bounds a surface $\Sigma$ that meets neither  $\Gamma_{L^{\prime} \cup U^{r-1}} \cup \cup_{i<r}D_i$, nor the path $\gamma_{r}$ from the  vertex of $\Gamma_{L \cup L^{\prime}}$ to $K_r$, nor the leaves of $G_1$.

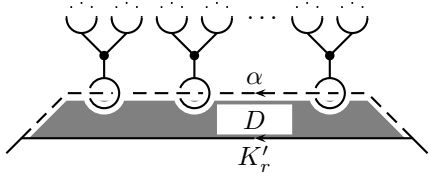
\begin{center}
\begin{pspicture}[.2](0,-1.1)(5,1.6) 
\pspolygon*[linecolor=gray,framearc=.3](0,-.4)(.5,.1)(4.5,.1)(5,-.4)
\psframe*[linecolor=white](2.5,-.35)(3.5,.05)
\rput(3,-.15){$D$}
\psline[linestyle=dashed,framearc=.3](-.1,-.4)(.5,.2)(3,.2)
\psline[linestyle=dashed,framearc=.3]{->}(5.1,-.4)(4.5,.2)(3,.2)
\psline{->}(5.3,-.6)(5.1,-.4)(3,-.4)
\psline(-.3,-.6)(-.1,-.4)(3,-.4)
\rput[t](3,-.5){$K^{\prime}_r$}
\psarc[border=2pt](1,.2){.2}{-160}{160}
\psarc(.7,1.2){.2}{180}{0}
\psarc[linestyle=dotted](.7,1.2){.2}{0}{180}
\psarc(1.3,1.2){.2}{180}{0}
\psarc[linestyle=dotted](1.3,1.2){.2}{0}{180}
\psline(.7,1)(1,.7)(1.3,1)
\psline{*-}(1,.7)(1,.4)
\psarc[border=2pt](2.2,.2){.2}{-160}{160}
\psarc(1.9,1.2){.2}{180}{0}
\psarc[linestyle=dotted](1.9,1.2){.2}{0}{180}
\psarc(2.5,1.2){.2}{180}{0}
\psarc[linestyle=dotted](2.5,1.2){.2}{0}{180}
\psline(1.9,1)(2.2,.7)(2.5,1)
\psline{*-}(2.2,.7)(2.2,.4)
\rput(3.1,1.2){$\dots$}
\rput(3,.4){$\alpha$}
\psarc[border=2pt](4,.2){.2}{-160}{160}
\psarc(3.7,1.2){.2}{180}{0}
\psarc[linestyle=dotted](3.7,1.2){.2}{0}{180}
\psarc(4.3,1.2){.2}{180}{0}
\psarc[linestyle=dotted](4.3,1.2){.2}{0}{180}
\psline(3.7,1)(4,.7)(4.3,1)
\psline{*-}(4,.7)(4,.4)
\end{pspicture}
\end{center}

\begin{lemma}
The graph $G_1$ and the surface $\Sigma$ can be modified so that $\Sigma$ does not meet $G_1$ at all outside the meridional leaves of $L$,
and the following set of assumptions\\ $H_2(G_1,\Gamma_{K_1,\dots,K_{r-1}}, \Gamma_{K_r \cup L^{\prime}},\Sigma)$

\begin{itemize}
\item $\Sigma$ meets neither  $\Gamma_{U^{(r-1)} \cup L^{\prime}} \cup \cup_{i<r}D_i$, nor $\gamma_{r}$, nor the leaves of $G_1$
\item $H(G_1,\Gamma_{K_1,\dots,K_{r-1}}, \Gamma_{K_r \cup L^{\prime}})$ is satisfied except 
 that
components of $G_1$ are allowed to have no meridians of $L\setminus K_r$ provided that they have a 
meridian of $K_r$.
\end{itemize}is still satisfied.
\end{lemma}

\bp We need to remove the intersections of $\Sigma$ with the edges of $G_1$.
By isotopy, without loss, assume that no edge adjacent to
a meridional leaf of $K_r$ intersects $\Sigma_r$ (push the intersection on the two other edges, if necessary). 
Similarly, assume that if a component contains only one meridian of $L$, the edge adjacent to this
component does not meet $\Sigma_r$. Now, the intersections of the edges adjacent to
non-meridional leaves can be removed by tubing $\Sigma_r$ around the part of the $Y$-graph that contains
the corresponding leaf. Here {\em tubing\/} means replacing a small disk of $\Sigma_r$ in a neighborhood of an intersection point with an edge by the closure of its complement in the boundary of a regular neighborhood of
the part of the $Y$-graph that contains
the corresponding leaf, as below. 

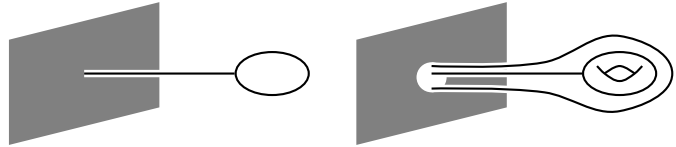
\begin{center}
\begin{pspicture}[.2](-.5,0)(4,2) 
\pspolygon*[linecolor=gray](0,0)(2,.5)(2,1.9)(0,1.4)
\psline[border=1pt](1,.95)(3,.95)
\psellipse(3.5,.95)(.5,.3)
\end{pspicture}
\begin{pspicture}[.2](-.5,0)(4.5,2) 
\pspolygon*[linecolor=gray](0,0)(2,.5)(2,1.9)(0,1.4)
\pscircle*[linecolor=white](1,.9){.2}
\psline[border=1pt](1,.95)(3,.95)
\psellipse(3.5,.95)(.5,.3)
\pscurve[border=1pt](1,.75)(2.5,.7)(3.5,.45)(4.2,.95)(3.5,1.45)(2.5,1.1)(1,1.05)
\pscurve(3.3,.95)(3.5,1.05)(3.7,.95)
\pscurve(3.2,1.05)(3.3,.95)(3.5,.85)(3.7,.95)(3.8,1.05)
\end{pspicture}
\end{center}

Thus, we are under the assumptions $H(E)$ that the only edge intersections occur on edges adjacent to a meridional leaf of some $K_j$,
with $j<r$, of components with at least two meridional leaves of $L$.
Define the complexity $c(\Sigma; G_1)$ of such a situation as follows.
Define the complexity $c_e(Y)$ of a component $Y$ of $G_1$ as the number of intersection points of its edges with $\Sigma$.
Define the complexity $c(\Sigma; G_1)$ as the pair
(maximal complexity $c_e$ of the components, number of components with this
complexity)
ordered by the lexicographic order. 

Now, to prove the lemma, it is enough to prove that there exists a pair $(\Sigma; G_1)$
with lower complexity such that
$H_2(G_1,\Gamma_{K_1,\dots,K_{r-1}}, \Gamma_{K_r \cup L^{\prime}},\Sigma)$ and $H(E)$ are satisfied.

Consider a component $Y$ of $G_1$ with maximal complexity, and its edge $e$ with the maximal number of intersection points.
By hypothesis, $e$ is adjacent with a meridional leaf $\ell$ of some component $K_i$ with $i<r$.
Remove the point of $e \cap \Sigma$ that is closest to $\ell$ as follows.
By our assumptions, $\Sigma$ intersects a neighborhood of $Y$ in the gray part of the following 
picture, where the intersection point that will be removed is at the top right corner. Perform the modification of Lemma~\ref{lemY5} so that the resulting three graphs are like in the following picture with respect to the positions of the possible intersections
with $\Sigma$. 

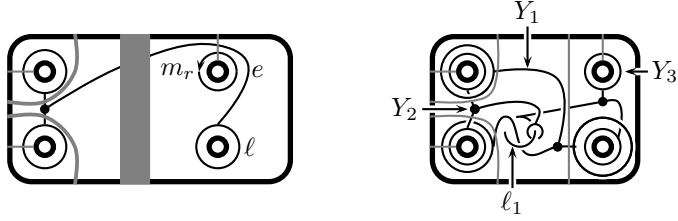
\begin{center}
\begin{pspicture}[.4](-.1,-.6)(4,2.4)
\psframe[linewidth=2pt, framearc=.3](0,0)(3.8,2)
\pscircle[linewidth=2pt](.5,.5){.15}
\pscircle[linewidth=2pt](.5,1.5){.15}
\pscircle[linewidth=2pt](2.8,.5){.15}
\pscircle[linewidth=2pt](2.8,1.5){.15}
\psarc{->}(2.8,1.5){.25}{0}{180}
\psarc(2.8,1.5){.25}{180}{0}
\pscircle(.5,.5){.3}
\pscircle(.5,1.5){.3}
\pscircle(2.8,.5){.3}
\psline(.5,.8)(.5,1.2)
\pscurve{*-}(.5,1)(1.6,1.6)(2.8,1.85)(3.2,1.6)(2.8,.8)
\psline[linecolor=gray](2.8,1.65)(2.8,2)
\psline[linecolor=gray](0,.5)(.35,.5)
\psline[linecolor=gray](0,1.5)(.35,1.5)
\psline[linecolor=gray, linewidth=.4](1.7,0)(1.7,2)
\pscurve[linecolor=gray, linewidth=.05](0,1.1)(.5,1.1)(.9,1.5)(.9,2)
\pscurve[linecolor=gray, linewidth=.05](0,.9)(.5,.9)(.9,.5)(.9,0)
\rput[r](2.5,1.5){$m_r$}
\rput[l](3.25,1.5){$e$}
\rput[l](3.15,.5){$\ell$}
\end{pspicture} 
\begin{pspicture}[.4](-1.5,-.6)(4,2.4)
\psframe[linewidth=2pt, framearc=.3](0,0)(2.8,2)
\pscurve(2.3,1.1)(2.5,1.1)(2.5,.6)
\psline(2.3,1.1)(1.2,.9)
\pscurve[border=1pt](.6,1)(1.45,1)(1.4,.8)
\pscircle[border=1pt](2.3,.5){.4}
\pscircle[linewidth=2pt](.5,1.5){.15}
\pscircle[linewidth=2pt](2.3,.5){.15}
\pscircle[linewidth=2pt](2.3,1.5){.15}
\pscurve[border=1pt]{-*}(.85,1.5)(1.7,1.4)(1.7,.5)
\pscircle(.5,.5){.35}
\psline(.6,1)(.5,1.25)
\psarc[border=1pt](.5,1.5){.35}{180}{320}
\psarc(.5,1.5){.35}{-60}{180}
\psline[border=1pt]{*-}(.6,1)(.5,.75)
\pscurve(.85,.5)(1,.9)(1.3,.4)(1.7,.5)(1.95,.5)
\pscircle(2.3,.5){.4}
\psarc[border=1pt](1.2,.7){.2}{180}{5}
\psarc(1.2,.7){.2}{45}{105}
\psarc(1.4,.7){.1}{-65}{240}
\pscircle(2.3,1.5){.25}
\pscircle(2.3,.5){.25}
\pscircle(.5,.5){.25}
\pscircle(.5,1.5){.25}
\pscircle[linewidth=2pt](.5,.5){.15}
\psline{-*}(2.3,1.25)(2.3,1.1)
\psline[linecolor=gray](2.3,1.65)(2.3,2)
\psline[linecolor=gray](0,.5)(.35,.5)
\psline[linecolor=gray](0,1.5)(.35,1.5)
\psline[linecolor=gray](1.85,0)(1.85,2)
\pscurve[linecolor=gray](0,1.1)(.5,1.1)(.8,1.2)(.9,1.5)(.9,2)
\pscurve[linecolor=gray](0,.9)(.5,.9)(.8,.8)(.9,.5)(.9,0)
\psline[border=1pt]{->}(1.3,2.1)(1.3,1.55)
\rput[b](1.3,2.15){$Y_1$}
\psline[border=1pt]{->}(-.1,1)(.5,1)
\rput[r](-.15,1){$Y_2$}
\psline[border=1pt]{->}(2.9,1.5)(2.6,1.5)
\rput[l](2.95,1.5){$Y_3$}
\psline[border=1pt]{->}(1.1,-.05)(1.1,.5)
\rput[t](1.1,-.1){$\ell_1$}
\end{pspicture}
\end{center}
Let $Y_1$ be the graph that replaces $Y$ with one edge intersection removed.
Let $Y_3$ be the graph with a meridional leaf of $K_r$, a leaf parallel to $\ell$, and another trivial leaf
$\ell_1$, and let $Y_2$ be the other one with one meridian of $\ell_1$.
Remove all the intersection of $Y_3$ with $\Sigma$ outside its meridional leaf of $K_r$ by tubing. If $Y_2$ has only one meridional leaf of 
$L$, then remove its intersections as before, too. Otherwise, don't change it, it has two meridional leaves, and its
complexity $c_e$ is lower than $c_e(Y)$. 
Slide the meridional leaf of $K_r$ in $Y_3$ so that it is around the arc $\alpha$ of $K_r$.
Thus, the obtained graph and the modified $\Sigma$ together
satisfy 
$H_2(G_1,\Gamma_{K_1,\dots,K_{r-1}}, \Gamma_{K_r \cup L^{\prime}},\Sigma)$ and $H(E)$, and have lower complexity.
The lemma is proved.
\eop

By Lemma~\ref{lemysur}, $K_r \setminus \alpha$ is obtained from $\partial D \setminus \alpha$ by surgery on a $Y$--link $G_2$ in the neighborhood of 
$\Sigma \setminus D$ such that any component of $G_2$ contains exactly one meridian of $\partial D$. Let $U_r = \partial D$. Thus, $K_r$ is obtained from $U_r$ by surgery on $G_1 \cup G_2$, 
$G_1 \cup G_2$ $\mu$-laces the trivial  
link $U^{(r)}$ of $M \setminus L^{\prime}$,
$(G_1 \cup G_2,\Gamma_{U^{(r)} \cup L^{\prime}})$ presents $\Gamma_{L \cup L^{\prime}}$ in $M$.
Let us now modify $G=G_1 \cup G_2$ so that the last three conditions 
of $H(G,\Gamma_L,\Gamma_{L^{\prime}})$ are satisfied in addition to the previous ones.

\noindent $\bullet$ {\em Cutting the leaves so that the only leaves that 
link $L^{\prime}$ algebraically
are meridional leaves of  
$L^{\prime}$}

Use Move $Y4$ of \cite{ggp} (Lemma~\ref{lemY4}) to cut the leaves of $G_2$ that are not meridians of $K_r$ so that they are either $0$-framed meridians of $L^{\prime}$ or they do not link $L^{\prime}$ at all. 
Indeed, this move allows us to cut the leaves into leaves that are
homologically trivial in the complement of $L^{\prime}$, and meridians of the
components of $L^{\prime}$ without
creating further intersections of $G_2$ with the disk $D$.
Define the complexity of a leaf as the minimal number of leaves in such a
decomposition minus one.
Define the complexity of a $Y$-graph as the sum of the complexities of its
leaves.
Finally define the complexity of a $Y$-link as the pair
(maximal complexity of the components, number of components with this
complexity)
ordered by the lexicographic order. The leaves can be cut in order to make
this complexity decrease without
creating further intersections of $G$ with $D$.

\noindent $\bullet$ {\em Sliding the handles so that no component of $G$ contains more than one meridional leaf of a
given component
of $L^{\prime}$.}

Now, we wish to remove the $Y$-components with a meridional leaf of $K_r$ and two meridional leaves of the same component
$J$ of $L^{\prime}$.
By Lemma~\ref{lemysur}, a surgery with respect to such a graph $G_3$ corresponds to a band sum with the boundary of a genus one Seifert
surface as below.

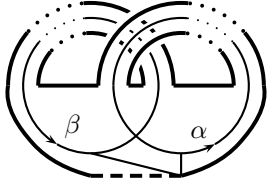
\begin{center}
\begin{pspicture}[.4](-.2,.2)(3.6,2.8)
\psarc[linewidth=1.5pt](1.1,1.6){1.1}{0}{95}
\psarc[linewidth=1.5pt,linestyle=dotted](1.1,1.6){1.1}{95}{140}
\psarc[linewidth=1.5pt](1.1,1.6){1.1}{140}{180}
\psarc[linewidth=1.5pt](1.1,1.6){.7}{0}{95}
\psarc[linewidth=1.5pt,linestyle=dotted](1.1,1.6){.7}{95}{140}
\psarc[linewidth=1.5pt](1.1,1.6){.7}{140}{180}
\psarc[linestyle=dotted](1.1,1.6){.9}{95}{140}
\psarc{->}(1.1,1.6){.9}{140}{240}
\psarc(1.1,1.6){.9}{240}{95}
\psarc[linewidth=1.5pt,border=2pt](2.3,1.6){1.1}{0}{40}
\psarc[linewidth=1.5pt,border=2pt,linestyle=dotted](2.3,1.6){1.1}{40}{85}
\psarc[linewidth=1.5pt,border=2pt](2.3,1.6){1.1}{85}{180}
\psarc[linewidth=1.5pt,border=2pt](2.3,1.6){.7}{0}{40}
\psarc[linewidth=1.5pt,border=2pt,linestyle=dotted](2.3,1.6){.7}{40}{85}
\psarc[linewidth=1.5pt,border=2pt](2.3,1.6){.7}{85}{180}
\psarc[border=1pt](2.3,1.6){.9}{85}{180}
\psarc{->}(2.3,1.6){.9}{180}{300}
\psarc(2.3,1.6){.9}{-60}{40}
\psarc[border=1pt,linestyle=dotted](2.3,1.6){.9}{40}{85}
\psecurve[linewidth=1.5pt](1.1,2.7)(0,1.6)(.33,.83)(1.1,.4)(1.3,.4)
\psecurve[linewidth=1.5pt](2.3,2.7)(3.4,1.6)(3.07,.83)(2.3,.4)(2.1,.4)
\psline[linewidth=1.5pt,linestyle=dashed](1.1,.4)(2.3,.4)
\rput[rb](2.65,.9){$\alpha$} 
\rput[lb](.75,.9){$\beta$} 
\psline(2.3,.4)(1.1,.7)
\psline(2.3,.4)(2.3,.7)
\psline[linewidth=1.5pt](.4,1.6)(1.2,1.6)
\psline[linewidth=1.5pt](1.6,1.6)(1.8,1.6)
\psline[linewidth=1.5pt](2.2,1.6)(3,1.6)
\end{pspicture}
\end{center}
where $\alpha$ and $\beta$ are meridians of $J$.
In this figure, a right-handed (resp.left-handed) Dehn twist of the surface along the simple curve $c(\alpha)$ freely homotopic to $\alpha$
transforms $\beta$ into $\alpha \beta$ (resp. $\alpha^{-1} \beta$) and does not change $\alpha$.
Therefore, the $Y$-graph $G_3$ is equivalent to a $Y$-graph whose leaves are a meridian of $K_r$, the meridian $\alpha$,
and the curve among $\alpha \beta$ and $\alpha^{-1} \beta$ that is null-homologous.

\noindent $\bullet$ {\em Realizing the algebraic cancellations to the Milnor invariants $\mu(K_r,K_s,K_t)$ where $K_s$ and $K_t$ are components of $L^{\prime}$.}

First recall from Lemma~\ref{lemintclas} that $\mu(K_r,K_s,K_t)$ is the sum of the contributions $\varepsilon  \eta$ of the $Y$-graphs of type 
$(m_r,-\varepsilon m_s, \eta m_t)$ or $(m_r,\varepsilon m_t, \eta m_s)$ where $\varepsilon$ and  $\eta$ belong to $\{-1,1\}$.
Second, exchange the order of the $Y$-graphs that link $U_r$ so that all the graphs that contribute with a sign opposite
to the Milnor invariant are followed by a graph that contributes with the sign of the Milnor invariant.
In order to exchange two $Y$-graphs that link $U_r$, perform the following sequence of operations. 
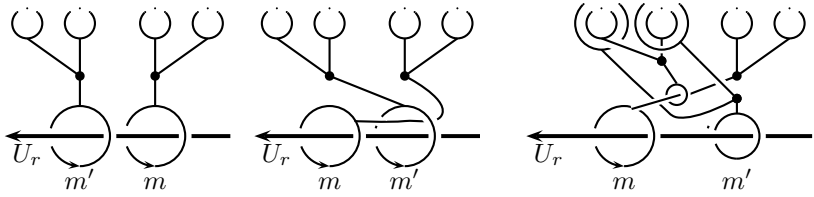
\begin{center}
\begin{pspicture}[.4](0,-.4)(3.2,2.4) 
\psline[linewidth=1.5pt]{<-}(0,.5)(3,.5)
\psarc[border=2pt](1,.5){.4}{-90}{165}
\psarc{->}(1,.5){.4}{-160}{-90}
\psarc[border=2pt](2,.5){.4}{-90}{160}
\psarc{->}(2,.5){.4}{-160}{-90}
\psline(1,.9)(1,1.8)
\psarc(1,2){.2}{135}{45}
\psarc[linestyle=dotted](1,2){.2}{45}{135}
\psarc(.3,2){.2}{135}{45}
\psarc[linestyle=dotted](.3,2){.2}{45}{135}
\psline(2,.9)(2,1.8)
\psarc(2,2){.2}{135}{45}
\psarc[linestyle=dotted](2,2){.2}{45}{135}
\psarc(2.7,2){.2}{135}{45}
\psarc[linestyle=dotted](2.7,2){.2}{45}{135}
\rput[lt](.1,.4){$U_r$}
\psline{*-}(2,1.3)(2.7,1.8)
\psline{*-}(1,1.3)(.3,1.8)
\rput[b](2,-.2){$m$}
\rput[b](1,-.2){$m^{\prime}$}
\end{pspicture}
\begin{pspicture}[.4](0,-.4)(3.2,2.4) 
\psline[linewidth=1.5pt]{<-}(0,.5)(3,.5)
\psarc[border=2pt](1,.5){.4}{-90}{165}
\psarc{->}(1,.5){.4}{-160}{-90}
\pscurve(2,1.3)(2.5,.8)(2,.7)(1.33,.7)
\psarc[border=2pt](2,.5){.4}{-90}{140}
\psarc{->}(2,.5){.4}{-160}{-90}
\psarc(2,.5){.4}{160}{165}
\psline(2,.9)(1,1.3)(1,1.8)
\psarc(1,2){.2}{135}{45}
\psarc[linestyle=dotted](1,2){.2}{45}{135}
\psarc(.3,2){.2}{135}{45}
\psarc[linestyle=dotted](.3,2){.2}{45}{135}
\psline(2,1.3)(2,1.8)
\psarc(2,2){.2}{135}{45}
\psarc[linestyle=dotted](2,2){.2}{45}{135}
\psarc(2.7,2){.2}{135}{45}
\psarc[linestyle=dotted](2.7,2){.2}{45}{135}
\rput[lt](.1,.4){$U_r$}
\psline{*-}(2,1.3)(2.7,1.8)
\psline{*-}(1,1.3)(.3,1.8)
\rput[b](1,-.2){$m$}
\rput[b](2,-.2){$m^{\prime}$}
\end{pspicture}
\begin{pspicture}[.4](-1.1,-.4)(3.2,2.4) 
\psarc(.2,2){.35}{135}{45}
\psline[border=2pt]{*-}(1,1.5)(.2,1.8)
\psline(1,1.5)(1,1.8)
\psarc[border=1pt](1,2){.35}{135}{45}
\psline[linewidth=1.5pt]{<-}(-.8,.5)(3,.5)
\psarc[border=2pt](.5,.5){.4}{-90}{165}
\psarc{->}(.5,.5){.4}{-160}{-90}
\psline(2,1.3)(1.35,1.05)
\pscircle[border=1pt](1.2,1.05){.15}
\pscurve(.2,1.65)(1,.9)(1.2,.75)(2,1)
\psline[border=1pt](1.25,1.05)(.6,.86)
\psarc[border=2pt](2,.5){.3}{-90}{165}
\psarc[border=2pt](2,.5){.3}{-160}{-90}
\psarc(2,.5){.4}{160}{165}
\psarc(1,2){.2}{135}{45}
\psarc[linestyle=dotted](1,2){.2}{45}{135}
\psarc(.2,2){.2}{135}{45}
\psarc[linestyle=dotted](.2,2){.2}{45}{135}
\psarc[linestyle=dotted](1,2){.35}{45}{135}
\psarc[linestyle=dotted](.2,2){.35}{45}{135}
\psline(2,1.3)(2,1.8)
\psarc(2,2){.2}{135}{45}
\psarc[linestyle=dotted](2,2){.2}{45}{135}
\psarc(2.7,2){.2}{135}{45}
\psarc[linestyle=dotted](2.7,2){.2}{45}{135}
\rput[lt](-.7,.4){$U_r$}
\psline{*-}(2,1.3)(2.7,1.8)
\psline[border=1.5pt]{*-}(2,1)(1.5,1.5)
\psline(1.5,1.5)(1.24,1.76)
\psline(2,1)(2,.8)
\psline(1.2,1.2)(1,1.5)
\rput[b](.5,-.2){$m$}
\rput[b](2,-.2){$m^{\prime}$}
\end{pspicture}
\end{center}
First slide the meridian $m$ of one of them inside the other one $m^{\prime}$. Next use move $Y_4$ (Lemma~\ref{lemY4}) to cut $m^{\prime}$ into $m^{\prime}$ and a leaf
that links the edge going to $m$. It is enough to slide inside components that contribute to $\mu(K_r,K_s,K_t)$. Thus,
we do not lose properties of our graphs, (and otherwise we could just perform the surgery).

Last, transform a pair of $K_r$--adjacent $Y$-graphs with opposite contributions to $\mu(K_r,K_s,K_t)$ into a family of 
$Y$-graphs that do not individually contribute to $\mu(K_r,K_s,K_t)$.
To do this, see the effect of the surgery along the two $K_r$--adjacent $Y$-graphs as a band sum with the boundary of a genus two surface
$\Sigma$
whose $1$--handles are $\alpha_1$, $\beta_1$, $\alpha_2$, $\beta_2$, and are meridians of $K_s$ and $K_t$.

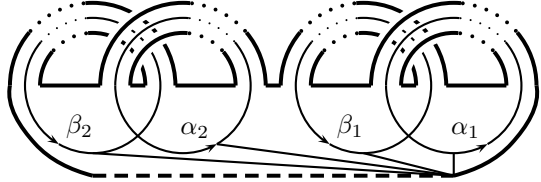
\begin{center}
\begin{pspicture}[.4](0,.2)(5.4,2.9)
\psarc[linewidth=1.5pt](1.1,1.6){1.1}{0}{95}
\psarc[linewidth=1.5pt,linestyle=dotted](1.1,1.6){1.1}{95}{140}
\psarc[linewidth=1.5pt](1.1,1.6){1.1}{140}{180}
\psarc[linewidth=1.5pt](1.1,1.6){.7}{0}{95}
\psarc[linewidth=1.5pt,linestyle=dotted](1.1,1.6){.7}{95}{140}
\psarc[linewidth=1.5pt](1.1,1.6){.7}{140}{180}
\psarc[linestyle=dotted](1.1,1.6){.9}{95}{140}
\psarc{->}(1.1,1.6){.9}{140}{240}
\psarc(1.1,1.6){.9}{240}{95}
\psarc[linewidth=1.5pt,border=2pt](2.3,1.6){1.1}{0}{40}
\psarc[linewidth=1.5pt,border=2pt,linestyle=dotted](2.3,1.6){1.1}{40}{85}
\psarc[linewidth=1.5pt,border=2pt](2.3,1.6){1.1}{85}{180}
\psarc[linewidth=1.5pt,border=2pt](2.3,1.6){.7}{0}{40}
\psarc[linewidth=1.5pt,border=2pt,linestyle=dotted](2.3,1.6){.7}{40}{85}
\psarc[linewidth=1.5pt,border=2pt](2.3,1.6){.7}{85}{180}
\psarc[border=1pt](2.3,1.6){.9}{85}{180}
\psarc{->}(2.3,1.6){.9}{180}{300}
\psarc(2.3,1.6){.9}{-60}{40}
\psarc[border=1pt,linestyle=dotted](2.3,1.6){.9}{40}{85}
\psecurve[linewidth=1.5pt](1.1,2.7)(0,1.6)(.33,.83)(1.1,.4)(1.3,.4)
\psline[linewidth=1.5pt,linestyle=dashed](1.1,.4)(5.9,.4)
\rput[rb](2.65,.9){$\alpha_2$} 
\rput[lb](.75,.9){$\beta_2$} 
\psline(5.9,.4)(1.1,.7)
\psline(5.9,.4)(2.75,.82)
\psline[linewidth=1.5pt](.4,1.6)(1.2,1.6)
\psline[linewidth=1.5pt](1.6,1.6)(1.8,1.6)
\psline[linewidth=1.5pt](2.2,1.6)(3,1.6) 
\psarc[linewidth=1.5pt](4.7,1.6){1.1}{0}{95}
\psarc[linewidth=1.5pt,linestyle=dotted](4.7,1.6){1.1}{95}{140}
\psarc[linewidth=1.5pt](4.7,1.6){1.1}{140}{180}
\psarc[linewidth=1.5pt](4.7,1.6){.7}{0}{95}
\psarc[linewidth=1.5pt,linestyle=dotted](4.7,1.6){.7}{95}{140}
\psarc[linewidth=1.5pt](4.7,1.6){.7}{140}{180}
\psarc[linestyle=dotted](4.7,1.6){.9}{95}{140}
\psarc{->}(4.7,1.6){.9}{140}{240}
\psarc(4.7,1.6){.9}{240}{95}
\psarc[linewidth=1.5pt,border=2pt](5.9,1.6){1.1}{0}{40}
\psarc[linewidth=1.5pt,border=2pt,linestyle=dotted](5.9,1.6){1.1}{40}{85}
\psarc[linewidth=1.5pt,border=2pt](5.9,1.6){1.1}{85}{180}
\psarc[linewidth=1.5pt,border=2pt](5.9,1.6){.7}{0}{40}
\psarc[linewidth=1.5pt,border=2pt,linestyle=dotted](5.9,1.6){.7}{40}{85}
\psarc[linewidth=1.5pt,border=2pt](5.9,1.6){.7}{85}{180}
\psarc[border=1pt](5.9,1.6){.9}{85}{180}
\psarc{->}(5.9,1.6){.9}{180}{300}
\psarc(5.9,1.6){.9}{-60}{40}
\psarc[border=1pt,linestyle=dotted](5.9,1.6){.9}{40}{85}
\psecurve[linewidth=1.5pt](5.9,2.7)(7,1.6)(6.67,.83)(5.9,.4)(5.7,.4)
\rput[rb](6.25,.9){$\alpha_1$} 
\rput[lb](4.35,.9){$\beta_1$} 
\psline(5.9,.4)(4.7,.7)
\psline(5.9,.4)(5.9,.7)
\psline[linewidth=1.5pt](3.4,1.6)(3.6,1.6)
\psline[linewidth=1.5pt](4,1.6)(4.8,1.6)
\psline[linewidth=1.5pt](5.2,1.6)(5.4,1.6)
\psline[linewidth=1.5pt](5.8,1.6)(6.6,1.6)
\end{pspicture}
\end{center}

We are in one of the following situations for the homology classes of the curves:
Either $[\alpha_1]=[\alpha_2]$ and $[\beta_1]=-[\beta_2]$,
or $[\alpha_1]=-[\alpha_2]$ and $[\beta_1]=[\beta_2]$,
or $[\alpha_1]=[\beta_2]$ and $[\beta_1]=[\alpha_2]$,
or $[\alpha_1]=-[\beta_2]$ and $[\beta_1]=-[\alpha_2]$.

Consider the following simple closed curves $c(\alpha_2)$, $c(\beta_1)$, $c(\beta_2)$, $c(\beta_1\alpha_2)$ and $c(\beta_1\beta_2)$ 
whose homology classes are $[\alpha_2]$, $[\beta_1]$, $[\beta_2]$, $[\beta_1\alpha_2]$ and $[\beta_1\beta_2]$, respectively.

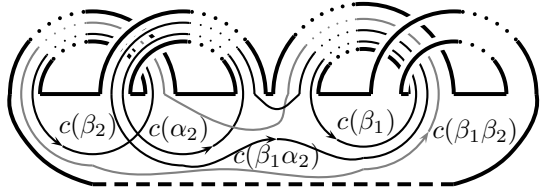
\begin{center}
\begin{pspicture}[.4](0,.2)(5.4,2.9)
\psarc[linecolor=gray](1.1,1.6){.95}{0}{95}
\psarc[linewidth=1.5pt](1.1,1.6){1.1}{0}{95}
\psarc[linewidth=1.5pt,linestyle=dotted](1.1,1.6){1.1}{95}{140}
\psarc[linewidth=1.5pt](1.1,1.6){1.1}{140}{180}
\psarc[linewidth=1.5pt](1.1,1.6){.7}{0}{95}
\psarc[linewidth=1.5pt,linestyle=dotted](1.1,1.6){.7}{95}{140}
\psarc[linewidth=1.5pt](1.1,1.6){.7}{140}{180}
\psarc[linestyle=dotted](1.1,1.6){.8}{95}{140}
\psarc{->}(1.1,1.6){.8}{140}{240}
\psarc(1.1,1.6){.8}{240}{95}
\psarc[linewidth=1.5pt,border=2pt](2.3,1.6){1.1}{0}{40}
\psarc[linewidth=1.5pt,border=2pt,linestyle=dotted](2.3,1.6){1.1}{40}{85}
\psarc[linewidth=1.5pt,border=2pt](2.3,1.6){1.1}{85}{180}
\psarc[linewidth=1.5pt,border=2pt](2.3,1.6){.7}{0}{40}
\psarc[linewidth=1.5pt,border=2pt,linestyle=dotted](2.3,1.6){.7}{40}{85}
\psarc[linewidth=1.5pt,border=2pt](2.3,1.6){.7}{85}{180}
\psarc[border=1pt](2.3,1.6){.8}{85}{180}
\psarc{->}(2.3,1.6){.8}{180}{300}
\psarc(2.3,1.6){.8}{-60}{40}
\psarc[border=1pt,linestyle=dotted](2.3,1.6){.8}{40}{85}
\psecurve[linewidth=1.5pt](1.1,2.7)(0,1.6)(.33,.83)(1.1,.4)(1.3,.4)
\psline[linewidth=1.5pt,linestyle=dashed](1.1,.4)(5.9,.4)
\rput[rb](2.65,.95){$c(\alpha_2)$} 
\rput[lb](.65,1){$c(\beta_2)$} 
\psline[linewidth=1.5pt](.4,1.6)(1.2,1.6)
\psline[linewidth=1.5pt](1.6,1.6)(1.8,1.6)
\psline[linewidth=1.5pt](2.2,1.6)(3,1.6) 
\psarc[linewidth=1.5pt](4.7,1.6){1.2}{0}{95}
\psarc[linewidth=1.5pt,linestyle=dotted](4.7,1.6){1.2}{95}{140}
\psarc[linewidth=1.5pt](4.7,1.6){1.2}{140}{180}
\psarc[linewidth=1.5pt](4.7,1.6){.6}{0}{95}
\psarc[linewidth=1.5pt,linestyle=dotted](4.7,1.6){.6}{95}{140}
\psarc[linewidth=1.5pt](4.7,1.6){.6}{140}{180}
\psarc[linestyle=dotted](4.7,1.6){.7}{95}{140}
\psarc{->}(4.7,1.6){.7}{140}{240}
\psarc(4.7,1.6){.7}{240}{95}

\rput[l](5.62,1.1){$c(\beta_1\beta_2)$} 
\psarc[linecolor=gray,linestyle=dotted](4.7,1.6){1}{95}{140}
\psarc[linecolor=gray](4.7,1.6){1}{140}{180}
\psarc[linecolor=gray](4.7,1.6){1}{-30}{95}
\psarc[linecolor=gray]{->}(4.7,1.6){1}{-90}{-30}
\psarc[linecolor=gray,linestyle=dotted](1.1,1.6){.95}{95}{140}
\psarc[linecolor=gray](1.1,1.6){.95}{140}{270}
\psecurve[linecolor=gray](.15,1.6)(1.1,.65)(2.9,.6)(4.7,.6)(5.7,1.6)
\psecurve[linecolor=gray](1.1,2.55)(2.05,1.6)(3.55,1.2)(3.7,1.6)(4.7,2.6)
\rput[t](3.55,.95){$c(\beta_1\alpha_2)$} 
\psarc[linestyle=dotted](4.7,1.6){.85}{95}{140}
\psarc(4.7,1.6){.85}{140}{180}
\psarc(4.7,1.6){.85}{270}{95}
\psarc(2.3,1.6){.95}{180}{270}
\psarc(2.3,1.6){.95}{0}{40}
\psarc[border=1pt,linestyle=dotted](2.3,1.6){.95}{40}{85}
\psarc[border=1pt](2.3,1.6){.95}{85}{180}
\psecurve{->}(1.35,1.6)(2.3,.65)(3.55,1)(4.7,.75)
\psecurve(2.3,.65)(3.55,1)(4.7,.75)(5.55,1.6)
\psecurve(2.3,2.55)(3.25,1.6)(3.55,1.4)(3.85,1.6)(4.7,2.45)

\psarc[linewidth=1.5pt,border=2pt](5.9,1.6){1.1}{0}{40}
\psarc[linewidth=1.5pt,border=2pt,linestyle=dotted](5.9,1.6){1.1}{40}{85}
\psarc[linewidth=1.5pt,border=2pt](5.9,1.6){1.1}{85}{180}
\psarc[linewidth=1.5pt,border=2pt](5.9,1.6){.7}{0}{40}
\psarc[linewidth=1.5pt,border=2pt,linestyle=dotted](5.9,1.6){.7}{40}{85}
\psarc[linewidth=1.5pt,border=2pt](5.9,1.6){.7}{85}{180}
\psecurve[linewidth=1.5pt](5.9,2.7)(7,1.6)(6.67,.83)(5.9,.4)(5.7,.4) 
\rput[lb](4.35,1.05){$c(\beta_1)$} 
\psline[linewidth=1.5pt](3.4,1.6)(3.5,1.6)
\psline[linewidth=1.5pt](4.1,1.6)(4.8,1.6)
\psline[linewidth=1.5pt](5.2,1.6)(5.3,1.6)
\psline[linewidth=1.5pt](5.9,1.6)(6.6,1.6)
\end{pspicture}
\end{center}

For a curve $c$, let $\tau_c$ denote the right-handed Dehn twist around this curve.
Recall the action of $\tau$ on homology classes $\tau_c(x) = x + \langle c , x\rangle_{\Sigma} c$.
Then the homeomorphism $\tau_{c(\alpha_2)}^{-1}\tau_{c(\beta_1)}^{-1}\tau_{c(\beta_1\alpha_2)}$ of $\Sigma$
transforms $\alpha_2$ and $\beta_1$ to conjugate curves, where the conjugation paths are in the neighborhood of the genus $2$ surface and avoids the disks $D_i$, for $i \leq r$, 
and it transforms $\alpha_1$ and $\beta_2$ into curves homologous to $\alpha_1 \alpha_2^{-1}$ and $\beta_1\beta_2$.
Therefore using this boundary-preserving homeomorphism in the first case allows us to transform the surgery
on the initial pair of $Y$-graphs into a surgery on a pair of $Y$-graphs such that each of the graphs has a homologically
trivial leaf and two meridional leaves.
In the second case, use $\tau_{c(\alpha_2)}\tau_{c(\beta_1)}\tau_{c(\beta_1\alpha_2)}^{-1}$.
Use $\tau_{c(\beta_2)}^{-1}\tau_{c(\beta_1)}^{-1}\tau_{c(\beta_1\beta_2)}$ and $\tau_{c(\beta_2)}\tau_{c(\beta_1)}\tau_{c(\beta_1\beta_2)}^{-1}$ in the third and fourth cases, respectively to achieve 
a similar reduction.
\eop

\section{Proof of the formulae for algebraically split links}
\setcounter{equation}{0}
\label{secprooffas}
 
We prove the surgery formulae of Theorem~\ref{thmfas} and Theorem~\ref{thmfasmu}
following a strategy that was used in \cite{ggp} to compare the filtration
of the space of $\ZZ$--spheres associated to algebraically split links to the filtration associated to $Y$-links.

According to Proposition~\ref{propclaspmu}, it is enough to prove these theorems for links that are presented by
pairs $(G,U^{(r)})$ where $G$ is a $Y$-link
that $\mu$-laces the trivial link $U^{(r)}$ of $M$, that is for $U^{(r)} \subset M_G$, where $U^{(r)}$ is equipped with surgery coefficients $p_1/q_1$, $p_2/q_2$, \dots, $p_r/q_r$.

$$[M_G;U^{(r)}] = \sum_{H \subset G} (-1)^{\sharp H} [M;H \cup U^{(r)}]$$
where $$[M;H \cup U^{(r)}]= \sum_{J \subset H, I \subset \{1,2,\dots,r\}}(-1)^{\sharp J + \sharp I}M_{J,(U_i;p_i/q_i)_{i \in I}} \sharp \left(\sharp_{j \in \{1,2,\dots,r\} \setminus I} L(p_j,-q_j)\right).$$
Then $[M_G;U^{(r)}] \stackrel{n}{\equiv} \sum_{H \subset G ;\sharp H \leq 2n} (-1)^{\sharp H} [M;H \cup U^{(r)}]$.\\
If there exists $i$, such that $U_i$ does not link any leaf of $H$, then $[M;H \cup U^{(r)}]=0$.\\
If there exists $i$, such that $U_i$ links only one leaf of $H$, then let $Y_1$ be the component of this leaf.  
$$[M;H \cup U^{(r)}]=- [M_{Y_1}; H \setminus Y_1 \cup U^{(r)}].$$

Recall that the surgery on $Y_1$ is a surgery associated with a genus one surface
bounded by some $K_i$ as in Subsection~\ref{sublagboun}. Then the inverse surgery
of this subsection transforms $U_i$ into $-K_i$ and since it can be realized as a genus one
cobordism, it can also be realized by a surgery on a $Y$-graph that laces $U_1$
and that sits in the complement of $G$. Let $Y_1^{-1}$ be such a graph. We can assume that its leaves are a meridian of $U_i$ and two leaves parallel to the two other leaves of $Y_1$. Compare with Lemma~\ref{lemY6}.

Then $[M_{Y_1}; H \setminus Y_1 \cup U^{(r)}]= [M_{Y_1^{-1}}; H \setminus Y_1 \cup U^{(r)}]$
and $$[M;H \cup U^{(r)}] = \frac{1}{2} [M;H \cup Y_1^{-1} \cup U^{(r)}].$$

As long as there is a component $U_j$ that bounds a disk $D_j$ intersecting $H \cup Y_1^{-1} \cup \dots \cup Y_k^{-1}$
once (and necessarily) inside a meridional leaf of some component $Y_{k+1}$ of $H$, add $Y_{k+1}^{-1}$, and write
$$[M;H \cup Y_1^{-1} \cup \dots \cup Y_k^{-1} \cup U^{(r)}]=
\frac{1}{2} [M;H \cup Y_1^{-1} \cup \dots \cup Y_k^{-1} \cup Y_{k+1}^{-1} \cup U^{(r)}].$$
$$[M;H \cup U^{(r)}]=\frac{1}{2^{k+1}}[M;H \cup Y_1^{-1} \cup \dots \cup Y_k^{-1} \cup Y_{k+1}^{-1} \cup U^{(r)}].$$
Finally,
$[M_G;U^{(r)}]$ is a rational combination of terms of the form $[M; H^{\prime} \cup U^{(r)}]$ where each $U_i$ links at least
two leaves of $H^{\prime}$. To be more specific, the considered $H^{\prime}$ are of the form $H \cup H_1^{-1}$,
where $H$ is a sublink of $G$, and $H_1^{-1}$ is a link made of inverses of the components of a sublink $H_1$
of $H$. In particular, the leaves of a component of $H_1^{-1}$ have the same constraints as the leaves
of a component of $G$.
Since a leaf of $H^{\prime}$ links at most one $U_i$, 
such a $H^{\prime}$ has at least $2r$ leaves linking the $U_i$. In particular, if $2r>6n$, $[M_G;U^{(r)}] \stackrel{n}{\equiv} 0$.
\begin{itemize}
\item Under the hypotheses of Theorem~\ref{thmfas}, assume $2r=6n$.
Up to elements in $\CF_{2n+1}$, $[M_G;U^{(r)}]$ is a rational combination of terms 
of the form $[M; H^{\prime} \cup U^{(r)}]$ where each $U_i$ links exactly
two leaves of $H^{\prime}$, and each leaf of $H^{\prime}$ is a meridional leaf of some $U_i$.
More precisely, let $G_3$ be the sublink of $G$ made of the components that have three meridional leaves, we have
$$[M_G;U^{(r)}]\stackrel{n}{\equiv} \sum_{H}(-1)^{\sharp H}[M;H\cup U^{(r)}]$$
where the sum runs over the $H$ that read as the disjoint union of two $Y$-links
$H_1$ and $H_2$ of $G_3$ such that for any component $U_i$ of $U^{(r)}$, either there
is one meridional leaf of $U_i$ in $H_1$ and no meridional leaf of $U_i$ in $H_2$, or 
there
is no meridional leaf of $U_i$ in $H_1$ and there are two meridional leaves of $U_i$ in $H_2$.
Let $\CH$ denote the set of the $(H_1,H_2)$ where $H_1 \cup H_2$ is a decomposition as above of such a graph.
\item Under the hypotheses of Theorem~\ref{thmfasmu}, at most two thirds of the leaves of the $H^{\prime}$
link the $U_i$ once, and the leaves of the other third do not link the $U_i$ at all.
Therefore, if $2r>4n$, $[M_G;U^{(r)}]$ belongs to $\CF_{2n+1}$. If $r=2n$,
up to elements in $\CF_{2n+1}$, $[M_G;U^{(r)}]$ is a rational combination of terms 
of the form $[M; H^{\prime} \cup U^{(r)}]$ where each $U_i$ links exactly
two leaves of $H^{\prime}$, and in each component of $H^{\prime}$, there are 
two meridional leaves of $U^{(r)}$ and a null-homologous leaf.
More precisely, let $G_2$ be the sublink of $G$ made of the components that have two meridional leaves, we have
$$[M_G;U^{(r)}]\stackrel{n}{\equiv} \sum_{H}(-1)^{\sharp H}[M;H\cup U^{(r)}]$$
where the sum runs over the $H$ that read as the disjoint union of two $Y$-links
$H_1$ and $H_2$ of $G_2$ such for any component $U_i$ of $U^{(r)}$, either there
is one meridional leaf of $U_i$ in $H_1$ and no meridional leaf of $U_i$ in $H_2$, or 
there
is no meridional leaf of $U_i$ in $H_1$ and there are two meridional leaves of $U_i$ in $H_2$.
Let $\CH$ denote the set of the $(H_1,H_2)$ where $H_1 \cup H_2$ is a decomposition as above of such a graph.
\end{itemize}

\noindent In both cases
$$[M_G;U^{(r)}]\stackrel{n}{\equiv} \sum_{(H_1,H_2)\in \CH}\left(\frac{-1}{2}\right)^{\sharp H_1} (-1)^{\sharp H_2}[M;H_1 \cup H_1^{-1} \cup H_2 \cup U^{(r)}]$$
where $H_1 \cup H_1^{-1} \cup H_2$ has $2n$ components (and therefore $(-1)^{\sharp H_2}=1$). Apply Theorem~\ref{thmflag} to compute it.
The tripod associated to a surgery on an oriented $Y$-graph whose leaves are $\ell_1,\ell_2,\ell_3$ was computed in
Lemma~\ref{lemtriplag} (thanks to Proposition~\ref{propborlag}).
It is 
\begin{pspicture}[0.4](-.05,-.1)(.9,.7) 
\psline{-}(0.05,.3)(.45,.6)
\psline{*-}(0.05,.3)(.45,.3) 
\psline{-}(0.05,.3)(.45,0)
\rput[l](.55,0){$\ell_3$}
\rput[l](.55,.3){$\ell_2$}
\rput[l](.55,.6){$\ell_1$}
\end{pspicture}
while the tripod
associated to an inverse of such a graph is
\begin{pspicture}[0.4](-.05,-.1)(1,.7) 
\psline{-}(0.05,.3)(.45,.6)
\psline{*-}(0.05,.3)(.45,.3) 
\psline{-}(0.05,.3)(.45,0)
\rput[l](.55,0){$\ell_{1 \parallel}$}
\rput[l](.55,.3){$\ell_{2 \parallel}$}
\rput[l](.55,.6){$\ell_{3 \parallel}$}
\end{pspicture}
where the parallels are taken with respect of the parallelization of the leaves.
Later, we shall consider twice the tripods of the components of $H_1$ and remove the $(-1)^{\sharp H_1}$.
Recall the formula of Lemma~\ref{lemvarlk}
$$lk_{M_{(U_i;p_i/q_i)}}(J,K)=lk_M(J,K)-\frac{q_i}{p_i}lk_M(U_i,J)lk_M(U_i,K).$$
If for some $i$, a contraction does not pair two curves linking $U_i$, then its contribution
to $[M;H^{\prime} \cup (U \setminus U_i)]$ and its contribution to $[M_{U_i};H^{\prime} \cup (U \setminus U_i)]$
will be the same. Therefore, it won't contribute to $[M;H^{\prime} \cup U ]$. Thus since there are exactly
two leaves $m_i$ and $m^{\prime}_i$ linking $U_i$ in each $H^{\prime}$, the only pairings that will contribute 
will pair these pairs together, and
the corresponding remaining linking number will be $\frac{q_i}{p_i}lk_M(U_i,m_i)lk_M(U_i,m^{\prime}_i)$.

\begin{itemize}
\item Under the hypotheses of Theorem~\ref{thmfas},
there is one contributing pairing for every $(H_1,H_2)\in \CH$. It can be seen as an edge-labelled Jacobi diagram
$\Gamma(H_1,H_2)$
together with a bijection $b$ from the set of its vertices to the set of components of $H_1 \cup H_1^{-1}\cup H_2$
that maps a vertex $v$ with adjacent edges labelled by $i,j,k$
to a component $b(v)$ of $G$ of type $(\varepsilon_i m_i, \varepsilon_j m_j,\varepsilon_k m_k)$ where $\varepsilon_i$, $\varepsilon_j$, $\varepsilon_k$
are in $\{-1,1\}$, or to the inverse of such a component.
Equip $\Gamma(H_1,H_2)$ with an orientation. 
Then if the orientation of a vertex $v$ as above is induced by the cyclic order $(i,j,k)$, assign it the sign $(-\varepsilon_i \varepsilon_j \varepsilon_k)$, and assign it $\varepsilon_i \varepsilon_j \varepsilon_k$ otherwise.
Define $\mbox{sign}(\Gamma(H_1,H_2);b)$ as the products of the signs of the vertices.
Then 
$$Z_n\left([M_G;U^{(r)}]\right)= \prod_{i=1}^{3n}\frac{q_i}{p_i}\sum_{(H_1,H_2)\in \CH}\frac{1}{2^{\sharp H_1}}\mbox{sign}(\Gamma(H_1,H_2);b)[\Gamma(H_1,H_2)]$$

Now, let $f=f(b)$ be the map from $V(\Gamma(H_1,H_2))$ to the set of components of $G$ obtained from a bijection $b$ as above by setting
$$\begin{array}{lll}
f(b)(v)&=b(v)& \mbox{if}\; b(v) \; \mbox{is a component of }\;
H_1 \cup H_2\\
&=Y_i& \mbox{if}\; b(v)=Y_i^{-1}.
\end{array}$$
There are $2^{\sharp H_1}$ $b$ such that $f(b)=f$, 
and, if $\sharp \mbox{Aut}_e(\Gamma)$ is the set of automorphisms of $\Gamma$ that induce the Identity on $E(\Gamma)$, there are $\sharp \mbox{Aut}_e(\Gamma)$ $b$ that define the same pairing. Since an automorphism that preserves the edges pointwise may only exchange vertices inside components $\tata$, $\sharp \mbox{Aut}_e(\Gamma) = 2^{\theta(\Gamma)}$.

Orient $G$ arbitrarily.
Let $\Gamma \in D_{e,n}$. Equip $\Gamma$ with an arbitrary orientation.
Let $G(\Gamma)$ be the set of maps $g$ from $V(\Gamma)$ to the set of components of $G$ that map a vertex $v$ with adjacent edges labelled by $i,j,k$, with respect to the order induced by the orientation,
to a component $g(v)$ of $G$ of type $(\varepsilon_i m_i, -\varepsilon_j m_j,\varepsilon_k m_k)$ or $(\varepsilon_i m_i, \varepsilon_k m_k, \varepsilon_j m_j)$. Define $\mbox{sign}(g,v)=\varepsilon_i \varepsilon_j \varepsilon_k$ for such a vertex.
Define $\mbox{sign}(\Gamma;g)$ as the product of the signs associated to the vertices.
Then
$$Z_n\left([M_G;U^{(r)}]\right)= \prod_{i=1}^{3n}\frac{q_i}{p_i}
\sum_{\Gamma \in D_{e,n}, g \in G(\Gamma)}\frac{\mbox{sign}(\Gamma;g)}{2^{\theta(\Gamma)}}[\Gamma].$$
Now, Lemma~\ref{lemintclas} easily leads to the conclusion of the proof of Theorem~\ref{thmfas}.

\item Orient $G$ arbitrarily.
Under the hypotheses of Theorem~\ref{thmfasmu}, a contributing pairing for $(H_1,H_2)\in \CH$
is a $2/3$--labelled Jacobi diagram $\Gamma$, equipped with a bijection from $V(\Gamma)$
to the set of components of $H_1 \cup H_1^{-1}\cup H_2$ that maps a vertex with two adjacent edges labelled by $i$ and $j$ to a component of type
$(\varepsilon_i m_i, \varepsilon_j m_j,f)$ or $(\varepsilon_j m_j, -\varepsilon_i m_i,f)$. For a fixed $2/3$--labelled Jacobi diagram $\Gamma$, there are $\sharp \mbox{Aut}_{2/3}(\Gamma)$ 
bijections from $V(\Gamma)$ to the set of components of $H_1 \cup H_1^{-1}\cup H_2$ that will correspond to the same pairing.

Let $\Gamma \in D_{2/3,n}$.
Equip $\Gamma$ with an orientation. Let $G(\Gamma)$ be the set of maps $g$ from $V(\Gamma)$
to the set of components of $G$ that map a vertex $v$ whose adjacent edges are labelled by $(i,j,\mbox{nothing})$ (with respect to the orientation of $\Gamma$)
to a component of type
$(\varepsilon_i m_i, \varepsilon_j m_j,f)$ or $(\varepsilon_j m_j, -\varepsilon_i m_i,f)$ of $G$. When $g \in G(\Gamma)$ is fixed, assign the framed oriented curve $\varepsilon_i\varepsilon_j f$ to the unlabelled edge of each $v\in V(\Gamma)$ as above. Then assign to each edge of $\Gamma$ the linking number of
the two curves assigned to its half-edges (change a curve $f$ into its parallel $f_{\parallel}$, if the two curves coincide) and define $lk(\Gamma;g)$ as the product over the edges of $\Gamma$
of the associated linking numbers.
$$Z_n\left([M_G;U^{(r)}]\right)= \prod_{i=1}^{2n}\frac{q_i}{p_i}
\sum_{\Gamma \in D_{2/3,n}, g \in G(\Gamma)}\frac{lk(\Gamma;g)}{\sharp \mbox{Aut}_{2/3}(\Gamma)}[\Gamma].$$
Now, Lemma~\ref{lemintclas} easily leads to the conclusion of the proof of Theorem~\ref{thmfasmu} when the Seifert surfaces are associated to a presentation
of the link by a graph that $\mu$-laces the unlink as in Proposition~\ref{propclaspmu}.
Fortunately, this is enough to conclude the proof of Theorem~\ref{thmfasmu}
thanks to the following proposition~\ref{propasmu} that ensures that the right-hand side of 
the equality of Theorem~\ref{thmfasmu} does not depend on the choice of the Seifert surfaces.
\end{itemize}
\eop

Let $n \in \NN$.
Let $D_{2/3,o,n}$ be the set of $2/3$-labelled unoriented Jacobi diagrams whose labelled edges 
are oriented. Forgetting the edge orientations
transforms an element $\Gamma$ of $D_{2/3,o,n}$ into an element $f(\Gamma)$ of $D_{2/3,n}$, and an element of $D_{2/3,n}$ comes from $\frac{2^{2n}}{\sharp \mbox{Aut}_{2/3}(\Gamma)}$ elements of
$D_{2/3,o,n}$.

Let $L=(K_i;p_i/q_i)_{i \in \{1,2,3,\dots, 2n\}}$ be a  framed $2n$--component algebraically split link in a rational homology sphere $\;M$.
Assume that for any $\{i,j,k\} \subset \{1,2,3,\dots, 2n\}$, $\mu(K_i,K_j,K_k)=0$.
Let $(F_i^-)_{i \in \{1,2,3,\dots, 2n\}} \cup (F_i^+)_{i \in \{1,2,3,\dots, 2n\}}$ be a collection of transverse surfaces such that, for any $i$, $F_i^-$ and $F_i^+$ are two Seifert surfaces of $K_i$ that do not meet
the $K_j$ for $j \neq i$. 

Let $\Gamma \in D_{2/3,o,n}$. Orient $\Gamma$. In such a $\Gamma$ the half-edges of the labelled edges inherit a label from the edge orientation. Namely, Edge $i$ goes from $i^-$ to $i^+$.

For any vertex of $\Gamma$, whose half-edges are labelled by $(i^{\varepsilon},j^{\eta},$ nothing$)$ with respect to the cyclic order induced by the orientation, assign the intersection curve $F_i^{\varepsilon} \cap F_j^{\eta}$ to its unlabelled half-edge.
To any unlabelled edge that is now equipped with two intersection curves associate the linking number of these curves.
Then define $\ell_{\Gamma}((F_i^-,F_i^+)_{i=1,\dots 2n})$ as the product over all the unlabelled edges of $\Gamma$ of the corresponding linking numbers. 
Note that $\ell_{\Gamma}((F_i^-,F_i^+)_{i=1,\dots 2n})[\Gamma]$ does not depend on the orientation of $\Gamma$.

When $F_i^+$ is a parallel copy of $F_i^-$, then 
$$\sum_{\Gamma \in D_{2/3,n}} \frac{\ell((F_i^-)_{i=1, \dots 2n};\Gamma)}{\sharp \mbox{Aut}_{2/3}(\Gamma)}[\Gamma]=\sum_{\Gamma \in D_{2/3,o,n}} \ell_{\Gamma}((F_i^-,F_i^+)_{i=1,\dots 2n})\frac{1}{2^{2n}}[\Gamma]$$ 

\begin{proposition}
\label{propasmu}
With the notation and hypotheses above
$$\sum_{\Gamma \in D_{2/3,o,n}} \ell_{\Gamma}((F_i^-,F_i^+)_{i=1,\dots 2n})\frac{1}{2^{2n}}[\Gamma]$$ 
 is independent of the choice of the surfaces $(F_i^-,F_i^+)_{i=1,\dots 2n}$ in the complement of $\cup_{j\neq i}K_j$, it only depends on $L$.
\end{proposition}
\bp We study the effect of changing a surface $F_i^{\varepsilon}$ to another
Seifert surface $F^{\prime}$ of $K_i$ disjoint from
the $K_j$ for $i \neq j$, and transverse to the other ones. Obviously, for any
$\Gamma$, the only modified ingredient is the linking number associated with the unlabelled
edge $e$ that shares a vertex with $i^{\varepsilon}$ that reads
$$\pm lk(F_i^{\varepsilon} \cap S_1, S_2 \cap S_3)$$
where $S_1$, $S_2$ and $S_3$ are the three other surfaces associated to the three
other labelled half-edges containing the vertices of $e$.

Let us compute the variation of such a linking number.
Recall that $H_2(M \setminus \cup_{j=1,2,\dots,2n}K_j)$ is generated by the homology classes
of the boundaries $\partial N(K_j)$ of the tubular neighborhoods of the $K_j$, for $j\neq i$.
Therefore the immersed oriented closed surface
$(F_i^{\varepsilon} \cup - F^{\prime})$ cobounds a $3$-dimensional chain $C$
with some copies $\partial N(K_j)$.
In particular, if $S_1$ is a Seifert surface for $K_{j(1)}$, the boundary of
$C\cap S_1$ is the union of $(F^{\prime} \cap S_1 -F_i^{\varepsilon} \cap S_1)$ and some copies of $K_{j(1)}$.
Since all the Milnor triple linking numbers vanish, $lk(K_{j(1)},S_2 \cap S_3)=0$,
and $$lk(F^{\prime} \cap S_1-F_i^{\varepsilon} \cap S_1, S_2 \cap S_3)
=\pm \langle C\cap S_1, S_2 \cap S_3 \rangle=\pm \langle C, S_1 \cap S_2 \cap S_3 \rangle.$$
Now, consider the two elements of $D_{2/3,o,n}$ obtained from $\Gamma$ by changing
the neighborhood of $e$ in $\Gamma$ as in the following figure.

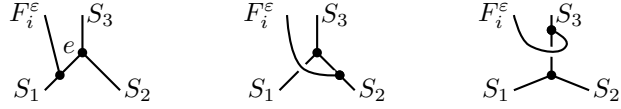
\begin{center}
\begin{pspicture}[.2](-1,-.2)(2,1.4)
\psline{-*}(0,1)(.2,.2)
\psline{*-}(.5,.5)(.5,1)
\psline{-}(1,0)(.5,.5)
\psline{-}(0,0)(.5,.5)
\rput[br](.4,.5){$e$}
\rput[r](-.05,1){$F_i^{\varepsilon}$}
\rput[r](-.05,0){$S_1$}
\rput[l](1.05,0){$S_2$}
\rput[l](.55,1){$S_3$}
\end{pspicture}
\begin{pspicture}[.2](-1,-.2)(2,1.4)
\psline{*-}(.5,.5)(.5,1)
\psline{-}(1,0)(.5,.5)
\psline{-}(0,0)(.5,.5)
\pscurve[border=1.5pt]{-*}(.1,1)(.3,.3)(.8,.2)
\rput[r](-.05,0){$S_1$}
\rput[l](1.05,0){$S_2$}
\rput[r](-.05,1){$F_i^{\varepsilon}$}
\rput[l](.55,1){$S_3$}
\end{pspicture}
\begin{pspicture}[.2](-1,-.2)(2,1.4)
\psline{*-}(.5,.2)(.5,1)
\psline{-}(1,0)(.5,.2)
\psline{-}(0,0)(.5,.2)
\pscurve[border=1.5pt]{-*}(0,1)(.2,.6)(.7,.6)(.5,.8)
\rput[r](-.05,1){$F_i^{\varepsilon}$}
\rput[r](-.05,0){$S_1$}
\rput[l](1.05,0){$S_2$}
\rput[l](.55,1){$S_3$}
\end{pspicture} 
\end{center}

(Actually, since the current definition of Jacobi diagrams does not allow looped edges, some of the above graphs may not be Jacobi diagrams. In order to make this proof work, allow Jacobi diagrams with looped edges, and set them to be zero in $\CA_n(\emptyset)$, so that the IHX relations involving such graphs are still valid and these graphs do not contribute to the sum of the statement.)
 
Assume without loss, that the orientations 
of the three graphs coincide outside the neighborhood of $e$ and are induced by the figure at the shown vertices.
Then the coefficients of these three elements of $D_{2/3,o,n}$ are perturbed in the same way. (Note that we did not need to take care about the above signs, they are well-defined in each step, and the result only depends on the cyclic order
of $S_1$, $S_2$, $S_3$.) Since the sum of the corresponding oriented graphs
vanishes in $\CA_n(\emptyset)$ and since all the graphs of $D_{2/3,o,n}$ can be grouped in three-element sets as above, the sum of the statement is independent of the surfaces. 
\eop

Similarly, we can show the following proposition.

\begin{proposition}
\label{propmilfour}
Let $L=(K_0,K_1,K_2,K_3)$ be a rationally algebraically split link
whose three-component sublinks have Milnor triple linking number $0$ in a rational homology sphere $M$. Let
$a$, $b$ and $c$ be three real numbers such that $a+b+c=0$.
Let $\Sigma_i$ be a Seifert surface for $K_i$ in the exterior of $L\setminus K_i$.  Then
$$\nu_{abc}(K_0,K_1,K_2,K_3) = $$
$$a lk(\Sigma_0 \cap \Sigma_1, \Sigma_2 \cap \Sigma_3) + b lk(\Sigma_0 \cap \Sigma_2, \Sigma_3 \cap \Sigma_1)
+c lk(\Sigma_0 \cap \Sigma_3, \Sigma_1 \cap \Sigma_2)$$
does not depend on the surfaces $\Sigma_i$ that satisfy the given assumption.
The invariant $\nu_{abc}$ satisfies the following properties.
\begin{itemize}
\item It is invariant under self-crossing changes of the components of $L$.
\item If $M=S^3$, $\nu_{abc}$ is the following combination of the Milnor invariants defined in \cite{mil}, 
$$\nu_{abc}=b \mu(10,23) -c\mu(01,23).$$
$\mu(01,23)=\nu_{1,0,-1}$ and $\mu(10,23)=\nu_{-1,1,0}$.
\end{itemize}
\end{proposition}
\bp
The proof of Proposition~\ref{propasmu} shows that $\nu_{abc}$ does not depend on the surfaces and that it is therefore well-defined.
Let us prove that $\nu_{abc}$ does not vary under self-crossing changes and is therefore a homotopy invariant of these four-component links.
To study the effect of a self-crossing change on $K_0$ inside a ball $B$, choose the surfaces $\Sigma_i$ for $i>0$ so that they intersect $B$ as parallel tubes around one strand of $K_0$.
Then their intersections like $\Sigma_2 \cap \Sigma_3$ will not meet $B$, and will also bound a surface $\Sigma_{23}$ in the exterior of $K_0$ and $K_1$ that intersects $B$ as parallel tubes around the same strand of $K_0$.
Now, $\Sigma_1 \cap \Sigma_{23}$ does not meet $B$, and then $$lk(\Sigma_0 \cap \Sigma_1, \Sigma_2 \cap \Sigma_3)=\pm lk(K_0,\Sigma_1 \cap \Sigma_{23})$$ does not vary under the considered crossing change of $K_0$.
 
According to \cite{mil}, if the ambient $3$-manifold is $S^3$, there is a bijection from the set 
of homotopy classes of four-component algebraically split links $L$
whose three-component sublinks have Milnor triple linking number $0$ to $\ZZ \oplus \ZZ$ that maps $L$ to $(\mu(01,23)(L),\mu(10,23)(L))$.

Furthermore, if $(K_0,K_1,K_2)$ is the trivial three-component link with meridians $\alpha_0,\alpha_1,\alpha_2$, and if the homotopy class of 
$K_{01}$ in the exterior of $(K_0,K_1,K_2)$ reads (with the notation of \cite{mil}),
$$\alpha_{2}^{k_0k_1}=\alpha_0\alpha_1 \alpha_2 (\alpha_0\alpha_1)^{-1}
(\alpha_0 \alpha_2^{-1} \alpha_0^{-1})
\alpha_2(\alpha_1 \alpha_2^{-1} \alpha_1^{-1})
 =[\alpha_0,[\alpha_1,\alpha_2]]$$
then $\mu(01,23)(K_0,K_1,K_2,K_{01})=1$ and $\mu(10,23)(K_0,K_1,K_2,K_{01})=0$.
More generally, if the homotopy class of 
$K_{3}$ reads $[\alpha_0,[\alpha_1,\alpha_2]]^{\mu_{01}}[\alpha_1,[\alpha_0,\alpha_2]]^{\mu_{10}}$, then
$\mu(01,23)(K_0,K_1,K_2,K_3)=\mu_{01}$ 
and $\mu(10,23)(K_0,K_1,K_2,K_{3})=\mu_{10}$.
The link presented by the following clasper 
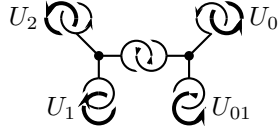
\begin{center}
\begin{pspicture}[.2](-1.2,-.7)(2,1.6) 
\psarc[linewidth=1.5pt]{->}(1,0){.2}{0}{92}
\psarc[border=2pt](1,.2){.2}{90}{270}
\psarc[border=2pt](1,.2){.2}{-90}{90}
\psarc[linewidth=1.5pt,border=2pt]{->}(1,0){.2}{90}{360}
\psarc{->}(.5,.7){.2}{0}{180}
\psarc{->}(.3,.7){.2}{180}{360}
\psarc[border=2pt](.3,.7){.2}{0}{180}
\psarc[border=2pt](.5,.7){.2}{180}{0}
\psarc[linewidth=1.5pt](1.5,1.2){.2}{-2}{180}
\psarc(1.3,1.2){.2}{180}{0}
\psarc[border=2pt]{->}(1.3,1.2){.2}{0}{180}
\psarc[linewidth=1.5pt,border=2pt]{->}(1.5,1.2){.2}{180}{360}
\psarc[linewidth=1.5pt](-.2,0){.2}{180}{90}
\psarc[border=2pt](-.2,.2){.2}{90}{270}
\psarc[border=2pt](-.2,.2){.2}{-90}{90}
\psarc[linewidth=1.5pt,border=2pt]{-<}(-.2,0){.2}{90}{180}
\psarc[linewidth=1.5pt](-.2,0){.2}{170}{190}
\psarc(-.5,1.2){.2}{0}{180}
\psarc[linewidth=1.5pt](-.7,1.2){.2}{178}{0}
\psarc[linewidth=1.5pt,border=2pt]{->}(-.7,1.2){.2}{0}{180}
\psarc[border=2pt]{->}(-.5,1.2){.2}{180}{360}
\psline(.7,.7)(1,.7)(1.3,1)
\psline{*-}(1,.7)(1,.4)
\psline{*-}(-.2,.7)(.1,.7)
\psline(-.2,.4)(-.2,.7)(-.5,1)
\rput[l](1.3,0){$U_{01}$}
\rput[r](-.5,0){$U_{1}$}
\rput[l](1.8,1.2){$U_{0}$}
\rput[r](-1,1.2){$U_{2}$}
\end{pspicture}
\end{center}
has the same properties than $(K_0,K_1,K_2,K_{01})$ and, according to Lemma~\ref{lemintclas},
$$\nu_{abc}((K_0,K_1,K_2,K_{01}))=-c.$$
More generally, if the homotopy class of 
$K_{3}$ reads $[\alpha_0,[\alpha_1,\alpha_2]]^{\mu_{01}}[\alpha_1,[\alpha_0,\alpha_2]]^{\mu_{10}}$, then 
$$\nu_{abc}(K_0,K_1,K_2,K_{3})=(b\mu_{10}- c\mu_{01})(K_0,K_1,K_2,K_{3}).$$
\eop

\section{On the polynomial form of the knot surgery formula: proofs and remarks}
\setcounter{equation}{0}
\label{secproofpol}

\noindent{\sc Proof of Theorem~\ref{thmpol}:}
Since the theorem easily follows from Theorem~\ref{thmfboun} for $n=1$, we assume $n \geq 2$.
First assume $(p,q)=(1,0)$. A $\frac{1}{r}$-surgery on $K$ is equivalent to $|r|$ $\mbox{sign}(r)$-surgeries on parallel copies on $K$.
These parallel copies form an $|r|$-component boundary link $L$ bounding parallel copies of $F$.
We have 
$$M(K;\frac{1}{r})= \sum_{J \subset L} (-1)^{\sharp J} [M;J]$$
Up to elements of
$\mbox{Ker}(Z_n)$, we only consider the sublinks $J$ of $L$
with at most $n$ components, according to Theorem~\ref{thmfboun}. There are $$\left(\begin{array}{c}|r|\\{j}\end{array}\right)=\frac{|r|(|r|-1)\dots(|r|-j+1)}{j!}$$
sublinks $J$ of $L$ with $j$ components and they are all isomorphic to the boundary
link $L_j$ whose components are framed by $\mbox{sign}(r)$. This shows that
$$\begin{array}{lll}
Z_n(M(K;\frac{1}{r}))-Z_n(M)&=\sum_{i=1}^n Y_{n,0}^{(i)}(K \subset M) r^i &\mbox{if}\; r \geq 0\\
&=\sum_{i=1}^n Y_{n,0}^{(i)-}(K \subset M) r^i &\mbox{if}\; r \leq 0
\end{array}$$
where 
$$Y_{n,0}^{(n)-}=Y_{n,0}^{(n)}=\frac{(-1)^n}{n!}Z_n\left([M;L_n]\right)$$ is given by Theorem~\ref{thmfboun}.
Now, we prove that the two polynomial expressions, the one for $r>0$ and the one for $r<0$, coincide. Applying the above result to $M(K;\frac{1}{r_0})$ with $r_0 <-n$,
implies that for any $r\geq r_0$, 
$$Z_n(M(K;\frac{1}{r}))-Z_n(M(K;\frac{1}{r_0}))=\sum_{i=1}^n Y_{n,0}^{(i)}(K \subset M(K;\frac{1}{r_0})) (r-r_0)^i.$$ 
The above result also implies that $Z_n(M(K;\frac{1}{r}))-Z_n(M(K;\frac{1}{r_0}))$
is
$$\begin{array}{ll}\sum_{i=1}^n Y_{n,0}^{(i)}(K \subset M) r^i +Z_n(M)-Z_n(M(K;\frac{1}{r_0}))&\mbox{if}\; r \geq 0\\
\sum_{i=1}^n Y_{n,0}^{(i)-}(K \subset M) r^i +Z_n(M)-Z_n(M(K;\frac{1}{r_0}))&\mbox{if}\; r \leq 0. \end{array}$$
Therefore the coefficients of the two polynomials coincide. This proves the existence of the polynomial expression with its given leading term
for $(p,q)=(1,0)$. Applying this result in $M(K;\frac{p}{q})$ and using the fact
that 
a $\frac{p}{q+rp}$-surgery on $K$ is equivalent to a $p/q$-surgery on $K$
and a $1/r$ surgery on a parallel copy on $K$ gives a similar polynomial expression for $Z_n(M(K;\frac{p}{q+rp}))-Z_n(M(K;\frac{p}{q}))$
with the same leading coefficient since, according to Theorem~\ref{thmfboun},  $Z_n([M;L_n])=Z_n([M(K;\frac{p}{q});L_n])$.
Now, up to polynomials in $r$ of degree less than $(n-1)$,
$$Z_n(M(K;\frac{p}{q+rp}))-Z_n(M)
=\frac{(-1)^n}{n!}Z_n\left([M;L_n]\right) r^n $$ $$+
\left(\frac{n(1-n)}{2}\frac{(-1)^n}{n!}Z_n\left([M;L_n]\right) + \frac{(-1)^{n-1}}{(n-1)!}Z_n\left([M(K;\frac{p}{q});L_{n-1}]\right)\right) r^{n-1}$$
$$=Y_{n,0}^{(n)}(K \subset M) (r^n + \frac{n q}{p} r^{n-1}) 
+ Y_{n,q/p}^{(n-1)}(K \subset M) r^{n-1}.$$
Thus $$Y_{n,q/p}^{(n-1)}(K \subset M) + n  \frac{(-1)^nq}{n!p}Z_n\left([M;L_n]\right)$$
$$=\frac{(-1)^{n-1}}{(n-1)!}Z_n\left([M(K;\frac{p}{q});L_{n-1}]\right)
+\frac{1}{2}\frac{(-1)^{n-1}(n-1)}{(n-1)!}Z_n\left([M;L_n]\right).$$
Then $$Y_{n,q/p}^{(n-1)}(K \subset M) - Y_{n,0}^{(n-1)}(K \subset M)$$
$$= \frac{(-1)^{n-1}}{(n-1)!} \left(\frac{q}{p} Z_n\left([M;L_n]\right)
+Z_n\left([M(K;\frac{p}{q});L_{n-1}]-[M;L_{n-1}]\right)\right)$$
where $$Z_n\left([M(K;\frac{p}{q});L_{n-1}]-[M \sharp L(p,-q);L_{n-1}]\right)=-\frac{q}{p}Z_n([M;L_n])$$
by Theorem~\ref{thmfboun}, and by additivity of $p^c(Z_n)=Z_n^c$ under connected sum, since $n\geq 2$,
$$Z_n^c([M \sharp L(p,-q);L_{n-1}])=Z_n^c[M ;L_{n-1}].$$
Therefore $Y_{n,q/p}^{(n-1)^c}(K \subset M) = Y_{n,0}^{(n-1)^c}(K \subset M)$.\\
The behaviour of $Y_{n,q/p}^{(i)}(K \subset M)$ under an orientation change of $M$ comes from the fact that $Z_n(-M)=(-1)^nZ_n(M)$, and the other assertions are easy to observe.
\eop

\begin{remark}
It is easy to see that  $\langle \langle\; \bigsqcup_{i \in \{1,\dots,n\}}  I(F^i) \;\rangle \rangle$ is an invariant of the knot.
First, it does not depend on the symplectic bases chosen for the Seifert surfaces because $H_1(F)$ may be identified to $H_1(F)^{\ast}$ via $(x \mapsto \langle x,. \rangle)$, and therefore
the tensor $\left(\sum_i x_i \otimes y_i-\sum_i y_i \otimes x_i\right)$ may be identified with the intersection form
of the surface that lives in $H_1(F)^{\ast} \otimes H_1(F)^{\ast}$.
Now, $\langle \langle \bigsqcup_{i \in \{1,\dots,n\}}  I(F^i) \;\;\rangle \rangle$ is invariant under the addition of a hollow handle. (See \cite[p. 27]{go} or \cite{KeL}
 for a reference for the fact that for two Seifert surfaces of a knot $K$, there exists a third Seifert surface of $K$ that is obtained from the two former ones by adding hollow handles.) Indeed let $m$ be a meridian of a one-handle whose boundary is the union of the hollow handle and two disks, and let $\ell$ be a dual curve for it with respect to the intersection form of the stabilized surface $F$. Since the innermost copy of $m$ does not link any curve of the other copies of $F$, the pair $(m,\ell)$ does not contribute to the pairing. Now, the next innermost meridian does not link any other curve either...
In such a way, it is easily seen that the pairs $(m,\ell)$ can be forgotten and this shows that  $\langle \langle \;\bigsqcup_{i \in \{1,\dots,n\}}  I(F^i) \;\rangle \rangle$ is invariant under a stabilization of $F$ by addition of a hollow handle.
\end{remark}

\noindent{\sc Proof of Proposition~\ref{propzsing}:}
Let $K_0=K_{\emptyset}$ be the positive desingularisation of $K^s$. Let $U^{(k)}$ be the trivial link that bounds a disjoint union of disks $D_i$ such that each $D_i$ meets $K^s$
exactly at one double point, and $\partial D_i$ does not algebraically link $K_0$, so that
each desingularisation of $K^s$ is obtained from $K_0$ by surgery on a subset
of $L=\{(U_i;-1)\}_{i \in \{1,\dots,k\}}$.
Then $$\sum_{i=0}^n Y_{n,q/p}^{(i)}(K^s \subset M) (r+\frac{q}{p})^i=Z_n([M(K_0;\frac{p}{q+rp});L]).
$$

Each $U_i$ bounds a genus one surface $\Sigma_i$ in $M \setminus K_0$
obtained from $D_i$ by tubing $K_0$, say in the $K^{\prime}_i$ part, where we fix the choice of the $K^{\prime}_i$  so that for any pair $\{i,j\}$, $K^{\prime}_i \cap K^{\prime}_j$ is connected.

Let us prove that such a choice is indeed possible for the $K^{\prime}_i$.
Fixing the choice of $K^{\prime}_i$ amounts to choosing an interval of the circle between the two preimages of the double point $i$. If some of the two possible intervals for a double point $i$ does not contain a pair of preimages for another double point,
pick such an interval. In the next steps, if some of the two intervals for a double point $i$ only contains pairs of preimages for another double point together with their associate already chosen intervals, then pick such an interval. It is easy to see that this process will stop when all the $K^{\prime}_i$ are chosen so that for any pair $\{i,j\}$, $K^{\prime}_i \cap K^{\prime}_j$ is connected.

Now, assume that the diameters of the tubes are all constant and different and that the tube for $U_j$
is thinner than the tube for $U_i$, if $K^{\prime}_j$ contains the two preimages of the double point $i$. 
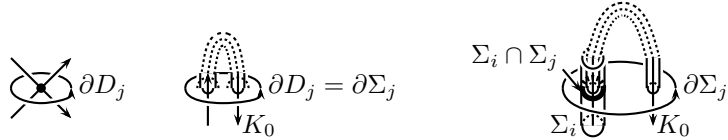
\begin{center}
\begin{pspicture}[.2](-.4,-.4)(1.4,1.4) 
\psecurve(.25,.05)(-.15,.25)(.25,.45)(.65,.25)(.25,.05)
\psline[border=1pt]{->}(-0.15,-0.15)(0.65,0.65)
\psline[border=1pt]{->}(-0.15,0.65)(0.65,-0.15)
\psdot(0.25,.25)
\psecurve[border=1pt]{->}(.25,.45)(-.15,.25)(.25,.05)(.65,.25)(.25,.45)
\rput[l](.75,.25){$\partial D_j$} \end{pspicture}
\begin{pspicture}[.2](-.9,-.4)(2,1.4) 
\psecurve(.25,.05)(-.25,.25)(.25,.45)(.75,.25)(.25,.05)
\psline[border=1pt,doubleline=true,doublesep=5pt](0.05,0.25)(0.05,0.45)
\psline[border=1pt]{->}(0.05,-0.25)(0.05,0.45)
\psecurve[linestyle=dashed,dash=1pt 1pt, doubleline=true,doublesep=4pt](0.05,0.05)(0.05,0.45)(.25,.9)(0.45,0.45)(0.45,0.05)
\psecurve[linestyle=dashed,dash=1pt 1pt](0.05,0.05)(0.05,0.45)(.25,.9)(0.45,0.45)(0.45,0.05)
\psline[border=1pt,doubleline=true,doublesep=4pt](0.45,0.45)(0.45,0.25)
\psline[border=1pt]{->}(0.45,0.45)(0.45,-0.25)
\psecurve[border=1pt]{->}(.25,.45)(-.25,.25)(.25,.05)(.75,.25)(.25,.45)
\rput[l](.85,.25){$\partial D_j=\partial \Sigma_j$} 
\rput[l](.5,-.25){$K_0$}
\psarc(0.05,0.25){2pt}{-180}{0}
\psarc(0.45,0.25){2pt}{-180}{0}
\psarc[linestyle=dashed,dash=1pt 1pt](0.05,0.25){4pt}{0}{180}
\psarc[linestyle=dashed,dash=1pt 1pt](0.45,0.25){4pt}{0}{180}
\end{pspicture}
\begin{pspicture}[.2](-3,-.4)(1.5,1.4) 
\psecurve(.45,-.1)(-.35,.25)(.45,.6)(1.15,.25)(.45,-.1)
\psarc[linestyle=dashed,dash=1pt 1pt](0.05,0.65){4pt}{0}{180}
\psline[border=1pt,doubleline=true,doublesep=8pt](.05,-0.25)(0.05,0.65)
\psline[border=1pt,doubleline=true,doublesep=4pt](.05,.25)(0.05,0.65)
\psline[border=1pt](0.05,-0.25)(0.05,0.65)
\psecurve[linestyle=dashed,dash=1pt 1pt,doubleline=true,doublesep=4pt](0.05,0.05)(0.05,0.65)(.45,1.3)(0.85,0.65)(0.85,0.05)
\psecurve[linestyle=dashed,dash=1pt 1pt](0.05,0.05)(0.05,0.65)(.45,1.3)(0.85,0.65)(0.85,0.05)
\psline[border=1pt,doubleline=true,doublesep=4pt](0.85,0.65)(0.85,0.25)
\psline[border=1pt]{->}(0.85,0.65)(0.85,-0.25)
\psecurve[border=1pt]{->}(.45,.6)(-.35,.25)(.45,-.1)(1.15,.25)(.45,.6)
\rput[l](1.25,.25){$\partial \Sigma_j$} 
\rput[l](.9,-.25){$K_0$}
\psarc(0.05,0.25){2pt}{-180}{0}
\psarc[linewidth=2pt](0.05,0.25){4pt}{-180}{0}
\psarc(0.85,0.25){2pt}{-180}{0}
\psarc[linestyle=dashed,dash=1pt 1pt](0.05,0.25){2pt}{0}{180}
\psarc[linewidth=2pt,linestyle=dashed,dash=1pt 1pt](0.05,0.25){4pt}{0}{180}
\psarc[linestyle=dashed,dash=1pt 1pt](0.85,0.25){2pt}{0}{180}
\psarc(0.05,-0.25){4pt}{-180}{0}
\psarc[border=1pt](0.05,0.65){4pt}{-180}{0}
\psarc[linestyle=dashed,dash=1pt 1pt](0.05,-0.25){4pt}{0}{180}
\rput[r](-.15,-.25){$\Sigma_i$}
\rput[r](-.4,.65){$\Sigma_i \cap \Sigma_j$}
\psline[border=1pt]{->}(-.35,.5)(-.1,.2)
\end{pspicture}
\end{center}
Then $\Sigma_i \cap \Sigma_j$
is empty if the pair $(D_i \cap K_0)$ does not link the pair $(D_j \cap K_0)$, and it is a meridian curve of $K_0$ in $D_j$ otherwise. Therefore, the $\mu$-invariants of the three-component sublinks of $L$ in $M(K_0;\frac{p}{q+rp})$ are zero 
and Theorem~\ref{thmfasmu} can be used to compute
$Z_n([M(K_0;\frac{p}{q+rp});L])$.

In particular, if $k>2n$,
$Z_n([M(K_0;\frac{p}{q+rp});L])=0$. Since the linking numbers between two intersection curves will be $\pm \frac{q+rp}{p}$ or zero, if $k=2n$,
$Z_n([M(K_0;\frac{p}{q+rp});L])$ is a monomial in $\left(\frac{q+rp}{p}\right)^n$.
\eop

\noindent{\sc Proof of Proposition~\ref{propzsingun}:}
In this case, the link $U^{(k)}$ of the previous proof is a boundary link in 
$M(K_0;\frac{p}{q+rp})$ because the produced genus one surfaces are disjoint.
Then Theorem~\ref{thmfboun} can be applied.
 It implies the first part of the proposition. Use bases $(m_i, \ell_i)$ for the Seifert surfaces where $m_i$ is a $0$-framed meridian of $K_0$, and $\ell_i$ is 
a curve along the tube of $\Sigma_i$ and $D_i$ that is homotopic to
$K^{s\prime}_i$, and that does not link $K_0$. In $M(K_0;\frac{p}{q+rp})$, the linking number of two meridians is
$\pm \frac{q+rp}{p}$, the linking number of a meridian and a longitude is $0$ or $\pm 1$ while the linking number of two longitudes is their linking number in $M$. Note that if one tube for $\Sigma_j$ goes inside another one for $\Sigma_i$ (if $K^{s,i}_j \neq K^{s\prime}_j$), and if $K^{\prime}_j$ is the positive desingularisation of $K^{s\prime}_j$ then $lk(\ell_i,\ell_j)=lk(\ell_i,K^{\prime}_j)=-lk(\ell_i,K_0-K^{\prime}_j)=-\ell_{ij}(K^s)$. There are at most $n$ pairings of meridians. Furthermore, since there is at least one innermost meridian that cannot be paired with a longitude, there is at least one pairing of meridians. Now, the number of pairs of meridians coincides with the number of pairs of longitudes in a pairing. 
\eop

As an example, we compute $Y_{2}^{(i)c}(K^s)$ where $K^s$ is a singular knot with two unlinked double points.

\begin{proposition}
\label{propvarztwo}
Let $K^s$ be a singular knot with two unlinked double points.
$$\sum_{I \subset \{1,2\}}(-1)^{\sharp I}Z^c_2(M(K_I;\frac{1}{r}))=$$
$$\frac{1}{4}\left(\left( 5 \ell_{12}(K^s)^2 + 2\ell_{11}(K^s)\ell_{22}(K^s) \right)r^2
-\ell_{12}(K^s)r\right)\tatata.$$
\end{proposition}
\bp Use the strategy and the notation of the proofs of the two previous propositions. Choose Seifert surfaces of the two knots of the crossing changes
with disjoint tubes whose longitudes $\ell_1$ and $\ell_2$ are homotopic to 
$K^{s,2}_1$ and $K^{s,1}_2$, respectively, so that
$$lk(\ell_1,\ell_2)=lk(\ell_1,\ell_2^+)=lk(\ell_1^+,\ell_2)=\ell_{12}(K^s)$$
$$lk(m_1,m_2)=-r=lk(m_i,m_i^+)$$
$$lk(\ell_i,\ell_i^+)=-\ell_{ii}(K^s)$$
$$lk(m_i,\ell_i^+)=lk(m_1,\ell_2)=lk(m_2,\ell_1)=0$$
$$lk(m_i^+,\ell_i)=1.$$
Then, according to Theorem~\ref{thmfboun},
$$\sum_{I \subset \{1,2\}}(-1)^{\sharp I}Z^c_2(M(K_I;\frac{1}{r}))=\frac{1}{4} p^c\left(\langle \langle \begin{pspicture}[0.5](-.2,-.5)(1.2,1.2) 
\psline{-}(0.05,0)(.35,.15)
\psline{*-*}(0.05,0)(.05,.8) 
\psline{-}(0.05,0)(.35,-.15)
\psline{-}(0.05,.8)(.35,.95)
\psline{-}(0.05,.8)(.35,.65)
\rput[lt](.5,-.15){$m_1$}
\rput[l](.5,.15){$\ell_1$}
\rput[l](.5,.65){$\ell_1^+$}
\rput[lb](.5,.95){$m_1^+$}
\end{pspicture}
\begin{pspicture}[0.5](-.2,-.5)(1.2,1.2) 
\psline{-}(0.05,0)(.35,.15)
\psline{*-*}(0.05,0)(.05,.8) 
\psline{-}(0.05,0)(.35,-.15)
\psline{-}(0.05,.8)(.35,.95)
\psline{-}(0.05,.8)(.35,.65)
\rput[lt](.5,-.15){$m_2$}
\rput[l](.5,.15){$\ell_2$}
\rput[l](.5,.65){$\ell_2^+$}
\rput[lb](.5,.95){$m_2^+$}
\end{pspicture} \rangle \rangle \right).$$

Note that $m_1$ and $m_2$ must be paired to another meridian.
Then the right-hand side of the equality can be rewritten as
 $$\frac{r^2}{4} p^c\left(\langle \langle  \begin{pspicture}[0.3](-.35,0)(2.1,1) 
\pscurve(1.85,0)(2.1,.4)(1.85,.8)
\psline{*-*}(1.85,0)(1.85,.8) 
\pscurve(0.05,0)(-.2,.4)(.05,.8)
\psline{*-*}(0.05,0)(.05,.8) 
\psline{-}(0.05,0)(.35,.15)
\psline{-}(1.85,.8)(1.55,.65)
\psline{-}(1.85,0)(1.55,.15)
\psline{-}(0.05,.8)(.35,.65)
\rput[l](.5,.15){$\ell_1$}
\rput[l](.5,.65){$\ell_1^+$}
\rput[r](1.5,.15){$\ell_2$}
\rput[r](1.5,.65){$\ell_2^+$}
\end{pspicture}
+\begin{pspicture}[0.3](-.2,0)(2,1) 
\pscurve(1.85,0)(.95,-.15)(.05,0)
\psline{*-*}(1.85,0)(1.85,.8) 
\pscurve(0.05,.8)(.95,.95)(1.85,.8)
\psline{*-*}(0.05,0)(.05,.8) 
\psline{-}(0.05,0)(.35,.15)
\psline{-}(1.85,.8)(1.55,.65)
\psline{-}(1.85,0)(1.55,.15)
\psline{-}(0.05,.8)(.35,.65)
\rput[l](.5,.15){$\ell_1$}
\rput[l](.5,.65){$\ell_1^+$}
\rput[r](1.5,.15){$\ell_2$}
\rput[r](1.5,.65){$\ell_2^+$}
\end{pspicture}
+ \begin{pspicture}[0.3](-.2,0)(2,1) 
\pscurve(1.85,0)(.95,-.15)(.05,0)
\psline{*-*}(1.85,0)(1.85,.8) 
\pscurve(0.05,.8)(.95,.95)(1.85,.8)
\psline{*-*}(0.05,0)(.05,.8) 
\psline{-}(0.05,0)(.35,.15)
\psline{-}(1.85,.8)(1.55,.65)
\psline{-}(1.85,0)(1.55,.15)
\psline{-}(0.05,.8)(.35,.65)
\rput[l](.5,.15){$\ell_1$}
\rput[l](.5,.65){$\ell_1^+$}
\rput[r](1.5,.15){$\ell_2^+$}
\rput[r](1.5,.65){$\ell_2$}
\end{pspicture}
 \rangle \rangle \right) -\frac{r}{4}\ell_{12}(K^s) \tatata.$$
Indeed, either two pairs of meridians are paired together.  
This leads to the
quadratic contribution in $r$ above, or there is one pair of meridians, it is necessarily $(m_1,m_2)$ and in this case $m_1^+$ must be paired with $\ell_1$ and $m_2^+$ must be paired with $\ell_2$. Then $\ell_1^+$ and  $\ell_2^+$ must be paired together, and this provides the linear contribution above.
\eop

\section{Computation of the Casson-Walker knot invariant}
\setcounter{equation}{0}
\label{secproofcasone}

Let $K$ be an order $O_K$ knot in a rational homology sphere $M$. Let $\widetilde{M\setminus K}$ be the infinite cyclic covering of $M\setminus K$. 
Denote the action of the homotopy class of the meridian of $K$ on $H_1(\widetilde{M\setminus K};\QQ)$ as the multiplication by $t$ so that a generator of $H_1(M\setminus K)/\mbox{Torsion}$ acts as the multiplication by 
$t^{1/O_K}$.
As in \cite[Chapter 2]{pup}, define the {\em Alexander polynomial\/} $\Delta(K)$ of $K$ as 
the order of the $\QQ[t^{\pm 1/O_K}]$-module $H_1(\widetilde{M\setminus K};\QQ)$ normalized so that 
$$\Delta(K)(1)=|\mbox{Torsion}(H_1(M\setminus K))|=\frac{|H_1(M)|}{O_K}\; \mbox{and} \;\Delta(K)(t^{1/O_K})=\Delta(K)(t^{-1/O_K}).$$

Then the formula of \cite[p 12-13]{pup} implies the following lemma.

\begin{lemma}
\label{lemcaschirknot}
For any knot $K$ such that $lk(K,K) \in \ZZ$ in a rational homology sphere $M$,
for any pair $(p,q)$ of coprime integers such that $q\neq 0$.
$$\lambda(M(K;p/q))-\lambda(M)=\frac{q}{p}\left(\frac{O_K}{|H_1(M)|}\frac{\Delta^{\prime\prime}(K)(1)}{2} - \frac{1}{24O_K^2} + \frac{1}{24}\right) +\lambda(L(p,-q)).$$
\end{lemma}
\bp Recall that $\lambda(M)=\frac{\overline{\lambda}(M)}{|H_1(M)|}$ where $|H_1(M)|$ is the cardinality of $H_1(M;\ZZ)$ and $\overline{\lambda}$ is the extension of $|H_1(M)|\lambda$ to oriented closed $3$-manifolds that is denoted by $\lambda$ in \cite{pup}.
For any knot $K$ in a rational homology sphere $M$, according to \cite[1.4.8,T2]{pup}, for $q>0$,
$$\lambda(M(K;p/q))-\lambda(M)=$$
$$\frac{q}{p} \left( \frac{O_K}{|H_1(M)|}\frac{\Delta^{\prime\prime}(K)(1)}{2} - \frac{1}{24O_K^2} -\frac{p^2+1}{24q^2}\right) +\frac{\mbox{sign}(pq)}{8} + \frac{s(p-qlk(K,K),q)}{2}$$
where the Dedekind sum $s(p-qlk(K,K),q)$ is defined in \cite{rg} (and in \cite[1.4.5]{pup}). This formula makes clear that
$$\lambda(M(K;p/q))-\lambda(M)=\frac{q}{p}\left(\frac{O_K}{|H_1(M)|}\frac{\Delta^{\prime\prime}(K)(1)}{2} - \frac{1}{24O_K^2} + \frac{1}{24}\right) +f(p,q)$$
for some $f(p,q)$ that depends neither on the knot $K$ with self-linking number $0$ nor on its ambient manifold $M$.
Applying this formula to the trivial knot $U$ of $S^3$ concludes the proof of the lemma.
\eop

We now express $\Delta(K)$ from the Seifert form of a Seifert surface for $K$
in the following proposition.

\begin{proposition}
\label{propalexseif}
Let $K$ be a knot of order $d$,
with self-linking number $(-a/b) \in \QQ/\ZZ$, where $a$ and $b$ are coprime integers, in a rational homology sphere $M$. 
Let $N(K)$ be a tubular neighborhood of $K$.
Let $\Sigma$ be a surface in $M$ whose boundary is made of $(d/b)$ parallel copies of a primitive curve of $\partial N(K)$. Let ${\cal B}_s$ be a symplectic basis for $H_1(\Sigma)/H_1(\partial \Sigma)$, and let $$\Delta_{\Sigma}(\tau)=\mbox{det} [lk(\tau^{1/2} b^{\prime +} - \tau^{-1/2}b^{\prime -},b)]_{(b, b^{\prime}) \in {\cal B}_s^2}$$
where $b^{\prime +}$ (resp. $b^{\prime -}$) is a representative
of $b^{\prime}$ pushed away from $\Sigma$ in the direction of the positive 
(resp. negative) normal direction to $\Sigma$.
Then $$\frac{d}{|H_1(M)|}\Delta(K)=\Delta_{\Sigma}(t^{1/d})\frac{b}{d}\frac{t^{\frac{1}{2b}}-t^{-\frac{1}{2b}}}{t^{\frac{1}{2d}}-t^{-\frac{1}{2d}}}.$$
\end{proposition}
\bp
First assume that the self-linking number of $K$ is zero.
Let $N(K)$ be a tubular neighborhood of $K$.
There exists a genus $g$ surface $\Sigma$ in $M$ whose boundary is made of $d$ parallel copies of $K$.
Consider a collar $\Sigma \times [-1,1]$ in $M$ such that 
$$(\Sigma \times [-1,1]) \cap N(K) =\partial \Sigma \times [-1,1].$$

Let $Y=M \setminus (N(K) \cup \Sigma \times ]-1,1[)$.

The infinite cyclic covering $\tilde{X}$ of $M \setminus N(K)$ can be seen
as
$$\left(\coprod_{k\in \ZZ}h^k(Y)
\coprod \coprod_{k\in \ZZ}
h^k(\Sigma \times [-1,1])\right)/\cong$$ 
where $h$ is a generator of the group of automorphisms of the covering $\tilde{X}$
and $\cong$ provides the following identifications. 
$$h^k\left((\sigma \in \Sigma, 1) \in Y\right) \cong h^k\left((\sigma \in \Sigma, 1) \in (\Sigma \times [-1,1])\right)$$ $$h^k\left((\sigma \in \Sigma, -1) \in Y \right) \cong h^{k+1}\left((\sigma \in \Sigma, -1) \in (\Sigma \times [-1,1])\right).$$

Then it is easy to see that, if the action of $h$ on $H_1(\tilde{X};\QQ)$ is denoted as a multiplication by $\tau$, 
$$H_1(\tilde{X};\QQ)=\frac{H_1(Y;\QQ) \otimes \QQ[\tau,\tau^{-1}]}{(\oplus_{b\in {\cal B}}(\tau b^+-b^-) \QQ)\otimes \QQ[\tau,\tau^{-1}]},$$
as a $\QQ[\tau,\tau^{-1}]$-module,
where ${\cal B}$ is a basis of $H_1(\Sigma)$ and, for $b \in {\cal B}$, $b^+$ (resp. $b^-$) denotes the class of $b$ in $H_1(\Sigma \times \{1\})$ (resp. in $H_1(\Sigma \times \{-1\})$).

In particular, if ${\cal C}$ is a basis of $H^1(Y;\QQ)$, then 
$$\Delta(K)(\tau=t^{1/d})=\mbox{det}\left[\tau^{1/2} c(b^+) -\tau^{-1/2} c(b^-)\right]_{(c,b)\in {\cal C} \times {\cal B}}$$
up to a multiplication by a unit of $\QQ[\tau,\tau^{-1}]$.\\

\noindent{\em Computation of $H^1(Y;\QQ)$.\/}\\
Let $Z=M \setminus (\Sigma \times ]-1,1[)$.

The collar $\Sigma \times [-1,1]$ is a genus $(2g+d-1)$-handlebody whose
$H_1$ has a basis ${\cal B}$ made of the classes $\ell_1, \ell_2, \dots, \ell_{d-1}$ of $(d-1)$ boundary components of $\Sigma$, and a symplectic basis ${\cal B}_s$ for $H_1(\Sigma)/H_1(\partial \Sigma)$.
Therefore, $Z$ has the rational homology of a genus $(2g+d-1)$-handlebody and $H^1(Z;\QQ)$ is freely generated by the linking numbers
with the elements of ${\cal B}$.

Use the following exact sequence to compute $H^1(Y;\QQ)$

$$H^1(Z,Y;\QQ) \hookrightarrow H^1(Z;\QQ) \rightarrow H^1(Y;\QQ) \rightarrow
H^2(Z,Y;\QQ) \rightarrow 0.$$

The pair $(Z,Y)$ has the homology of the pair 
$(N(K), \partial N(K) \setminus (\partial \Sigma \times [-1,1]))$
where $\partial N(K) \setminus (\partial \Sigma \times [-1,1])$ is a disjoint union of $d$ annuli $A(\ell^{++}_i)$ whose cores are parallels $\ell^{++}_1, \ell^{++}_2, \dots, \ell^{++}_{d}$ of $K$, and such that $$\partial A(\ell^{++}_i) = \ell^+_{i} -\ell^-_{i+1}$$
(where $\ell^-_{d+1}=\ell^-_1$). In particular,
$$\begin{array}{lll}H_j(Z,Y)&=0\;&\mbox{if}\; j \neq 1,2\\
&=\oplus_{i=2}^d \ZZ c_i\;&\mbox{if}\; j =1\\
&=\oplus_{i=2}^d \ZZ B_i\;&\mbox{if}\; j =2 \end{array}$$
where $c_i$ is the class of a path from $\ell^{++}_{1}$ to $\ell^{++}_{i}$ in $N(K)$, and
$B_i$ is the class of an annulus whose boundary is $(\ell^{++}_{i}-\ell^{++}_{1})$.

The image of $H^1(Z,Y;\QQ) \hookrightarrow H^1(Z;\QQ)$ is freely generated by the algebraic intersections  $\langle .,-A(\ell^{++}_i)\rangle=lk(.,\ell^-_{i+1} -\ell^+_{i})$ for $i\in {2, \dots,d}$. 

For $i\geq 2$, consider a curve $e_i$ that goes from $\ell_{i}$ to $\ell_{i+1}$ in $\Sigma$ 
and that avoids the chosen geometric symplectic basis of $H_1(\Sigma)/H_1(\partial \Sigma)$, and consider a closed loop $\mu_i$ in $N(K) \cup \Sigma \times [-1,1]$ that equals
$e_i$ outside $N(K)$. 

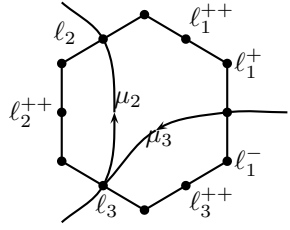
\begin{center}
\begin{pspicture}[0.4](-1,-.3)(4.2,3) 
\psline{-*}(2.2,.7)(2.2,1.35)
\psline{-*}(2.2,1.35)(2.2,2)
\psline{-*}(2.2,2)(1.65,2.325)
\psline{-*}(1.65,2.325)(1.1,2.65)
\psline{-*}(1.1,2.65)(.55,2.325)
\psline{-*}(.55,2.325)(0,2)
\psline{-*}(0,2)(0,1.35) 
\psline{-*}(0,1.35)(0,.7)
\psline{-*}(0,.7)(.55,.375)
\psline{-*}(.55,.375)(1.1,.05)
\psline{-*}(1.1,.05)(1.65,.375)
\psline{-*}(1.65,.375)(2.2,.7)
\rput[l](2.3,.7){$\ell_1^-$}
\rput[l](2.3,2){$\ell_1^+$}
\rput[lb](1.7,2.4){$\ell_1^{++}$}
\rput[r](-.1,1.35){$\ell_2^{++}$}
\rput[lt](1.7,.35){$\ell_3^{++}$}
\rput[r](.2,2.4){$\ell_2$}
\rput[t](.6,.25){$\ell_3$}
\psecurve{->}(-1.1,-1.088)(0,-.113)(.55,.375)(.7,1.35)(.55,2.325)
\psecurve(.55,.375)(.7,1.35)(.55,2.325)(0,2.813)(-1.1,3.788)
\rput[lb](.7,1.4){$\mu_2$}
\psecurve{->}(4.1,1.35)(3,1.35)(2.2,1.35)(1.25,1.11)(.55,.375)
\psecurve(2.2,1.35)(1.25,1.11)(.55,.375)(0,-.113)
\rput[t](1.3,1.1){$\mu_3$}
\end{pspicture}
\end{center}

Then $lk(\partial B_i=\ell^{++}_{i}-\ell^{++}_{1},\mu_j)=\delta_{ij}$.
Therefore the map $H^1(Y;\QQ) \rightarrow
H^2(Z,Y;\QQ)$ admits a section whose image is $\oplus_{i=2}^d\QQ lk(.,\mu_i)$.

Thus $$H^1(Y;\QQ)=\oplus_{b \in {\cal B}_s}\QQ lk(.,b) \oplus \oplus_{i=2}^d \QQ lk(.,\mu_i).$$
Since $lk(\ell_i^{\pm},b)=0$ for any $b \in {\cal B}_s$, up to units of $\QQ[t^{\pm 1/O_K}]$,
$$\Delta(K)=\Delta_{\Sigma}(\tau) \Delta(d)$$
with 
$$\Delta(d)=
\mbox{det} [lk(\tau^{1/2} \ell_i^+ - \tau^{-1/2}\ell_i^-,\mu_j)]_{(i,j) \in \{2,\dots,d\}^2}$$
where $\ell_i^+=\ell^{++}_{i}$ and  $\ell_i^-= \ell^{++}_{i-1}$.
\begin{sublemma}
$$\Delta(d)=\frac{\tau^{d/2}-\tau^{-d/2}}{d(\tau^{1/2}-\tau^{-1/2})}.$$
\end{sublemma}
\noindent{\sc Proof of the sublemma:}
By pushing $\mu_j$ along
the negative normal of the Seifert surface of $\ell^{++}_{1}$ (or $K$)
we see that $lk(\ell^{++}_{1},\mu_j)=-\frac{1}{d}$.

Set $z=\tau^{1/2}-\tau^{-1/2}$ and $\rho=\tau^{1/2}$,
 $\Delta(d)$ is the determinant of the following matrix
$[\Delta_{ij}]_{(i,j) \in \{2,\dots,d\}^2}$
where $$\Delta_{2j}=lk(\rho \ell_2^{++} - \rho^{-1}\ell_1^{++},\mu_j)= \rho\delta_{2j}-\frac{z}{d},$$
and for $i>2$ 
$$\Delta_{ij}=lk(\rho (\ell_i^{++}-\ell_2^{++}) - \rho^{-1}(\ell^{++}_{i-1}-\ell_1^{++}),\mu_j)=\rho(\delta_{ij}-\delta_{2j}) -\rho^{-1}\delta_{(i-1)j},$$
that is for $d=5$,
$$[\Delta_{ij}]=\left[\begin{array}{cccc}
\rho-\frac{z}{d}&-\frac{z}{d}&-\frac{z}{d}&-\frac{z}{d}\\
-\rho-\rho^{-1}&\rho&0&0\\
-\rho&-\rho^{-1}&\rho&0\\
-\rho&0&-\rho^{-1}&\rho\end{array}
\right],$$
$\Delta(2)=\frac{\rho+\rho^{-1}}{2}$ and $\Delta(3)=\frac{\tau+\tau^{-1}+1}{3}.$

In general the development with respect to the first column gives that
$$\Delta(d)=(\rho-\frac{z}{d}) \rho^{(d-2)} 
-\frac{z}{d}(\rho+\rho^{-1}) \left( \sum_{j=3}^d \rho^{(d-j)-(j-3)} \right)
-\frac{z}{d} \rho \sum_{i=4}^d \rho^{i-3} \sum_{j=i}^d \rho^{(d-j)-(j-i)}$$
where
$$\rho^{-1}\left( \sum_{j=3}^d \rho^{(d-j)-(j-3)} \right)=
\sum_{j=3}^d \rho^{(d-2j+2)}=\rho^{(2-d)} + \rho^{(4-d)} + \dots + \rho^{(d-4)}
.$$ Thus,
$$\Delta(d)=\rho^{(d-1)}- \frac{z}{d} \rho^{(d-2)} -\frac{ \rho^{(d-3)} - \rho^{(1-d)}}{d}
-\frac{z}{d} \rho^{-2} \sum_{i=3}^d \sum_{j=i}^d  \rho^{(d+2i-2j)}.$$
$$\sum_{i=3}^d \sum_{j=i}^d \rho^{(d+2i-2j)}=(d-2)\rho^{d} + (d-3)\rho^{d-2}+ (d-4)\rho^{d-4}+ \dots 
+ \rho^{(6-d)}$$
$$z \sum_{i=3}^d \sum_{j=i}^d \rho^{(d+2i-2j)}
=(d-2)\rho^{(d+1)} - \rho^{(d-1)} -\rho^{(d-3)} -\dots - \rho^{(5-d)}$$
$$d\Delta(d)=d\rho^{(d-1)}-\rho^{(d-1)} + \rho^{(1-d)} 
-(d-2)\rho^{(d-1)} + \rho^{(d-3)} +\rho^{(d-5)} +\dots + \rho^{(3-d)}$$
$$=\rho^{(d-1)}+ \rho^{(d-3)} +\rho^{(d-5)} +\dots + \rho^{(3-d)} + \rho^{(1-d)}=\frac{\tau^{d/2}-\tau^{-d/2}}{(\tau^{1/2}-\tau^{-1/2})}.$$
\eop

Back to the proof of Proposition~\ref{propalexseif},
since
$\Delta(K)(t=\tau^d)(1)=\frac{|H_1(M)|}{d}$,
$$\Delta(K)(t)=\frac{|H_1(M)|}{d}\frac{t^{1/2}-t^{-1/2}}{d(t^{1/(2d)}-t^{-1/(2d)})}\Delta_{\Sigma}(\tau).$$
and Proposition~\ref{propalexseif} is proved in the self-linking number $0$ case.
Let us now deduce the general case from this case.
Let $K$ be a knot with order $d$ and with self-linking number $(-a/b)$
where $a$ and $b$ are coprime. Let $m$ be a meridian of $K$, there exist a parallel $L$ of $K$ and a surface $\Sigma$ in $M \setminus K$ whose boundary is made of $(d/b)$ parallel copies of $am+bL$.
Then there exists a primitive curve $m_J$ such that $\langle m_J,am+bL \rangle=1$. Let $J$ be the knot with meridian $m_J$ and with complement $M\setminus K$.
This knot has order $(d/b)$ and self-linking number $0$. Its Alexander
polynomial is then given by the proposition. Furthermore, since it satisfies $\Delta(J)(1)=|\mbox{Torsion}(H_1(M \setminus K))|=\Delta(K)(1)$, $\Delta(J)(t_J=\tau^{d/b})=\Delta(K)(t_K=\tau^{d})$.
Then $\Delta(K)(t_K)=\Delta(J)(t_J=t_K^{1/b})$ and we are done.
\eop

Proposition~\ref{propalexseif} implies the following
lemma that together with Lemma~\ref{lemcaschirknot} proves Proposition~\ref{propcasknot} for $n=1$. 
We use the notation of Section~\ref{seccasstate}.
\begin{lemma}
\label{lemalexorder}
Under the assumptions of Proposition~\ref{propalexseif}, $$\frac{d}{|H_1(M)|}\frac{\Delta^{\prime\prime}(K)(1)}{2}=\frac{\langle \langle I(\Sigma) \rangle \rangle_{W_1}}{2d^2} + \frac{1}{24b^2} - \frac{1}{24d^2}.$$
\end{lemma}
\bp
First note that when $K$ is null-homologous, $O_K=d=b=1$.
Then since $\lambda= W_1 \circ Z_1$, Lemma~\ref{lemcaschirknot} together with Theorem~\ref{thmfboun} together imply that
$$\frac{1}{|H_1(M)|}\frac{\Delta^{\prime\prime}(K)(1)}{2}=\frac{\langle \langle I(\Sigma) \rangle \rangle_{W_1}}{2}.$$
Therefore, according to Proposition~\ref{propalexseif} (that is well-known in this case),
$$\frac{\Delta^{\prime\prime}_{\Sigma}(1)}{2}=\frac{\langle \langle I(\Sigma) \rangle \rangle_{W_1}}{2}.$$
Then since $\Delta_{\Sigma}(t)=\Delta_{\Sigma}(t^{-1})$, 
$$\Delta_{\Sigma}(\exp(u))=1 +\frac{\langle \langle I(\Sigma) \rangle \rangle_{W_1}}{2} u^2 +O(4)$$
where $O(4)$ stands for an element of $u^4\QQ[[u]]$,
and this formula remains true for any $\Sigma$
as in the statement of Proposition~\ref{propalexseif}.
Since
$$\exp(u)^{\frac{1}{2d}}-\exp(u)^{-\frac{1}{2d}}=\frac{u}{d}(1+\frac{u^2}{24d^2} +O(4)),$$ it is easy to conclude.
\eop

Now that Proposition~\ref{propcasknot} is shown for $n=1$,
let us prove it for $n=2$.
By the formula that is recalled in the beginning of the proof of Lemma~\ref{lemcaschirknot},
$$\lambda(M(K;p/q))-\lambda(M)= \frac{q}{p}\frac{O_K}{|H_1(M)|}\frac{\Delta^{\prime\prime}(K)(1)}{2} +f(p,q,lk(K,K),O_K)$$
for some $f(p,q,lk(K,K),O_K)$ that only depends on $p$, $q$, $lk(K,K)$, $O_K$, and that therefore does not change under surgery on a knot $K_2$ that does not link $K$ algebraically,
so that 
$$\sum_{I \subset \{1,2\}} (-1)^{\sharp I} \lambda \left(M_{(K_i;p_i/q_i)_{i \in I}}\sharp \sharp_{j \in \{1,2\} \setminus I} L(p_j,-q_j)\right)= \sum_{I \subset \{1,2\}} (-1)^{\sharp I} \lambda \left(M_{(K_i;p_i/q_i)_{i \in I}}\right)$$
$$=\frac{q_1}{p_1}\left(\frac{O_{K_1}\Delta^{\prime\prime}(K_1\subset M(K_2;p_2/q_2))(1)}{2|H_1(M(K_2;p_2/q_2))|}-\frac{O_{K_1}\Delta^{\prime\prime}(K_1\subset M)(1)}{2|H_1(M)|}\right)$$
$$=\frac{q_1}{2p_1O_{K_1}^2}\left(\langle \langle I(\Sigma_1) \subset M(K_2;p_2/q_2) \rangle \rangle_{W_1}
-\langle \langle I(\Sigma_1) \subset M \rangle \rangle_{W_1} \right)$$
according to Lemma~\ref{lemalexorder}.
Therefore, Proposition~\ref{propcasknot} for $n=2$ follows from the following lemma.
\begin{lemma} 
\label{lemvarlambdaprime}
Under the assumptions of Proposition~\ref{propcasknot},
$$\langle \langle I(\Sigma_1) \subset M(K_2;p/q) \rangle \rangle_{W_1}
-\langle \langle I(\Sigma_1) \subset M \rangle \rangle_{W_1} 
=-\frac{2q}{d_2^2p} lk\left(\Sigma_1 \cap \Sigma_2, (\Sigma_1 \cap \Sigma_2)_{\parallel}\right).$$
\end{lemma}
\bp
Let $(x_i,y_i)_{i \in \{1,\dots,g\}}$ be a symplectic basis for $H_1(\Sigma_1)/H_1(\partial \Sigma_1)$.
Because of the variation of linking numbers after surgery recalled in Lemma~\ref{lemvarlk}, the variation of the expression of $\langle \langle I(\Sigma_1) \rangle \rangle_{W_1}$ given before Proposition~\ref{propcasknot} reads
$$\langle \langle I(\Sigma_1) \subset M(K_2;p/q) \rangle \rangle_{W_1}
-\langle \langle I(\Sigma_1) \subset M \rangle \rangle_{W_1} =$$
$$2 \frac{q^2}{p^2} \sum_{(j,k) \in \{1,2,\dots, g\}^2} lk(x_j,K_2)lk(K_2,x_k^+)lk(y_j,K_2)lk(K_2,y_k^+)$$
$$-2 \frac{q^2}{p^2} \sum_{(j,k) \in \{1,2,\dots, g\}^2}lk(x_j,K_2)lk(K_2,y_k^+)lk(y_j,K_2)lk(K_2,x_k^+)$$
$$-2\frac{q}{p} \sum_{(j,k) \in \{1,2,\dots, g\}^2}\left(lk(x_j,K_2)lk(K_2,x_k^+)lk(y_j,y_k^+)-lk(x_j,K_2)lk(K_2,y_k^+)lk(y_j,x_k^+)\right)$$
$$-2\frac{q}{p}\sum_{(j,k) \in \{1,2,\dots, g\}^2}\left(lk(x_j,x_k^+)lk(y_j,K_2)lk(K_2,y_k^+)-lk(x_j,y_k^+)lk(y_j,K_2)lk(K_2,x_k^+)\right)$$
where the quadratic part in $q/p$ is obviously zero.
On the other hand, when $c \in H_1(\Sigma_1)$, 
$$\langle c, \Sigma_1 \cap \Sigma_2 \rangle_{\Sigma_1} =d_2 lk(c,K_2).$$
Therefore in $H_1(\Sigma_1)$, 
$$\Sigma_1 \cap \Sigma_2=d_2\sum_{i=1}^g (lk(x_i,K_2)y_i -lk(y_i,K_2)x_i)$$
and 
$$lk(\Sigma_1 \cap \Sigma_2, (\Sigma_1 \cap \Sigma_2)^{+})$$
$$=d_2^2\sum_{(j,k) \in \{1,2,\dots, g\}^2}
lk\left(lk(x_j,K_2)y_j -lk(y_j,K_2)x_j, lk(x_k,K_2)y_k^+ -lk(y_k,K_2)x_k^+\right).$$
\eop

Then Proposition~\ref{propcasknot} is proved for $n=2$. Since $$lk_{M(K_3;p_3/q_3)}(\Sigma_1 \cap \Sigma_2, (\Sigma_1 \cap \Sigma_2)_{\parallel})
-lk_M(\Sigma_1 \cap \Sigma_2, (\Sigma_1 \cap \Sigma_2)_{\parallel})
=-\frac{q_3}{p_3}lk_M(\Sigma_1 \cap \Sigma_2,K_3)^2$$
$$=-\frac{q_3}{d_3^2p_3}\langle \Sigma_1, \Sigma_2, \Sigma_3 \rangle^2 $$ this in turn implies Proposition~\ref{propcasknot}   for $n=3$.
Now, since $lk(\Sigma_1 \cap \Sigma_2,K_3)$ does not vary under a surgery on a knot that does not link $K_1$, $K_2$ and $K_3$ algebraically, Proposition~\ref{propcasknot}
is also true for $n\geq 4$ and hence for all $n$.
\eop

\noindent{\sc Proof of Proposition~\ref{propcasknottwo}:}
Use that $\lambda^{\prime}(K^s)=\lambda^{\prime}(U,K^-)$ where $U$ is a trivial knot that surrounds the crossing change. (See the proofs of Propositions~\ref{propzsing} and~\ref{propvarztwo} in Section~\ref{secproofpol}.)
\eop

\section{Proofs of the statements on $\lambda_2$ and $w_3$}
\setcounter{equation}{0}
\label{secprooflambdatwo}

Theorem~\ref{thmpol} guarantees the existence of a polynomial surgery formula
$$\lambda_2(M(K;p/q))-\lambda_2(M)=\lambda_2^{\prime \prime}(K) (q/p)^2 + w_3(K)(q/p) +C(K;q/p) + \lambda_2(L(p;-q))$$
where $C(K;q/p)$ only depends on $q/p$ mod $\ZZ$ and $C(U;q/p)=0$.
Since $Z_2^c(-M)=Z_2^c(M)$, $w_3(K \subset M)=-w_3(K \subset (-M))$.

Furthermore, according to Proposition~\ref{propvarztwo}, if 
$K^s$ is a singular link with two unlinked double points, then $w_3(K^s)=-\frac{\ell_{12}(K^s)}{2}$ and $C(K^s;q/p)=0$.

The only unproved assertion of Theorem~\ref{thmw3} is that the knot invariants $C(K;q/p)$ read $c(q/p)\lambda^{\prime}(K)$ for knots that bound a surface whose $H_1$ vanishes in $H_1(M)$. The proof of this assertion will be given in this section.

Also note that for any knot $K$ in a rational homology sphere $M$,
$w_3(K \subset M)= w_3(K \subset M \sharp N)$
and $C(K \subset M;q/p)= C(K \subset M \sharp N;q/p)$.

Let $K^s$ be a singular
knot with one double point in a rational homology sphere.
Let $K^+$ and $K^-$ be its two desingularisations, and let $K^{\prime}$ and $K^{\prime \prime}$
be the two knots obtained from $K^s$ by smoothing the double point. Assume that
$K^{\prime}$ and $K^{\prime \prime}$ are null-homologous, set
 $$f(K^s)=\frac{\lambda^{\prime}(K^{\prime}) + \lambda^{\prime}(K^{\prime \prime})}{2} -\frac{\lambda^{\prime}(K^+) + \lambda^{\prime}(K^-)+lk^2(K^{\prime},K^{\prime \prime})}{4}.$$
 
Note that $f(K^s \subset M)=f(K^s \subset M \sharp N)$.

In order to prove Proposition~\ref{propvarwthree}, we shall successively prove the following lemmas. The two last ones Lemmas~\ref{lemredcalcpar} and \ref{lemcalcpar} obviously imply Proposition~\ref{propvarwthree}.

\begin{lemma}
\label{lemnotvarunderxing}Let $K^s$ be a singular
knot with one double point in a rational homology sphere. The invariants $C(K^s;q/p)$
 and $(w_3-f)(K^s)$ do not vary under a surgery on a knot that is null-homologous in the complement of $K^s$.
\end{lemma}

\begin{lemma}
\label{lemsurpresgr}
Let $\Gamma$ be a non-necessarily connected graph in a rational homology sphere
$M$,
such that every loop of $\Gamma$ is null-homologous in $M$.
Then there exist a graph $\Gamma_0$ in $S^3$, an algebraically split (rationally) framed link $L$ in $S^3$ whose components are null-homologous in $S^3 \setminus \Gamma_0$,
 and a rational homology sphere $N$, such that $(S^3(L),\Gamma_0)= (M,\Gamma) \sharp N$.
\end{lemma}

\begin{lemma}
\label{lemredcalcpar}
Let $K^s_n$ be the following singular knot
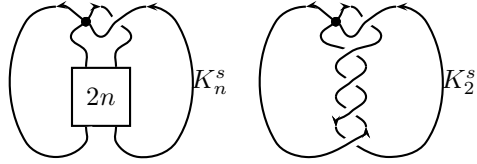
\begin{center}
\begin{pspicture}[.2](-.2,-.1)(3,2.2) 
\psecurve{->}(1.2,.6)(1.2,.4)(1.2,.2)(1.5,0)(2,.2)(2,1.8)(1.6,2)(1.2,1.8)
\psecurve(.8,1)(.8,1.2)(.8,1.4)(.6,1.6)(.8,1.8)(1,2)
\psecurve{->}(.6,1.6)(.8,1.8)(1,2)(1.2,1.8)
\psecurve[border=2pt](.8,1.8)(1,2)(1.2,1.8)(1.4,1.6)(1.2,1.4)(1.2,1.2)(1.2,1)
\psframe(.6,.4)(1.4,1.2)
\rput(1,.8){$2n$}
\psecurve{*->}(1,1.6)(.8,1.8)(.4,2)(.2,1.8)
\psecurve[border=2pt](.8,1.8)(.4,2)(0,1.8)(0,.2)(.4,0)(.8,.2)(.8,.4)(.8,.6)
\psecurve[border=2pt]{-*}(2,1.8)(1.6,2)(1.2,1.8)(1,1.6)(.8,1.8)(.4,2)
\rput[l](2.2,1){$K_n^s$}
\end{pspicture}
\begin{pspicture}[.2](-.2,-.1)(2.2,2.2) 
\psecurve{->}(1,.6)(.8,.4)(1,.2)(1.5,0)(2,.2)(2,1.8)(1.6,2)(1.2,1.8)
\psecurve(1,1)(1.2,1.2)(1,1.4)(.6,1.6)(.8,1.8)(1,2)
\psecurve(1,.2)(1.2,.4)(1,.6)(.8,.8)(1,1)
\psecurve{->}(.6,1.6)(.8,1.8)(1,2)(1.2,1.8)
\psecurve[border=2pt]{->}(.8,1.8)(1,2)(1.2,1.8)(1.4,1.6)(1,1.4)(.8,1.2)(1,1)(1.2,.8)(1,.6)(.8,.4)(1,.2)
\psecurve{*->}(1,1.6)(.8,1.8)(.4,2)(.2,1.8)
\psecurve[border=2pt]{->}(.8,1.8)(.4,2)(0,1.8)(0,.2)(.4,0)(1,.2)(1.2,.4)(1,.6)
\psecurve[border=2pt](1,.6)(.8,.8)(1,1)(1.2,1.2)(1,1.4)
\psecurve[border=2pt]{-*}(2,1.8)(1.6,2)(1.2,1.8)(1,1.6)(.8,1.8)(.4,2)
\rput[l](2.2,1){$K_2^s$}
\end{pspicture}
\end{center}
where \begin{pspicture}[.4](-.1,-.1)(.6,.6) 
\psframe(0,0)(.5,.5)
\rput(.25,.25){$2n$}
\end{pspicture} represents $|n|$ vertical juxtapositions of the motive 
\begin{pspicture}[.2](.65,.2)(1.35,.9)
\psecurve(1,.1)(1.15,.25)(1,.4)(.85,.55)(1,.7)
\psecurve[border=1pt](1,1)(.85,.85)(1,.7)(1.15,.55)(1,.4)(.85,.25)(1,.1)
\psecurve[border=1pt](1,.4)(.85,.55)(1,.7)(1.15,.85)(1,1)
\end{pspicture} if $n>0$ and
$|n|$ vertical juxtapositions of the motive 
\begin{pspicture}[.2](.65,.2)(1.35,.9)
\psecurve(1,.1)(.85,.25)(1,.4)(1.15,.55)(1,.7)
\psecurve[border=1pt](1,1)(1.15,.85)(1,.7)(.85,.55)(1,.4)(1.15,.25)(1,.1)
\psecurve[border=1pt](1,.4)(1.15,.55)(1,.7)(.85,.85)(1,1)
\end{pspicture}if $n<0$.
Then for any singular knot $K^s$ with one double point $p$, such that the two knots $K^{\prime}$ and $K^{\prime \prime}$ obtained from $K^s$ by smoothing $p$ are null-homologous,
 $$(w_3-f)(K^s)=(w_3-f)(K^s_{-lk(K^{\prime},K^{\prime \prime})}).$$
\end{lemma}
\begin{lemma}
\label{lemcalcpar}
For all $n \in \ZZ$, $(w_3-f)(K_n^s)=0$.
\end{lemma}

We shall next prove the following proposition that generalizes a Casson lemma 
from integral to rational homology spheres.
\begin{proposition}
\label{propgmcar}
Let $C$ be a real-valued invariant of null-homologous knots in rational homology spheres
such that 
\begin{itemize}
\item $C(K \subset M)= C(K \subset M \sharp N)$, 
\item $C(U)=0$, 
\item $C(K)$ does not vary under a surgery on a knot $J$ such that $(J,K)$ is a boundary link,
\item if $K^s$ is a singular knot with one double point, $C(K^s)$ does not vary under surgery on a knot that is null-homologous in the complement of $K^s$. 
\end{itemize}
Then there exists $c \in \RR$ such that 
\begin{itemize}
\item if $K^s$ is singular knot with one double point $p$, such that the two knots $K^{\prime}$ and $K^{\prime \prime}$ obtained from $K^s$ by smoothing $p$ are null-homologous,
then $C(K^s)=c lk(K^{\prime},K^{\prime \prime})$, and,
\item if $K$ bounds a surface whose $H_1$ maps to zero in $H_1(M)$, $C(K)=c \lambda^{\prime}(K)$.
\end{itemize}
\end{proposition}

Since the $C(.;p/q)$ satisfy the hypotheses of the proposition above
(thanks to Theorem~\ref{thmfboun} for the hypothesis on boundary links), this proposition will be sufficient to conclude the proof of Theorem~\ref{thmw3}.
\eop

Let us now prove all the lemmas and the proposition.

\noindent{\sc Proof of Lemma~\ref{lemnotvarunderxing}:}
Let $J$ be a null-homologous knot unlinked with $K^{\prime}$ and $K$. Let $F_J$ be a Seifert surface for $J$ that does not meet $K^s$, and let $(m, \ell)$ be the usual basis of the genus one surface obtained by tubing a trivial knot $V$
surrounding the double point of $K^s$, $m$ is a meridian of $K^-$, $\ell$ is homotopic to $K^{\prime}$ and $lk(\ell,K^-)=0$. By Theorem~\ref{thmfboun},
$$Z_2^c(M(J;\frac{p_J}{q_J})(K^s;\frac{p}{q} ) )- Z_2^c(M(K^s;\frac{p}{q})  )$$
$$= \frac{q_J}{4p_J} p^c \left( \langle \langle \begin{pspicture}[0.4](-.5,-.3)(1.3,1) 
\psline{-}(0,0.05)(.15,.35)
\psline{*-*}(0,0.05)(.8,.05) 
\psline{-}(0,0.05)(-.15,.35)
\psline{-}(.8,0.05)(.95,.35)
\psline{-}(.8,0.05)(.65,.35)
\rput[br](-.15,.5){$\ell$}
\rput[b](.15,.5){$m$}
\rput[b](.65,.5){$m^+$}
\rput[lb](.95,.5){$\ell^+$}
\end{pspicture} I(F_J) \subset M(K^-;\frac{p}{q})\;\rangle \rangle \right).$$

Since, according to Lemma~\ref{lemvarlk}, the pairing of $m$ and a curve $c$ in the contraction above will give rise to the coefficient
$(-q/p)lk(K,c)=-rlk(K,c)$, $C(K^s;q/p)$ does not vary under a $(p_J/q_J)$-surgery on $J$.

$$w_3(K^s \subset M(J;\frac{p_J}{q_J})) - w_3(K^s \subset M) = $$ $$\left(\frac{\partial}{\partial r}\right)_{r=0}  W_2\left( Z_2^c(K^s \subset M(J;\frac{p_J}{q_J})) - Z_2^c(K^s \subset M)\right)$$
where $m$ must be paired either with $m^+$ or with $I(F_J)$, and in the latter
case $m^+$ must be paired with $\ell$ in order to lead to a linear contribution in $r$.
$$w_3(K^s \subset M(J;\frac{p_J}{q_J})) - w_3(K^s \subset M) = $$
$$ =- \frac{q_J}{4p_J} \langle \langle \begin{pspicture}[0.4](-.5,-.1)(1,.8)  
\psline{*-*}(0,0.05)(.4,.05) 
\pscurve(0,0.05)(.2,.2)(.4,0.05)
\psline{-}(0,0.05)(-.15,.35)
\psline{-}(.4,0.05)(.55,.35)
\rput[b](-.15,.4){$\ell$}
\rput[b](.6,.4){$\ell^+$}
\end{pspicture} I(F_J) \subset M \rangle \rangle_{W_2} + \frac{q_j}{4p_J}
\langle \langle \begin{pspicture}[0.4](-.5,-.1)(1,.8)  
\psline{*-*}(0,0.05)(.4,.05) 
\pscurve(0,0.05)(.2,.2)(.4,0.05)
\psline{-}(0,0.05)(-.15,.35)
\psline{-}(.4,0.05)(.55,.35)
\rput[b](-.15,.4){$K$}
\rput[b](.6,.5){$\ell^+$}
\end{pspicture} I(F_J) \subset M \rangle \rangle_{W_2}.$$

Since $K=K^{\prime \prime} + \ell$, as far as the connected pairing with $I(F_J)$ is concerned,
$$\begin{pspicture}[0.15](-.6,-.15)(1.15,.6) 
\pscurve{*-*}(0,0.05)(.2,-.1)(.4,0.05) 
\pscurve(0,0.05)(.2,.2)(.4,0.05)
\psline{-}(0,0.05)(-.2,0.05)
\psline{-}(.4,0.05)(.6,0.05)
\rput[rb](-.25,-.05){$K$}
\rput[lb](.65,-.05){$K^+$} \end{pspicture}
=\begin{pspicture}[0.15](-.5,-.15)(1.1,.6) 
\pscurve{*-*}(0,0.05)(.2,-.1)(.4,0.05) 
\pscurve(0,0.05)(.2,.2)(.4,0.05)
\psline{-}(0,0.05)(-.2,0.05)
\psline{-}(.4,0.05)(.6,0.05)
\rput[rb](-.25,-.05){$\ell$}
\rput[lb](.65,-.05){$\ell^+$}
\end{pspicture} + 
\begin{pspicture}[0.15](-.8,-.15)(1.3,.6) 
\pscurve{*-*}(0,0.05)(.2,-.1)(.4,0.05) 
\pscurve(0,0.05)(.2,.2)(.4,0.05)
\psline{-}(0,0.05)(-.2,0.05)
\psline{-}(.4,0.05)(.6,0.05)
\rput[rb](-.25,-.05){$K^{\prime \prime}$}
\rput[lb](.65,-.05){$K^{\prime \prime +}$}
\end{pspicture}
+2 
\begin{pspicture}[0.15](-.5,-.15)(1.1,.6) 
\pscurve{*-*}(0,0.05)(.2,-.1)(.4,0.05) 
\pscurve(0,0.05)(.2,.2)(.4,0.05)
\psline{-}(0,0.05)(-.2,0.05)
\psline{-}(.4,0.05)(.6,0.05)
\rput[rb](-.25,-.05){$\ell$}
\rput[lb](.65,-.05){$K^{\prime \prime}$}
\end{pspicture}
$$
and
$$\begin{pspicture}[0.15](-.6,-.15)(1.1,.6) 
\pscurve{*-*}(0,0.05)(.2,-.1)(.4,0.05) 
\pscurve(0,0.05)(.2,.2)(.4,0.05)
\psline{-}(0,0.05)(-.2,0.05)
\psline{-}(.4,0.05)(.6,0.05)
\rput[rb](-.25,-.05){$K$}
\rput[lb](.65,-.05){$\ell^+$} \end{pspicture}
=
\begin{pspicture}[0.15](-.5,-.15)(1.1,.6) 
\pscurve{*-*}(0,0.05)(.2,-.1)(.4,0.05) 
\pscurve(0,0.05)(.2,.2)(.4,0.05)
\psline{-}(0,0.05)(-.2,0.05)
\psline{-}(.4,0.05)(.6,0.05)
\rput[rb](-.25,-.05){$\ell$}
\rput[lb](.65,-.05){$\ell^+$}
\end{pspicture} + 
\begin{pspicture}[0.15](-.8,-.15)(1.2,.6) 
\pscurve{*-*}(0,0.05)(.2,-.1)(.4,0.05) 
\pscurve(0,0.05)(.2,.2)(.4,0.05)
\psline{-}(0,0.05)(-.2,0.05)
\psline{-}(.4,0.05)(.6,0.05)
\rput[rb](-.25,-.05){$K^{\prime \prime}$}
\rput[lb](.65,-.05){$\ell^{+}$}
\end{pspicture}.
$$
Therefore,
$$\begin{pspicture}[0.15](-.6,-.15)(1.1,.6) 
\pscurve{*-*}(0,0.05)(.2,-.1)(.4,0.05) 
\pscurve(0,0.05)(.2,.2)(.4,0.05)
\psline{-}(0,0.05)(-.2,0.05)
\psline{-}(.4,0.05)(.6,0.05)
\rput[rb](-.25,-.05){$K$}
\rput[lb](.65,-.05){$\ell^+$} \end{pspicture}
= \frac12
\begin{pspicture}[0.15](-.5,-.15)(1.1,.6) 
\pscurve{*-*}(0,0.05)(.2,-.1)(.4,0.05) 
\pscurve(0,0.05)(.2,.2)(.4,0.05)
\psline{-}(0,0.05)(-.2,0.05)
\psline{-}(.4,0.05)(.6,0.05)
\rput[rb](-.25,-.05){$\ell$}
\rput[lb](.65,-.05){$\ell^+$}
\end{pspicture} - \frac12 
\begin{pspicture}[0.15](-.8,-.15)(1.3,.6) 
\pscurve{*-*}(0,0.05)(.2,-.1)(.4,0.05) 
\pscurve(0,0.05)(.2,.2)(.4,0.05)
\psline{-}(0,0.05)(-.2,0.05)
\psline{-}(.4,0.05)(.6,0.05)
\rput[rb](-.25,-.05){$K^{\prime \prime}$}
\rput[lb](.65,-.05){$K^{\prime \prime +}$}
\end{pspicture}
+ \frac12
\begin{pspicture}[0.15](-.6,-.15)(1.15,.6) 
\pscurve{*-*}(0,0.05)(.2,-.1)(.4,0.05) 
\pscurve(0,0.05)(.2,.2)(.4,0.05)
\psline{-}(0,0.05)(-.2,0.05)
\psline{-}(.4,0.05)(.6,0.05)
\rput[rb](-.25,-.05){$K$}
\rput[lb](.65,-.05){$K^+$} \end{pspicture}.
$$

$$w_3(K^s \subset M(J;\frac{p_J}{q_J})) - w_3(K^s \subset M) = $$
$$ - \frac{q_J}{8p_J} \langle \langle \begin{pspicture}[0.4](-.5,-.2)(.8,1)  
\psline{*-*}(0,0.05)(.4,.05) 
\pscurve(0,0.05)(.2,.2)(.4,0.05)
\psline{-}(0,0.05)(-.15,.35)
\psline{-}(.4,0.05)(.55,.35)
\rput[b](-.15,.5){$\ell$}
\rput[b](.55,.5){$\ell^+$}
\end{pspicture} I(F_J) \rangle \rangle_{W_2} - \frac{q_J}{8p_J}
\langle \langle \begin{pspicture}[0.4](-.5,-.2)(.8,1)  
\psline{*-*}(0,0.05)(.4,.05) 
\pscurve(0,0.05)(.2,.2)(.4,0.05)
\psline{-}(0,0.05)(-.15,.35)
\psline{-}(.4,0.05)(.55,.35)
\rput[b](-.15,.5){$K^{\prime \prime}$}
\rput[b](.55,.5){$K^{\prime \prime +}$}
\end{pspicture} I(F_J) \rangle \rangle_{W_2}
+ \frac{q_J}{8p_J}
\langle \langle \begin{pspicture}[0.4](-.5,-.2)(1,1)  
\psline{*-*}(0,0.05)(.4,.05) 
\pscurve(0,0.05)(.2,.2)(.4,0.05)
\psline{-}(0,0.05)(-.15,.35)
\psline{-}(.4,0.05)(.55,.35)
\rput[b](-.15,.5){$K$}
\rput[b](.55,.5){$K^+$}
\end{pspicture} I(F_J) \rangle \rangle_{W_2}.$$

Thus, according to Proposition~\ref{propcasknotvar}, since
$$\langle \langle \begin{pspicture}[0.3](-.7,-.3)(1.4,.6) 
\pscurve{*-*}(0,0.05)(.2,-.1)(.4,0.05) 
\pscurve(0,0.05)(.2,.2)(.4,0.05)
\psline{-}(0,0.05)(-.2,0.05)
\psline{-}(.4,0.05)(.6,0.05)
\rput[rb](-.25,-.05){$K$}
\rput[lb](.65,-.05){$K^+$} \end{pspicture} I(F_J) \;\rangle \rangle_{W_2}= \langle \langle \begin{pspicture}[0.3](-.6,-.3)(1.2,.6)  
\psline{-}(.35,0.05)(-.05,.05)
\rput[rb](-.1,0){$K$}
\rput[lb](.4,0){$K^+$}
\end{pspicture} I(F_J) \;\rangle \rangle_{W_1},$$
$$w_3(K^s \subset M(J;\frac{p_J}{q_J})) - w_3(K^s \subset M) $$
$$ =\frac{q_J}{2p_J} \left(   \lambda^{\prime}(J,K^{\prime}) +
\lambda^{\prime}(J,K^{\prime \prime})- \lambda^{\prime}(J,K)\right)$$
$$ =\frac{q_J}{p_J} \left(   \frac{\lambda^{\prime}(J,K^{\prime})}{2} +
\frac{\lambda^{\prime}(J,K^{\prime \prime})}{2}- \frac{\lambda^{\prime}(J,K^+)}{4}
- \frac{\lambda^{\prime}(J,K^-)}{4}\right)$$
$$=f(K^s \subset M(J;\frac{p_J}{q_J})) - f(K^s \subset M).$$
\eop

\noindent{\sc Proof of Lemma~\ref{lemsurpresgr}:}
After a possible connected sum with some lens spaces,
the $\QQ/\ZZ$--valued linking form of $M$ is diagonal \cite{wall}, and the generators of $H_1(M;\ZZ)$ can be represented by a link $L$ of algebraically unlinked curves $K_i$ that do not link $\Gamma$, algebraically.
Then for each $K_i$, there exists a surface $\Sigma_i$ in the exterior $(M \setminus \mbox{Int}{N}(L))$ of $L$  whose boundary is a connected essential curve of $\partial N(K_i)$, and that does not meet $\Gamma$.
Thus, $H^1(M \setminus \mbox{Int}{N}(L);\ZZ)$ is freely generated by the algebraic intersections with the $\Sigma_i$, and there exists
a surgery on $L$ that transforms $M$ into a homology sphere $H$.
The manifold $H$ can in turn be transformed into $S^3$ by surgery on a boundary link of $H$ bounding a disjoint union $F_H$ of surfaces in $H$ that can be assumed to be disjoint from the first surgery link and from the image of $\Gamma$ in $H$. This proves the lemma.
\eop

\noindent{\sc Proof of Lemma~\ref{lemredcalcpar}:}
Apply Lemma~\ref{lemsurpresgr} to $\Gamma=K^s$, then $K^{s,0}=\Gamma_0$. Note that $lk(K^{\prime}_0,K^{\prime\prime}_0)=lk(K^{\prime},K^{\prime\prime})$. 
Recall that $(w_3-f)(K^s \subset M)=(w_3-f)(K^s \subset M \sharp N)$. Thanks to Lemma~\ref{lemnotvarunderxing}, $(w_3-f)(K^s \subset M \sharp N)=(w_3-f)(K^{s,0} \subset S^3)$.
Now that the proof has been reduced to the case where $M=S^3$, recall that a
crossing change on $K^{\prime}$ or $K^{\prime \prime}$ may be realized by a surgery on a knot satisfying the hypotheses of Lemma~\ref{lemnotvarunderxing}, that changes neither $lk(K^{\prime},K^{\prime\prime})$ nor $(w_3-f)(K^s)$.
Unknotting $K^{\prime}$ first by crossing changes and next unknotting the parts of $K^{\prime \prime}$ between two consecutive intersection points 
with the disk bounded by $K^{\prime}$ transforms $K^s$ into $K^s_{-lk(K^{\prime},K^{\prime\prime})}$.
\eop

\noindent{\sc Proof of Lemma~\ref{lemcalcpar}:}
By the crossing change formula of Proposition~\ref{propcasknottwo}, $\lambda^{\prime}(K^+_{n})-\lambda^{\prime}(K^+_{n-1})=-1$, and
$\lambda^{\prime}(K^+_{n})=-n$.
Since $K^-_n$, $K^{\prime}_n$ and $K^{\prime \prime}_n$ are trivial, $f(K^s_{n})=-\frac{n(n-1)}{4}$.

On the other hand, since $w_3(K^-_n)=0$, $w_3(K^s_n)=w_3(K^+_n)$.
 The unlinked double crossing change formula of Theorem~\ref{thmw3} implies that
$$w_3(K^+_{n+2})-2w_3(K^+_{n+1}) +w_3(K^+_{n})=-\frac{1}{2}$$
Since $K^+_0$ is trivial, $w_3(K^+_{0})=0$, and since $K^+_1$ is the figure-eight knot that coincides with its mirror image, $w_3(K^+_{1})=0$, too. 
Then $w_3(K^s_n)=w_3(K^+_{n})=-\frac{n(n-1)}{4}$.
\eop

\noindent{\sc Proof of Proposition~\ref{propgmcar}:} 
Let $K^s$ be as in the hypotheses of Proposition~\ref{propvarwthree}. The proof of Lemma~\ref{lemredcalcpar} shows that
$C(K^s)=C(K^s_{-lk(K^{\prime},K^{\prime \prime})})$. Since $C(U)=0$, 
$C(K^s_n)=C(K^+_n)$. 
Since the hypotheses of the proposition imply that $C$ maps singular knots of $S^3$ with two unlinked double points to $0$,
$C(K^+_{n+2})-2C(K^+_{n+1}) +C(K^+_{n})=0$, and $C(K^+_{n})$ is affine with respect to $n$. Since $C(K^+_{0})=0$, $C(K^+_{n})$ is linear. Then there exists $c$ such that $C(K^s)=c lk(K^{\prime},K^{\prime \prime})$.

Let $K$ be a knot that bounds a Seifert surface $\Sigma$ whose $H_1$ maps to zero in $H_1(M)$.
Applying Lemma~\ref{lemsurpresgr} to the one-skeleton of $\Sigma$ allows us to reduce the proof that $C(K)=c \lambda^{\prime}(K)$ to the case of knots in $S^3$, thanks to the hypothesis on boundary links. Then this case
is easily proved with the crossing change formula.
\eop

\noindent {\sc Proof of Proposition~\ref{propgenusone}:}
Consider the genus one surface $\Sigma$ in $H$ and its symplectic basis $(a,b)$ below.
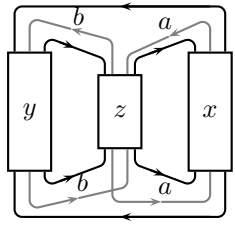
\begin{center}
\begin{pspicture}[0.4](-.2,-1.5)(3.2,1.5)
\psline[linearc=.1](2.9,.8)(2.9,1.4)(.1,1.4)(.1,.8)
\psline[linearc=.1](.5,.8)(.5,1)(1.3,.7)(1.3,.5)
\psline[linearc=.1](1.7,.5)(1.7,.7)(2.5,1)(2.5,.8)
\psline{->}(2.7,1.4)(1.5,1.4)
\psline{->}(.7,.925)(.9,.85)
\psline{->}(1.9,.775)(2.1,.85)
\psframe(1.2,-.5)(1.8,.5)
\psframe(0,-.8)(.6,.8)
\psframe(2.4,-.8)(3,.8)
\rput(.3,0){$y$}
\rput(1.5,0){$z$}
\rput(2.7,0){$x$}
\psline[linearc=.1](2.9,-.8)(2.9,-1.4)(.1,-1.4)(.1,-.8)
\psline[linearc=.1](.5,-.8)(.5,-1)(1.3,-.7)(1.3,-.5)
\psline[linearc=.1](1.7,-.5)(1.7,-.7)(2.5,-1)(2.5,-.8)
\psline{->}(2.7,-1.4)(1.5,-1.4)
\psline{->}(.7,-.925)(.9,-.85)
\psline{->}(1.9,-.775)(2.1,-.85)
\rput[b](2.1,-1.1){$a$}
\psline[linearc=.1,linecolor=gray]{->}(1.4,-.5)(1.4,-1.2)(2.05,-1.2)
\psline[linearc=.1,linecolor=gray](2.05,-1.2)(2.7,-1.2)(2.7,-.8)
\rput[b](2.1,1.1){$a$}
\psline[linearc=.1,linecolor=gray]{->}(2.7,.8)(2.7,1.25)(2.15,1)
\psline[linearc=.1,linecolor=gray](2.15,1)(1.6,.75)(1.6,.5)
\rput[b](1,-1.05){$b$}
\psline[linearc=.1,linecolor=gray]{->}(.3,-.8)(.3,-1.3)(.95,-1.15)
\psline[linearc=.1,linecolor=gray](.95,-1.15)(1.6,-1)(1.6,-.5)
\rput[b](.95,1.15){$b$}
\psline[linearc=.1,linecolor=gray]{->}(1.4,.5)(1.4,.95)(.85,1.1)
\psline[linearc=.1,linecolor=gray](.85,1.1)(.3,1.25)(.3,.8)
\end{pspicture}
\end{center}
$\langle a,b \rangle=1$, $lk(a,a^+)=\frac{x+z}{2}$,
$lk(b,b^+)=\frac{y+z}{2}$, $lk(a,b^+)=\frac{-1-z}{2}$, $lk(a^+,b)=\frac{1-z}
{2}$, $$\lambda^{\prime}(K(x,y,z))= \frac{(x+z)(y+z) +1  
-z^2}{4}=\frac{xy +yz +zx +1}{4}.$$
Note that $\lambda^{\prime}(\phi(X),\phi(Y)) =\lambda^{\prime}(\phi(Y),\phi(Z))  
=\lambda^{\prime}(\phi(Z),\phi(X))$.
In particular, both sides of the equality to be proved are symmetric under a cyclic permutation of  
$((X,x),(Y,y),(Z,z))$.
Using this cyclic symmetry, the formula for the pretzel knot $K(x,y,z)$ follows from the crossing change formula
starting with the trivial knot $K_{-1,1,1}$:
$$4w_3(K(x+2,y,z))-4w_3(K(x,y,z))=\lambda^{\prime}(K(x+2,y,z)) +  
\lambda^{\prime}(K(x,y,z))
+\left(\frac{y+z}{2}\right)^2.$$
$$16\left(w_3(K(x+2,y,z))-w_3(K(x,y,z))\right)=(2x+2)(y+z)+2 +4yz +y^2 +z^2.$$
$$32 w_3(K(x,y,z)) = 2x + 4xyz + xy^2 +xz^2 + x^2(y+z) +F(y,z).$$
Otherwise, the following lemma~\ref{lemvargenone} reduces the proof of Proposition~\ref{propgenusone} to the case where the knot $\phi(K(x,y,z))$ is in $S^3$, thanks to Lemma~\ref{lemsurpresgr}, and
next when the knot is a pretzel knot $K(x,y,z)$ by crossing changes on $X$ and $Y$.

\begin{lemma}
\label{lemvargenone}
Let $\phi$ be an embedding of $H$ in a rational homology sphere such that
$\phi(X)$ and $\phi(Y)$ are null homologous in the exterior of $\phi(H)$.
Let $J$ be a knot in the exterior of $\phi(H)$ that links neither $\phi(X)$ nor $\phi(Y)$, then $$w_3(\phi(K(x,y,z)) \subset M(J;p/q)) - w_3(\phi(K(x,y,z)) \subset M)$$
$$=\frac{q}{2p} \left(3 \lambda^{\prime}(\phi(X),\phi(Y),J) -x \lambda^{\prime}(\phi(X),J) - y \lambda^{\prime}(\phi(Y),J) -z \lambda^{\prime}(\phi(Z),J)\right).$$
\end{lemma}
\noindent{\sc Proof of Lemma~\ref{lemvargenone}:}
According to Theorem~\ref{thmfboun}, 
if $F_J$ is a Seifert surface of $J$ in the complement of the genus one Seifert surface $\Sigma$ of $\phi(K(x,y,z))$ in $\phi(H)$,
$$w_3(\phi(K(x,y,z)) \subset M(J;p/q)) - w_3(\phi(K(x,y,z)) \subset M)=
\frac{q}{4p} \langle \langle\; I(\Sigma) \; I(F_J)\; \rangle \rangle_{W_2}$$
where $$I(\Sigma)=\begin{pspicture}[0.4](-.5,-.3)(1.3,1) 
\psline{-}(0,0.05)(.15,.35)
\psline{*-*}(0,0.05)(.8,.05) 
\psline{-}(0,0.05)(-.15,.35)
\psline{-}(.8,0.05)(.95,.35)
\psline{-}(.8,0.05)(.65,.35)
\rput[br](-.15,.5){$a$}
\rput[b](.15,.5){$b$}
\rput[b](.65,.5){$b^+$}
\rput[lb](.95,.5){$a^+$}\end{pspicture}.$$

Write $$\frac{q}{4p} \langle \langle\; I(\Sigma) \; I(F_J)\; \rangle \rangle_{W_2}=C_A +C_B$$
where $C_A$ is the contribution of the pairings that pair
two univalent vertices of $I(\Sigma)$, and $C_B$ is the contribution of the pairings that pair all the univalent vertices of $I(\Sigma)$ to univalent vertices of $I(F_J)$.
$$C_A=\frac{q}{4p} \langle \langle\; 
\left(\frac{x+z}{2} \begin{pspicture}[0.15](-.7,-.15)(1.1,.6) 
\pscurve{*-*}(0,0.05)(.2,-.1)(.4,0.05) 
\pscurve(0,0.05)(.2,.2)(.4,0.05)
\psline{-}(0,0.05)(-.2,0.05)
\psline{-}(.4,0.05)(.6,0.05)
\rput[rb](-.25,-.05){$b$}
\rput[lb](.65,-.05){$b^+$}
\end{pspicture}\; 
+\frac{y+z}{2}\begin{pspicture}[0.15](-.7,-.15)(1.1,.6) 
\pscurve{*-*}(0,0.05)(.2,-.1)(.4,0.05) 
\pscurve(0,0.05)(.2,.2)(.4,0.05)
\psline{-}(0,0.05)(-.2,0.05)
\psline{-}(.4,0.05)(.6,0.05)
\rput[rb](-.25,-.05){$a$}
\rput[lb](.65,-.05){$a^+$}
\end{pspicture}
+z  \begin{pspicture}[0.15](-.7,-.15)(1.1,.6) 
\pscurve{*-*}(0,0.05)(.2,-.1)(.4,0.05) 
\pscurve(0,0.05)(.2,.2)(.4,0.05)
\psline{-}(0,0.05)(-.2,0.05)
\psline{-}(.4,0.05)(.6,0.05)
\rput[rb](-.25,-.05){$a$}
\rput[lb](.65,-.05){$b$}
\end{pspicture}
\right)I(F_J)\; \rangle \rangle_{W_2}.$$

From now on, we write $X$, $Y$ and $Z$ for $\phi(X)$, $\phi(Y)$ and $\phi(Z)$, respectively.
$$C_A=\frac{q}{4p} \langle \langle\; 
\left(\frac{x}{2} \begin{pspicture}[0.15](-.7,-.15)(1.1,.6) 
\pscurve{*-*}(0,0.05)(.2,-.1)(.4,0.05) 
\pscurve(0,0.05)(.2,.2)(.4,0.05)
\psline{-}(0,0.05)(-.2,0.05)
\psline{-}(.4,0.05)(.6,0.05)
\rput[rb](-.25,-.05){$X$}
\rput[lb](.65,-.05){$X^+$}
\end{pspicture}\; 
+\frac{y}{2}\begin{pspicture}[0.15](-.7,-.15)(1.1,.6) 
\pscurve{*-*}(0,0.05)(.2,-.1)(.4,0.05) 
\pscurve(0,0.05)(.2,.2)(.4,0.05)
\psline{-}(0,0.05)(-.2,0.05)
\psline{-}(.4,0.05)(.6,0.05)
\rput[rb](-.25,-.05){$Y$}
\rput[lb](.65,-.05){$Y^+$}
\end{pspicture}
+\frac{z}{2}  \begin{pspicture}[0.15](-1.8,-.15)(2.2,.6) 
\pscurve{*-*}(0,0.05)(.2,-.1)(.4,0.05) 
\pscurve(0,0.05)(.2,.2)(.4,0.05)
\psline{-}(0,0.05)(-.2,0.05)
\psline{-}(.4,0.05)(.6,0.05)
\rput[rb](-.25,-.1){$(X+Y)$}
\rput[lb](.65,-.1){$(X+Y)^+$}
\end{pspicture}
\right)I(F_J)\; \rangle \rangle_{W_2}.$$

Thus, according to Proposition~\ref{propcasknotvar}, $$C_A=-\frac{q}{p}\left(\frac{x \lambda^{\prime}(X,J)}{2} + \frac{y \lambda^{\prime}(Y,J)}{2} +\frac{z \lambda^{\prime}(Z,J)}{2}\right).$$

Let us now compute the contribution of the pairings that are bijections from 
the set of univalent vertices of $I(\Sigma)$ to the set of univalent vertices of $I(F_j)$.
For them, we may change $a$ to $Y$ and $b$ to $X$ and write
$$I(\Sigma)=\begin{pspicture}[0.4](-.5,-.3)(1.3,1) 
\psline{-}(0,0.05)(.15,.35)
\psline{*-*}(0,0.05)(.8,.05) 
\psline{-}(0,0.05)(-.15,.35)
\psline{-}(.8,0.05)(.95,.35)
\psline{-}(.8,0.05)(.65,.35)
\rput[br](-.15,.5){$X$}
\rput[b](.15,.5){$Y$}
\rput[b](.65,.5){$Y^+$}
\rput[lb](.95,.5){$X^+$}
\end{pspicture}$$
where the superscripts $+$ distinguish two copies of $X$ (or $Y$) whose linking numbers with the curves of $F_J$ are the same.

Let us compute the contribution $C_B$ of the pairings that are bijections from 
the set of univalent vertices of $I(\Sigma)$ to the set of univalent vertices of some
$$I(c,d,e,f)=\begin{pspicture}[0.4](-.5,-.3)(1.3,1) 
\psline{-}(0,0.05)(.15,.35)
\psline{*-*}(0,0.05)(.8,.05) 
\psline{-}(0,0.05)(-.15,.35)
\psline{-}(.8,0.05)(.95,.35)
\psline{-}(.8,0.05)(.65,.35)
\rput[br](-.15,.5){$c$}
\rput[b](.15,.5){$d$}
\rput[b](.65,.5){$e$}
\rput[l](.95,.6){$f$}
\end{pspicture}$$
to $\langle \langle\; I(c,d,e,f)\; I(\Sigma)  \;\rangle \rangle_{W_2}$.

Note the symmetry under the exchange of the pair $(X,Y)$ with the pair $(X^+,Y^+)$.

The contribution of the pairings that pair $c$ and $d$ to
$X$ and $X^+$ is
$$lk(c,X)lk(d,X)\langle \langle\; \begin{pspicture}[0.4](-.9,-.3)(1.3,1.2) 
\psline{-*}(0,0.05)(.2,.45)
\psline{-}(.2,.45)(.35,.75)
\psline{*-*}(0,0.05)(.8,.05) 
\psline(.2,.45)(-.2,.45)
\psline{-}(-.2,.45)(-.35,.75)
\psline{-*}(0,0.05)(-.2,.45)
\psline{-}(.8,0.05)(.95,.35)
\psline{-}(.8,0.05)(.65,.35)
\rput[br](-.4,.75){$Y$}
\rput[bl](.4,.75){$Y^+$}
\rput[b](.65,.5){$e$}
\rput[l](.95,.6){$f$}
\end{pspicture}  +
\begin{pspicture}[0.4](-.9,-.3)(1.3,1.2) 
\pscurve{-*}(0,0.05)(.1,.2)(-.2,.45)
\pscurve[border=1pt]{-*}(0,0.05)(-.1,.2)(.2,.45)
\psline{-}(.2,.45)(.35,.75)
\psline{*-*}(0,0.05)(.8,.05) 
\psline(.2,.45)(-.2,.45)
\psline{-}(-.2,.45)(-.35,.75)
\psline{-}(.8,0.05)(.95,.35)
\psline{-}(.8,0.05)(.65,.35)
\rput[br](-.4,.75){$Y$}
\rput[bl](.4,.75){$Y^+$}
\rput[b](.65,.5){$e$}
\rput[l](.95,.6){$f$}
\end{pspicture}\;\rangle \rangle_{W_2}$$
that is zero by the antisymmetry relation in the space of Jacobi diagrams.
Similarly, the contribution of the pairings that pair $c$ and $d$ to
$Y$ and $Y^+$ vanishes.

The contributions of the pairings that pair $d$ and $e$ to
$X$ and $X^+$ is
$$2lk(d,X)lk(e,X)lk(c,Y)lk(f,Y) W_2\left(\tatata + \begin{pspicture}[.2](-.2,-.1)(.8,.6)
\psline(.1,0)(.5,.4)
\psline[border=1pt](.5,0)(.1,.4)
\psline{*-*}(.1,0)(.5,0)
\psline{*-*}(.1,.4)(.5,.4)
\pscurve(.1,0)(0,.2)(.1,.4)
\pscurve(.5,0)(.6,.2)(.5,.4)
\end{pspicture}\right)$$
where $$W_2\left(\begin{pspicture}[.2](-.2,-.1)(.8,.6)
\psline(.1,0)(.5,.4)
\psline[border=1pt](.5,0)(.1,.4)
\psline{*-*}(.1,0)(.5,0)
\psline{*-*}(.1,.4)(.5,.4)
\pscurve(.1,0)(0,.2)(.1,.4)
\pscurve(.5,0)(.6,.2)(.5,.4)
\end{pspicture}\right)=W_2\left(\tetra\right)=1.$$

Therefore, the contribution to $\langle \langle\; I(c,d,e,f)\; I(\Sigma)  \;\rangle \rangle_{W_2}$ of the pairings that are bijections from 
the set of univalent vertices of $I(\Sigma)$ to the set of univalent vertices of $I(c,d,e,f)$
is $$\frac{3}{4} \langle \langle I(c,d,e,f) \begin{pspicture}[0.2](-.6,-.1)(1.2,.6)  
\psline{-}(.35,0.1)(-.05,.1)
\rput[rb](-.1,0){$X$}
\rput[lb](.4,0){$X^+$}
\end{pspicture}\begin{pspicture}[0.2](-.6,-.1)(1.2,.6)  
\psline{-}(.35,0.1)(-.05,.1)
\rput[rb](-.1,0){$Y$}
\rput[lb](.4,0){$Y^+$}
\end{pspicture}\rangle \rangle_{W_1}$$
Therefore, according to Proposition~\ref{propcasknotvar},
$$C_B=\frac{q}{4p} \frac{3}{4} \langle \langle I(F_J) \begin{pspicture}[0.2](-.6,-.1)(1.2,.6)  
\psline{-}(.35,0.1)(-.05,.1)
\rput[rb](-.1,0){$X$}
\rput[lb](.4,0){$X^+$}
\end{pspicture}\begin{pspicture}[0.2](-.6,-.1)(1.2,.6)  
\psline{-}(.35,0.1)(-.05,.1)
\rput[rb](-.1,0){$Y$}
\rput[lb](.4,0){$Y^+$} \end{pspicture} \rangle \rangle_{W_1} = 3\frac{q}{2p} \lambda^{\prime}(J,X,Y)$$
$$= \frac{3}{2}\left(\lambda^{\prime}((X,Y) \subset M(J;q/p)) - \lambda^{\prime}((X,Y) \subset M))\right).$$
\eop

\section{More about surgeries on general knots in rational homology spheres}
\setcounter{equation}{0}

Theorem~\ref{thmpol} describes the polynomial behaviour of $Z_n$ under surgeries
on null-homologous knots. It can easily be generalized to the case of non null-homologous knots $K$ a primitive satellite $\ell$ of which bounds a Seifert surface. Let $m_K$ be the meridian of such a knot $K$ such that
$\langle m_K ,\ell\rangle_{\partial N(K)}=O_K$.

A surgery curve $\mu$ on $\partial N(K)$ is determined by its coordinates $(p_K,q_K)$ in the symplectic basis $(m_K, \frac{1}{O_K} \ell)$ of $H_1(\partial N(K);\QQ)$
where $p_K= \frac{1}{O_K} \langle \mu ,\ell \rangle$ is the linking number
of $K$ and $\mu$, and $q_K=\langle m_K ,\mu \rangle$. The associate surgery coefficient is $p_K/q_K$.

\begin{theorem}
\label{thmpolprim}
Let $n \in \NN$.
Let $K$ be a knot of order $O_K$ in a rational homology sphere $M$ such that a primitive satellite $\ell$ of $K$ bounds a Seifert surface $F$.
Let $F^1,\dots, F^n$ be parallel copies of $F$. Let $p_K/q_K \in \QQ$ be a surgery coefficient for $K$.
Then
$$Z_n(M(K;\frac{p_K}{q_K}))-Z_n(M)=\sum_{i=0}^n Y_{n,q_K/(p_KO_K^2)}^{(i)}(K \subset M) (\frac{q_K}{p_K})^i$$
where $$Y_{n,q_K/(p_KO_K^2)}^{(n)}(K)= \frac{1}{n!2^n O_K^{2n}} \langle \langle \bigsqcup_{i \in \{1,\dots,n\}}  I(F^i) \;\;\rangle \rangle\;$$
$Y_{n,q_K/(p_KO_K^2)}^{(i)}$ only depends on $q_K/(p_KO_K^2)$ mod $\ZZ$,
and, if $n\geq 2$, $p_c(Y_{n,q_K/(p_KO_K^2)}^{(n-1)})= Y_{n,q_K/(p_KO_K^2)}^{(n-1)c}$ does not depend on $p_K$ and $q_K$.
Furthermore, if $m$ is a primitive satellite of $K$ such that $\langle m ,\ell\rangle_{\partial N(K)}=1$, and if $\hat{K} \subset \hat{M}$ is the knot with the same complement as $K$ whose meridian is $m$, then, if $n \geq 2$,
$$Y_{n}^{(n-1)c}(K \subset M)=\frac{1}{O_K^{2n-2}}Y_{n}^{(n-1)c}(\hat{K} \subset \hat{M})
+n \langle m, m_K \rangle  O_K Y_{n}^{(n)c}(K \subset M).$$
\end{theorem}
\bp Let $\mu=p_Km_K +(q_K/O_K)\ell$ be a surgery curve on $\partial N(K)$. Let $$\left(p=\langle \mu ,\ell \rangle = O_K p_K,q=\langle m ,\mu \rangle=p_K \langle m, m_K \rangle + q_K/O_K\right)$$ be the coordinates of $\mu$ in the symplectic basis $(m, \ell)$ of $H_1(N(K);\ZZ)$.
Note that changing $m$ to another curve such that $\langle m ,\ell\rangle_{\partial N(K)}=1$ leaves $p$ invariant and does not change the class of $\frac{q}{p}$ in $\QQ/\ZZ$.
When the other data are fixed, the mod $\ZZ$ congruence class of 
$$\frac{q}{p}=\frac{q_K}{p_K O_K^2} + \frac{\langle m, m_K \rangle}{O_K}$$
depends on the class of $\frac{q_K}{p_K O_K^2}$ in $\QQ/\ZZ$.
From the formula of Theorem~\ref{thmpol}
$$Z_n(\hat{M}(\hat{K};\frac{p}{q+rp}))-Z_n(\hat{M})=\sum_{i=0}^n Y_{n,q/p}^{(i)}(\hat{K} \subset \hat{M}) (r+\frac{q}{p})^i,$$
we deduce
$$Z_n(M(K;\frac{p_K}{q_K +r O_K^2 p_K}))-Z_n(M)$$
$$=\sum_{i=0}^n Y_{n,q/p}^{(i)}(\hat{K} \subset \hat{M}) (r+ \frac{q_K}{p_K O_K^2} + \frac{\langle m, m_K \rangle}{O_K})^i+Z_n(\hat{M}) -Z_n(M) $$
$$=\sum_{i=0}^n Y_{n,q_K/(p_K O_K^2)}^{(i)}(K \subset M) (rO_K^2+ \frac{q_K}{p_K })^i$$
where
$$Y_{n,q_K/(p_K O_K^2)}^{(n)}(K \subset M)=\frac{1}{n!2^nO_K^{2n}}\langle \langle \bigsqcup_{i \in \{1,\dots,n\}}  I(F^i) \;\;\rangle \rangle\;$$
and, if $n\geq 2$,
$$Y_{n,q_K/(p_K O_K^2)}^{(n-1)}(K \subset M)=\frac{1}{O_K^{2n-2}}Y_{n,q/p}^{(n-1)}(\hat{K} \subset \hat{M})
+n \langle m, m_K \rangle  O_K Y_{n}^{(n)}(K \subset M).$$
\eop

\begin{remarks}
A knot $K$ of order $O_K$ in a rational homology sphere has a primitive satellite that is null-homologous in its exterior if and only if the self-linking number of $K$ reads $d/O_K$ (mod $\ZZ$) where $d$ is coprime with $O_K$.

Like in the proof of Theorem~\ref{thmpolprim}, the case of knots without null-homologous primitive satellites
can be reduced to the case of knots of order $O_K >1$ with self-linking number
$0$.
This latter case is still unclear to me
(except for the degree 1 case that can be treated with the methods of the article).

Relationships between surgery formulae for various $q/p$ can be found
using some equivalences of surgeries. See~\cite{go5}.
\end{remarks}

\section{Questions}
\setcounter{equation}{0}

The statements of Theorems~\ref{thmfas} and \ref{thmfasmu} make sense for
rationally algebraically split links. Do they hold true in this case?

 How do the properties of surgery formulae generalize for surgeries on non null-homologous knots?
 
What is the graded space associated to the filtration of the rational vector space generated by rational homology spheres, defined using Lagrangian-preserving surgeries?

 The degree $n$ parts of the LMO invariant and the Kontsevich-Kuperberg-Thurston invariant coincide on the intersection
 of $\CF_n$ with the vector space generated by homology spheres.
 The configuration space invariant for knots in $S^3$ is obtained from the Kontsevich integral by an isomorphism that inserts a (possibly trivial) specific two-leg box $\beta$ on each chord of a chord diagram. See \cite{les5} for a more specific statement.
Do the LMO invariant and the Kontsevich-Kuperberg-Thurston invariant actually coincide?
 Is the Kontsevich-Kuperberg-Thurston invariant obtained from the LMO invariant by inserting the two-leg box $\beta$, $k$ (or $2k$ or $3k$) times on each degree $k$ component of a Jacobi diagram?

\end{document}